\documentclass[11pt]{article}

\usepackage[margin=2.5cm]{geometry}
\usepackage[english]{babel}
\usepackage{amsmath,amssymb,amsthm}
\usepackage{subfig}
\usepackage{picins}

\title{Cardinality of Rauzy classes}
\author{Vincent Delecroix}

\theoremstyle{plain}
\newtheorem{theorem}{Theorem}[section]
\newtheorem{lemma}[theorem]{Lemma}
\newtheorem{proposition}[theorem]{Proposition}
\newtheorem{corollary}[theorem]{Corollary}

\theoremstyle{definition}
\newtheorem{definition}[theorem]{Definition}

\theoremstyle{remark}

\numberwithin{equation}{section}
\numberwithin{figure}{section}

\def\FF{\mathbb{F}}
\def\NN{\mathbb{N}}
\def\ZZ{\mathbb{Z}}
\def\QQ{\mathbb{Q}}
\def\RR{\mathbb{R}}
\def\CC{\mathbb{C}}

\def\PP{\mathbb{P}}

\def\AAA{\mathcal{A}}

\def\ker{\operatorname{ker}}
\def\Arf{\operatorname{Arf}}

\newcommand\outgos[1]{\overline{#1}}
\newcommand\incoms[1]{\underline{#1}}

\usepackage{graphics}
\def\picinput#1{\includegraphics{#1.pdf}}

\begin{document}

\maketitle

\abstract{We prove an explicit formula for cardinalities of Rauzy classes of permutations introduced as part of a renormalization algorithm on interval exchange transformations. Our proof uses a geometric interpretation of permutations and Rauzy diagrams in terms of translation surfaces and moduli spaces.}

\tableofcontents

\section{Introduction}
\label{section:introduction}
Let $\lambda = (\lambda_1, \lambda_2, \ldots, \lambda_n) \in \RR_+^n$ be a vector with positive coordinates and $\pi \in S_n$ be a permutation. For $i=2,\dots,n+1$, we define $x_i = \sum_{j=1}^{i-1} \lambda_j$ and $y_i = \sum_{j=1}^{i-1} {\lambda_\pi^{-1}(j)}$ and note $x_1=y_1=0$ and $a=x_{n+1}=y_{n+1}$. The interval exchange transformation $T = T_{\lambda,\pi}$ with data $(\lambda,\pi)$ is the map defined on $[0,\ a)$ into itself by
\[
T(x) = x - x_i + y_{\pi(i)} \quad \text{on $[x_i,\, x_{i+1})$}.
\]
In other words, on the subinterval $[x_i,\, x_{i+1})$, the map $T$ acts as a translation by $y_{\pi(i)}-x_i$. An interval exchange transformation is bijective and right continuous. The map $T$ is an examples of measurable dynamical system as it preserves the Lebesgue measure on $[0,\ a)$.
\begin{figure}[!ht]
\centering
\picinput{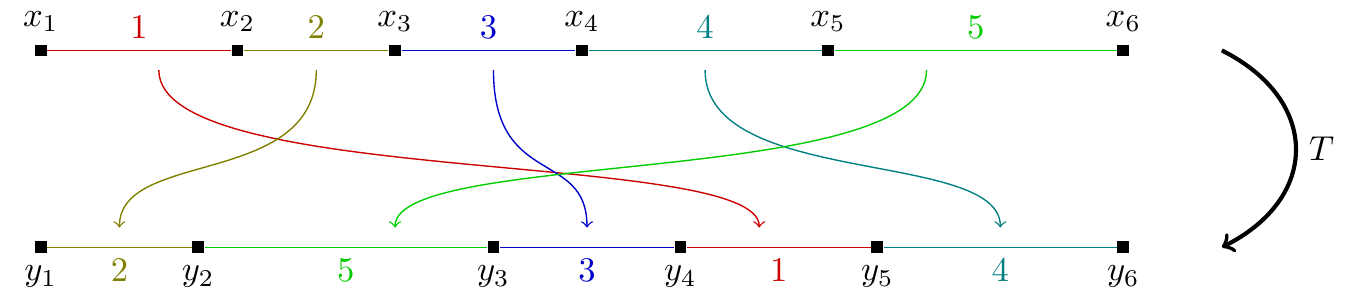}
\caption{An interval exchange transformation with permutation $\pi=\binom{1\, 2\, 3\, 4\, 5}{2\, 5\, 3\, 1\, 4}$.}
\end{figure}

If $\pi(\{1,2,\ldots,k\}) = \{1,2,\ldots,k\}$ for $k$ such that $1 \leq k < n$, then the two subintervals $[0,\ x_{k+1})$ and $[x_{k+1},\ a)$ are invariant under $T_{\lambda,\pi}$. We are interested in permutations that do not allow such a splitting.
\begin{definition} \label{def:irreducible_permutation}
A permutation $\pi \in S_n$ is \textit{irreducible} (or \textit{indecomposable}) if there is no $k$, $1 \leq k < n$, such that $\pi(\{1,2,\ldots,k\}) = \{1,2,\ldots,k\}$.
\end{definition}
We denote by $S_n^o$ the set of irreducible permutations in $S_n$. It was proved by M.~Keane \cite{Keane75} that if $\pi \in S^o_n$ then for Lebesgue almost all $\lambda \in \RR_+^n$ the interval exchange $T_{\lambda,\pi}$ is minimal. Later H.~Masur \cite{Masur82} and W.A.~Veech \cite{Veech82}, independently, proved that for Lebesgue almost all $\lambda \in \RR_+^n$ the interval exchange $T_{\lambda,\pi}$ is uniquely ergodic.

In order to study the dynamics of interval exchange transformations, \cite{Rauzy79} defines an induction procedure (named \textit{Rauzy induction}) on the space of interval exchange transformations. In other words, a map $R: S^o_n \times \RR_+^n \rightarrow S^o_n \times \RR_+^n$. There are two cases of induction depending whether $x_n < y_n$ (\textit{top induction}) or $x_n > y_n$ (\textit{bottom induction}). The induction is not defined if $x_n = y_n$. Let $T_{\lambda,\pi}$ be an interval exchange transformation with $x_n \not= y_n$ and $T_{\lambda',\pi'} = R \left( T_{\lambda,\pi} \right)$ the one obtained by Rauzy induction. The permutation $\pi'$ only depends on the type of the induction. Hence, there are two combinatorial operations $R_t:S_n^o \rightarrow S_n^o$ (top induction) and $R_b:S_n^o \rightarrow S_n^o$ (bottom induction) which corresponds to the operation on the permutation $\pi$ associated to the Rauzy induction.  The equivalence classes induced by the action of $R_t$ and $R_b$ on $S_n^o$ are called \textit{Rauzy classes}.

As far as we know, only the Rauzy class $\mathcal{R}^{sym}$ of the symmetric permutation $\pi^{sym} \in S_n^o$ defined by $\pi^{sym}(k) = n-k+1$ for $k=1,\ldots,n$ had been described in \cite{Rauzy79}. In particular Rauzy proved that its cardinality is $|\mathcal{R}^{sym}| = 2^{n-1}-1$. Motivated by the study of \cite{Rauzy79}, the aim of this article is to study the combinatorics of Rauzy classes in $S^o_n$ and establish a formula for their cardinalities.

\subsection*{Aknowledgment}
I wish to thank Corentin Boissy, Erwan Lanneau and Arnaldo Nogueira for their patient lectures and comments and Samuel Leli\`evre for his help on drawing pictures with the PGF Ti\textit{k}Z library for \LaTeX. All formulas in the article have been tested with the help of the mathematical software Sage \cite{Sage}.

\section{Main results}
We recall elements from Teichm\"uller theory which yield to a classification of Rauzy diagrams. Let $I = [0,a)$ be an interval and $T$ a map from $I$ into itself. Let $f:[0,a) \rightarrow \RR_+$ and $X$ be the quotient of $\{(x,y) \in [0,a) \times \RR_+;\ y \leq f(x) \}$ by the relation $(x,f(x)) \sim (T(x),0)$. The space $X$ together with the flow $\phi_t$ in the vertical direction is called a \emph{suspension} and $f$ the \emph{roof function}. The flow $\phi_t$ has the property that the first return map on the interval $I \times \{0\} \subset X$ is exactly the map $T$. W.~Veech \cite{Veech82} considered roof functions which are constant on each subinterval of an interval exchange transformation. The suspension $X$ obtained by this procedure is a translation surface and the flow $\phi_t$ corresponds to the geodesic flow on $X$ in the vertical direction. Translation surfaces are part of Teichm\"uller theory and will play an important role in the construction of our counting formulas.

A translation surface $S$ has a flat metric, except at a finite number of points where there are conical singularities whose angles are integer multiples of $2\pi$. If $S$ is a suspension of an interval exchange transformation $T$, the conical singularities of $S$ come from the singularities of $T$. Let $n_1 2\pi, n_2 2\pi, \ldots, n_k 2\pi$ be the list of angles of the conical singularities of $S$. We call the integer partition $p=(n_1,n_2,\ldots,n_k)$ the {\em profile} of $S$. The genus $g$ of $S$ is related to $p$ by $2g-2 = \sum_{i=1}^k (n_i - 1) = s(p) - l(p)$ where $s(p) = n_1 + \ldots + n_k$ is the \emph{sum} of $p$ and $l(p) = k$ its \emph{length}. We emphasize that the suspension associated to an interval exchange transformation is not unique but all of them have the same profile. Furthermore, suspensions obtained from permutations in the same Rauzy class have the same profile. Let $\Omega \mathcal{M}_g$ be the moduli space of translation surfaces of genus $g$ and $p =(n_1,n_2,\ldots,n_k)$ an integer partition such that $2g-2=\sum_{i=1}^k (n_i -1)$. The {\em stratum with profile $p$} denoted $\Omega \mathcal{M}_g(p)$ is the subset of $\Omega \mathcal{M}_g$ made out of surfaces whose profiles are $p$. On $\Omega \mathcal{M}_g$ acts the Teichm\"uller flow which preserves strata and for which the Rauzy-Veech induction on suspensions is viewed as a first return map. There is a bijection between extended Rauzy classes and connected components of strata $\Omega \mathcal{M}(p)$ \cite{Veech82} where an {\em extended Rauzy class} is an equivalence class of irreducible permutations under the action of $R_t$, $R_b$ and $s$, where $s$ is the operation which acts on $\pi \in S_n$ by $s(\pi)(k) = \pi^{-1}(n-k-1)$. C.~Boissy \cite{Boissy09} proved a bijection between Rauzy classes and connected components of strata with a choice of a part of the profile $p$. This choice corresponds to the marking of suspension induced by the left endpoint of the interval exchange transformation which is not affected during the Rauzy inductions. The combinatorial question of classifying (extended) Rauzy classes is hence translated into the geometric one of classifying connected components. M.~Kontsevich and A.~Zorich \cite{KontsevichZorich03} classified connected components of strata in terms of geometrical invariants: the \textit{spin parity} (an element of $\{0,1\}$) and the \textit{hyperellipicity}. A spin parity occurs when the profile $p$ has only odd parts and give rise to at least two distinct connected components. The term hyperellipticity stands for a serie of connected components that appear for the profiles $(2g-1,1^k)$ and $(g,g,1^k)$. This yields to a classification of (extended) Rauzy classes.

Our approach to count permutations in Rauzy classes relies in the above geometric interpretation of Rauzy classes. Let $\pi \in S_n^o$ be an irreducible permutation and $p_\pi$ the profile of a suspension of $\pi$. The profile does not reflect the structure of an embedded segment in a surface and we refine the notion. We say that $\pi$ has {\em marking} $m|a$ if the extremities of the interval corresponds to the same singularity $P$ in the suspension which has a conical angle $m 2 \pi$ and $a \in \{0,\ldots,m-1\}$ is such that $(2a+1)\pi$ is the angle between the left part and the right part of the interval measured from $P$. It has {\em marking} $m_l \odot m_r$ if the two extremities of the interval correspond to two different singularities of angles $m_l 2 \pi$ on the left and $m_r 2 \pi$ on the right. The data which consists of the profile and the marking is called the {\em marked profile} of the permutation $\pi$. We denote by $(m|a,p')$ (resp. $(m_l \odot m_r, p')$) a profile $p' \uplus (m)$ (resp. $p' \uplus (m_l,m_r)$) with marking $m|a$ (resp. $m_l \odot m_r$). Here $\uplus$ stands for the disjoint union of partitions considered as multisets. Our main theorem (see below) is a recurrence formula for the number of irreducible permutations with given marked profile.

We first consider standard permutations introduced in \cite{Rauzy79}.
\begin{definition} \label{def:standard_permutation}
A permutation $\pi \in S_n$ is {\em standard} if $\pi(1) = n$ and $\pi(n) = 1$.
\end{definition}
A standard permutation is in particular irreducible. Those permutations were used for dynamical purpose in \cite{NogueiraRudolph97} and \cite{AvilaForni07} in order to prove the weak mixing property of interval exchange transformations and in \cite{KontsevichZorich03},\cite{Zorich08} and \cite{Lanneau08} in the study of connected components of strata. In terms of moduli space of translation surfaces, a standard permutation corresponds to a so called {\em Strebel differential}.

Let $\mathbf{p}$ be a marked profile whose profile is $p$. We denote by $\gamma^{std}(\mathbf{p})$ the number of standard permutations with marked profile $\mathbf{p}$. Moreover, if $p$ has only odd terms, we define $\delta^{std}(\mathbf{p}) = \gamma_1^{std}(\mathbf{p}) - \gamma_0^{std}(\mathbf{p})$ where $\gamma_s^{std}(\mathbf{p})$ denotes the number of standard permutations with marked profile profile $\mathbf{p}$ and spin parity $s$. We prove explicit formulas for $\gamma^{std}(\mathbf{p})$ and $\delta^{std}(\mathbf{p})$. The formulas involve the numbers $z_p$, $c(p)$ and $d(p)$ which are defined next but we first introducte notations for partitions. Let $p = (n_1,n_2,\ldots,n_k)$ be an integer partition considered as a multiset (each part has a multiplicity equals to its number of occurences in the partition). We recall that the disjoint union is denoted $\uplus$. We have $s(p_1 \uplus p_2) = s(p_1) + s(p_2)$ and $l(p_1 \uplus p_2) = l(p_1) + l(p_2)$. If $q \subset p$ is a subpartition of $p$ we denote by $p \backslash q$ the unique partition $r$ such that $p = q \uplus r$.

We recall that the conjugacy classes of $S_n$ are in bijection with integer partitions of $n$. We denote by $z_p$ the cardinality of the conjugacy class associated to $p$. If $e_i$ is the number of occurences of $i$ in $p$, then
\[
z_p = \prod_{i=1}^n i^{e_i} e_i!\ i^{e_i}.
\]
If $p$ satisfies $s(p) + l(p) \equiv 0 \mod 2$ we define (the formula is due to G.~Boccara \cite{Boccara80})
\[
c(p) = \frac{2 (n-1)!}{n+1} \left( \sum_{q \subset (n_2,n_3,\ldots,n_k)} (-1)^{s(q)+l(q)} \binom{n}{s(q)}^{-1} \right)
\]
where the summation is above each subpartition of $(n_2,n_3,\ldots,n_k)$ with multiplicity in the sense that $(1,3)$ occurs twice in $(1,1,3)$. Moreover, if the partition $p$ has only odd parts we define $d(p) = (n-1)! / 2^g$ where $g=(s(p)-l(p))/2$.

Our proofs are based on surgeries of partitions which are used to obtain recurrence (with a geometric counterpart as in \cite{KontsevichZorich03} and \cite{EskinMasurZorich03}). If $m$ is a part of $p$ and $a \in \{0,\ldots,m-1\}$ we denote by $p_{m|a}$ the partition obtained from $p$ by removing $m$ and inserting the two parts $a$ and $m-a-1$ (if $a$ is $0$ or $m-1$ we replace $m$ by $m-1$). If $m_l$ and $m_r$ are two distinct parts of $p$ we denote by $p_{m_l \odot m_r}$ the partition obtained from $p$ by removing the parts $m_l$ and $m_r$ and inserting $m_l + m_r - 1$. We have $s(p_{m|a}) = s(p)-1$ and $s(p_{m_l \odot m_r}) = s(p)-1$ (notations $p_{m|a}$ and $p_{m_l \odot m_r}$ comes from \cite{Boccara80}).

\begin{theorem} \label{thm0:standard_perms}
Let $p$ be an integer partition such that $s(p)+l(p) \equiv 0 \mod 2$. Let $m$ be a part of $p$, $a \in \{1,\ldots,m-2\}$. Set $p' = p \backslash (m)$. Then, we have
\begin{equation} \label{eq:gamma_std}
\gamma^{std}(m|a,p') = \frac{c(p_{m|a})}{z_{p'}}
\qquad \text{and} \qquad
\delta^{std}(m|a,p') =
\left\{ \begin{array}{ll}
0 & \text{if $a \equiv 0 \mod 2$} \\
 \frac{d(p_{m|a})}{z_{p'}} & \text{otherwise}
\end{array} \right. \ .
\end{equation}

Assume that $p$ has only odd terms. Let $m_l$ and $m_r$ be two distinct parts of $p$. Set $p' = p \backslash (m_l,m_r)$. Then, we have
\begin{equation} \label{eq:delta_std}
\gamma^{std}(m_l \odot m_r,p') = \frac{d(p_{m_l \odot m_r})}{z_{p'}}
\qquad \text{and} \qquad
\delta^{std}(m_l \odot m_r,p') = \frac{d(p_{m_l \odot m_r})}{z_{p'}}.
\end{equation}
\end{theorem}
The numbers $c(p)$ and $d(p)$ can be interpreted as counting of labeled permutations and $z_{p'}$ as the cardinality of the group which exchanges the labels.

Let $\mathbf{p}$ be a marked profile. We define $\gamma^{irr}(\mathbf{p})$ (resp. $\delta^{irr}(\mathbf{p}) = \gamma_1^{irr}(\mathbf{p}) - \gamma_0^{irr}(\mathbf{p})$) the number of irreducible permutations with given marked profile (resp. the difference between the numbers of irreducible permutations with odd and even spin parity). The below theorem gives recursive formulas for the numbers $\gamma^{irr}$ and $\delta^{irr}$ which involve the numbers $\gamma^{std}$ and $\delta^{std}$.
\begin{theorem} \label{thm0:cardinality}
Let $p=(m_1, m_2, \ldots)$  be an integer partition such that $s(p) + l(p) \equiv 0 \mod 2$. Let $m \in p$ and $a \in \{0, \ldots, m-1\}$ then
\[
\gamma^{irr} \left(m|a,p'\right) = 
\gamma^{std} \left(m+2|m-a,p'\right) -
\sum_{\substack{%
    {\scriptscriptstyle p'_1 \uplus p'_2 = p'}\\
    {\scriptscriptstyle m_1 + m_2 = m-1}\\
    {\scriptscriptstyle a_1 + a_2=a-1}}}
  \gamma^{irr} \left(m_1|a_1,p'_1 \right) \ \gamma^{std} \left(m_2+2|m_2-a_2,p'_2 \right).
\]
Let $(m_l,m_r) \subset p$ then
\begin{align*}
\gamma^{irr} (m_l \odot m_r, p')=& \gamma^{std}((m_l+1) \odot (m_r+1),\, p') \\
&-\sum_{p'_1 \uplus p'_2 = p'} \quad
\sum_{\substack{%
      {\scriptscriptstyle k_1+k_2=m_l-1}\\
      {\scriptscriptstyle 0 \leq a < k_1}}}
{\textstyle \gamma^{irr}(k_1|a,\,p'_1) \ \gamma^{std}((k_2+1) \odot (m_r+1),\,p'_2)} \\
&- \sum_{p'_1 \uplus p'_2 = p'} \quad
\sum_{\substack{%
      {\scriptscriptstyle k_1+k_2=m_r-1}\\
      {\scriptscriptstyle 1 \leq a < k_1-1}}}
{\textstyle \gamma^{irr}(m_l \odot k_1,\,p'_1) \ \gamma^{std}(k_2+2|a,\,p'_2) } \\
&-\sum_{p'' \uplus (m) = p'} \quad 
\sum_{\substack{%
    {\scriptscriptstyle p_1 \uplus p_2 = p''}\\
    {\scriptscriptstyle k_1+k_2=m-1}}}
\gamma^{irr}(m_l \odot k_1,\,p_1) \ \gamma^{std}((k_2+1) \odot (m_r+1),\,p_2)
\end{align*}
Moreover, if $p$ has only odd parts we have
\[
\delta^{irr} \left(m|a,p'\right) = 
(-1)^a\ \delta^{std} \left(m+2|m-a,p'\right) -
\sum_{\substack{%
    {\scriptscriptstyle p'_1 \uplus p'_2 = p'}\\
    {\scriptscriptstyle m_1 + m_2 = m-1}\\
    {\scriptscriptstyle a_1 + a_2=a-1}}}
  (-1)^{a_2} \delta^{irr} \left(m_1|a_1,p'_1 \right) \ \delta^{std} \left(m_2+2|m_2-a_2,p'_2 \right).
\]
And if $(m_1,m_2) \subset p$ then
\begin{align*}
\delta^{irr} (m_l \odot m_r, p')=& \delta^{std}(m_l+m_r+1,\, p') \\
&+\sum_{p'_1 \uplus p'_2 = p'} \quad
\sum_{\substack{%
      {\scriptscriptstyle k_1+k_2=m_l-1}\\
      {\scriptscriptstyle 0 \leq a < k_1}}}
{\textstyle \delta^{irr}(k_1|a,\,p'_1) \ \delta^{std}(k_2+m_r+1,\,p'_2)} \\
&+ \sum_{p'_1 \uplus p'_2 = p'} \quad
\sum_{\substack{%
      {\scriptscriptstyle k_1+k_2=m_r-1}\\
      {\scriptscriptstyle 1 \leq a < k_1-1}}}
{\textstyle \delta^{irr}(m_l \odot k_1,\,p'_1) \ \delta^{std}(k_2+2|a,\,p'_2) } \\
&+\sum_{p'' \uplus (m) = p'} \quad 
\sum_{\substack{%
    {\scriptscriptstyle p_1 \uplus p_2 = p''}\\
    {\scriptscriptstyle k_1+k_2=m-1}}}
\delta^{irr}(m_1 \odot k_1,\,p_1) \  \delta^{std}(k_2+m_r+1,\,p_2)
\end{align*}
where $\delta^{std}(m,p') = \sum_{a =1}^{m-2} \delta^{std}(m|a,p')$.
\end{theorem}

Theorems~\ref{thm0:cardinality} and~\ref{thm0:standard_perms} do not treat the case of Rauzy classes associated to hyperelliptic components $\Omega \mathcal{M}^{hyp}_g (2g-1,1^k)$ and $\Omega \mathcal{M}^{hyp}_g(g,g,1^k)$ where $(1^k)$ denotes the partition that contains $k$ times the part $1$. The component $\Omega \mathcal{M}^{hyp}_g(2g-1)$ (resp. $\Omega \mathcal{M}^{hyp}_g(g,g)$) corresponds to the extended Rauzy class of the symmetric permutation of degree $2g$ (resp. $2g+1$). We know since \cite{Rauzy79} that the cardinality of the extended Rauzy class of the symmetric permutation of degree $n$ is $2^{n-1}-1$. To obtain the cardinality of each hyperelliptic class, we establish a general formula that relates the cardinality of an extended Rauzy class $\mathcal{R}$ associated to a profile $p$ to the one of $\mathcal{R}_0$ obtained from $\mathcal{R}$ by adding $k$ marked points. The extended Rauzy class $\mathcal{R}_0$ has profile $p \uplus (1^k)$.
\begin{theorem} \label{thm0:adding_marked_points}
Let $k$, $\mathcal{R}$ and $\mathcal{R}_0$ be as above. Let $r$ the number of standard permutations in $\mathcal{R}$ then
\[
|\mathcal{R}_0| = \binom{n+k+1}{k}\ |\mathcal{R}| + \binom{n+k}{k-1}\ n\,r.
\]
\end{theorem}
As a particular case of the above theorem, we obtain an explicit formula for the cardinalities of Rauzy classes associated to hyperelliptic components.
\begin{corollary}
Let $\mathcal{R} \subset S^o_{2g+k}$ (resp. $\mathcal{R} \subset S^o_{2g+k+1}$) be the extended Rauzy class associated to $\Omega \mathcal{M}^{hyp}(2g-1,1^k)$ (resp. $\Omega \mathcal{M}^{hyp}(g,g,1^k)$) and $n=2g$ (resp. $n=2g+1$), then
\[
|\mathcal{R}| = \binom{n+k+1}{k}\ (2^{n-1}-1) + \binom{n+k}{k-1}\ n.
\]
\end{corollary}

\vspace{0.6cm}

The paper is organized as follows. In Section~\ref{section:iet_and_translation_surfaces}, we review the definitions of Rauzy classes and extended classes. We describe the Rauzy classes of the symmetric permutation $\pi^{sym} \in S_n^o$ defined by $\pi(k) = n-k+1$ (Section~\ref{subsubsection:symmetric_permutation}) and the permutation of rotation class $\pi^{rot} \in S^o_n$ defined by $\pi^{rot}(1) = n$, $\pi^{rot}(n) = 1$ and $\pi^{rot}(k) = k$ for $k=2,\ldots,n-1$ (Section~\ref{subsubsection:Rauzy_diagrams_of_rotations}). We recall the classification of Rauzy classes and extended Rauzy classes in terms of connected components of strata of the moduli space of Abelian differential. In particular, we obtain a formula for cardinalities of Rauzy classes in terms of the numbers $\gamma^{irr}(\textbf{p})$ and $\delta^{irr}(\textbf{p})$. In section~\ref{section:counting_standard_permutations}, we study standard permutations in order to proove Theorem~\ref{thm0:standard_perms}. In section~\ref{section:from_standards_to_all} we see how standard permutations can be used to describe the set of all permutations and prove Theorems~\ref{thm0:cardinality} and \ref{thm0:adding_marked_points}.

\subsection*{Proofs overview}
Now, we explain our strategy to compute cardinalities of Rauzy diagrams. 

First we formulate a definition of Rauzy classes in terms of invariants of permutations in Section~\ref{section:iet_and_translation_surfaces} (see in particular Theorem~\ref{thm:classification_of_Rauzy_classes}). This reformulation follows from the work of \cite{Veech82}, \cite{Boissy09} and the classification of connected components of strata of Abelian differentials in \cite{KontsevichZorich03}. Using this geometric definition, we are able to express cardinalities of Rauzy classes in terms of the numbers $\gamma^{irr}(\mathbf{p})$ and $\delta^{irr}(\mathbf{p})$ which counts irreducible permutations with given profile $p$ (see Corollary~\ref{corollary:formula_for_the_cardinality_of_a_Rauzy_class}).

The computation of the numbers $\gamma^{irr}(\mathbf{p})$ and $\delta^{irr}(\mathbf{p})$ is done in two steps. Both steps use geometrical surgeries used in the classification of connected components of strata \cite{KontsevichZorich03} and \cite{Lanneau08}. The first step consists in studying standard permutations. We consider the numbers $c(p)$ and $d(p)$ of labeled permutations and get a recurrence in terms of partitions of $n-1$ for both of them (Theorems~\ref{thm:std_perm_strata_recurrence} and \ref{thm:std_perm_spin_recurrence}). We then prove that the recurence can be solved into explicit formulas (Theorems~\ref{thm:Boccara_formula_for_c} and \ref{thm:formula_for_d}). These explicit formula corresponds to the formula given in the above introduction. The link between standard permutations and the number of labeled standard permutations as in Theorem~\ref{thm0:standard_perms} is proved in Corollary~\ref{corollary:c_and_gamma} and \ref{corollary:d_and_delta}.

The second step consists in proving Theorem~\ref{thm0:cardinality} which express the numbers $\gamma^{irr}(\mathbf{p})$ (resp. $\delta^{irr}(\textbf{p})$) in terms of $\gamma^{std}(\mathbf{p})$ (resp. $\delta^{std}(\mathbf{p})$). We use a simple construction: to a standard permutation $\pi \in S_n$ we associate the permutation $\tilde{\pi} \in S_{n-2}$ obtained by ``removing its ends''. Formally $\tilde{\pi}(k) = \pi(k+1)-1$ for $k=1,\ldots,n-2$. The operation $\pi \mapsto \tilde{\pi}$ gives a (trivial) combinatorial bijection between standard permutations in $S_n$ and all permutations in $S_{n-2}$. As the permutations obtained by this operation are not necessarily irreducible we define Rauzy classes of reducible permutations. To any permutation we can associate a profile and a spin invariant (see Sections~\ref{subsubsection:profile} and \ref{subsubsection:spin_invariant}). As each permutation is a unique concatenation of irreducible permutations, we study how are related the invariants of a permutation to the invariants of its irreducible components (this is done in Lemmas~\ref{lemma:profile_concatenation} and \ref{lemma:spin_parity_concatenation}). In geometric terms, a reducible permutation corresponds to an ordered list of surfaces in which each surface is glued to the preceding and the next one at a singularity. The operation $\pi \to \tilde{\pi}$ can be analyzed as a surgery operation and the invariants of $\tilde{\pi}$ depend only on the ones of $\pi$ (Proposition~\ref{prop:profile_of_degeneration} and \ref{prop:spin_of_degeneration}). Theorem~\ref{thm0:cardinality} follows from an inclusion-exclusion counting for irreducible permutations among all permutations.

\section{Permutations, interval exchange transformations and translation surfaces}
\label{section:iet_and_translation_surfaces}
In this section we define the Rauzy induction of interval exchange transformations on the space of parameters $\RR^n_+ \times S^o_n$. We study in particular the two combinatorial operations on irreducible permutations $R_t,R_b: S^o_n \rightarrow S^o_n$ which define Rauzy classes. Next we recall the relation between translation surfaces and interval exchange transformations. Our aim is to give another definition of Rauzy classes and extended Rauzy classes (Definition~\ref{def:rauzy_class_combinatorial}) as well as a classification in terms of invariants of a permutation: the \emph{profile} which is an integer partition, the \emph{hyperellipticity} and the \emph{spin parity} which is an element of $\{0,1\}$ (Theorem~\ref{thm:classification_of_Rauzy_classes}).

We recall that if $p$ be an integer partition then we denote by $\gamma^{irr}(p)$ (resp. $\gamma^{irr}_1(p)$ and $\gamma^{irr}_0(p)$) the number of irreducible permutations with profile $p$ (resp. profile $p$ and spin parity $1$ and $0$). We set $\delta^{irr}(p) = \gamma^{irr}_1(p) - \gamma^{irr}_0(p)$. The cardinality of every Rauzy class, but the ones which are associated to components of strata which contain an hyperelliptic component, depend only on the numbers $\gamma^{irr}(p)$ and $\delta^{irr}(p)$ (see~Corollary~\ref{corollary:formula_for_the_cardinality_of_a_Rauzy_class}).

The next two sections of this paper are devoted to the computations of $\gamma^{irr}(p)$ and $\delta^{irr}(p)$. The explicit formulas for the cardinalities of hyperelliptic Rauzy classes are given in Corollary~\ref{cor:cardinality_of_hyperelliptic_rauzy_diagrams}.

\subsection{Rauzy induction and Rauzy classes}
\subsubsection{Labeled permutations}
\label{subsubsection:labeled_permutations}
We introduce a labeled version of permutations which comes from \cite{MarmiMoussaYoccoz05} and \cite{Bufetov06} inspired from \cite{Kerchoff85} (see also \cite{Boissy10}). Many constructions are easier to formulate with this definition.
\begin{definition}
A \emph{labeled permutation on a finite set $\mathcal{A}$} is a couple of bijections $\pi_t,\pi_b : \mathcal{A} \rightarrow \{1,\ldots,n\}$ where $n$ is the cardinality of $\mathcal{A}$. The elements of $\mathcal{A}$ are called the \emph{labels} of $(\pi_t,\pi_b)$ and $\mathcal{A}$ the \emph{alphabet}.
\end{definition}
In order to distinguish labeled permutations from permutations we will sometimes call them \emph{reduced permutations} instead of permutations. The number $n$ is called the \emph{length} of the permutation. To a labeled permutation we associate a reduced one by the map $(\pi_t,\pi_b) \mapsto \pi_b \circ \pi_t^{-1}$. We also consider the natural section given by $\pi \mapsto (id,\pi)$ for which the alphabet of the labeled permutation $(id,\pi)$ is $\{1,\, 2,\, \ldots, \,n\}$.

A labeled permutation $\pi = (\pi_t,\pi_b)$ is written as a table with two lines
\[
\pi = \left( \begin{array}{cccc}
\pi_t^{-1}(1) & \pi_t^{-1}(2) & \ldots & \pi_t^{-1}(n)  \\
\pi_b^{-1}(1) & \pi_b^{-1}(2) & \ldots & \pi_b^{-1}(n)  
\end{array} \right)
\]
The \emph{top line} (resp. \emph{bottom line}) of $\pi$ is the ordered list of labels $\pi_t^{-1}(i)$ for $i=1,\ldots,n$ (resp. $\pi_b^{-1}(i)$ for $i=1,\ldots,n$). For a reduced permutation $\pi$ we use the section defined above and write
\[
\pi = \left( \begin{array}{cccc}
1 & 2 & \ldots & n \\
\pi^{-1}(1) & \pi^{-1}(2) & \ldots & \pi^{-1}(n)
\end{array} \right)
\quad \text{or simply} \quad
\left( \pi^{-1}(1)\ \pi^{-1}(2)\ \ldots\ \pi^{-1}(n) \right).
\]
The above notation coincides with the notation of $\pi^{-1}$ in group theory. With our notation, the label $i$ is at the position $\pi(i)$ on the bottom line. The difference of notation will not cause any problem as we never use the composition of permutations that arises from interval exchange transformations. The only operation considered here is the concatenation (see Section~\ref{subsection:concatenation_irreducibility}).

The definitions of standard and irreducible permutations extend to labeled permutations.
\begin{definition} \label{def:irreducible_and_standard_labeled_permutation}
We say that $(\pi_t,\pi_b)$ is \emph{irreducible} (resp. \emph{standard}) if $\pi_b \circ \pi_t^{-1} \in S_n$ is irreducible (resp. standard).
\end{definition}
Our aim is to count reduced permutations, however in Section~\ref{section:counting_standard_permutations} we will mainly deal with labeled ones. In \cite{Boissy10}, C.~Boissy analyze the difference between reduced and labeled permutations.

\subsubsection{Rauzy induction and Rauzy classes}
\label{subsubsection:combinatorial_rauzy_classes}
\begin{figure}[!ht]
\centering
\picinput{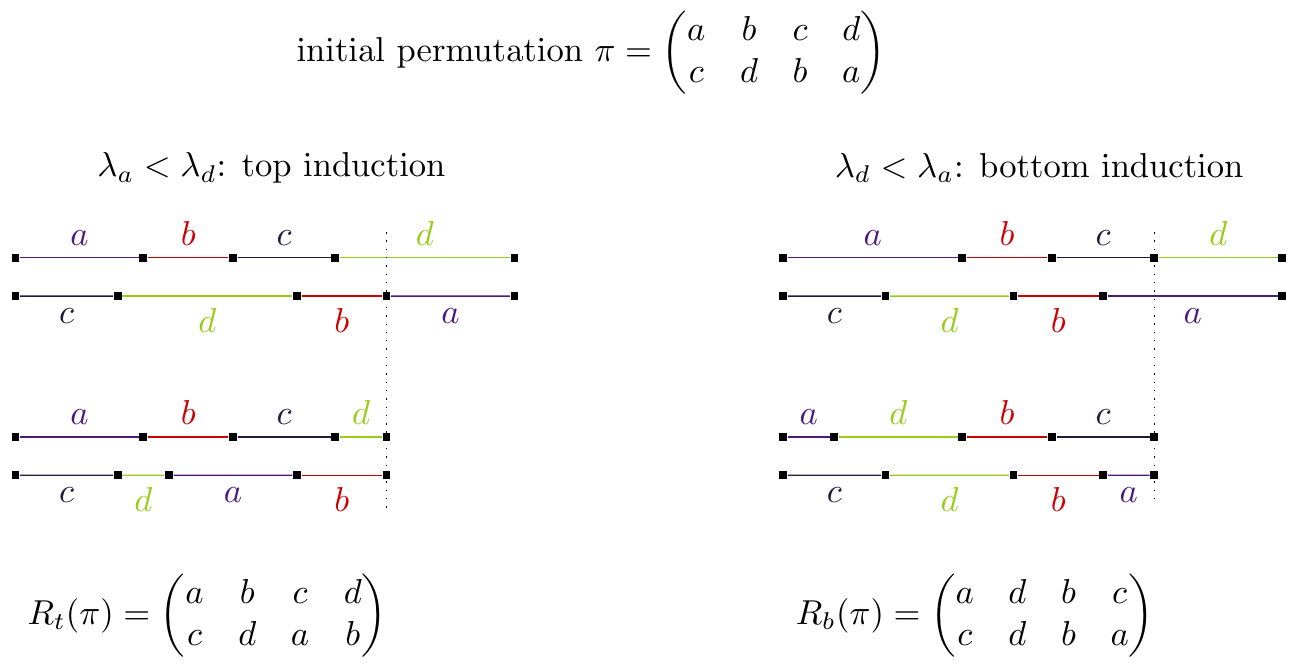}
\caption{The two alternatives of the Rauzy induction}
\label{fig:alternatives_of_Rauzy_induction}
\end{figure}

Let $T = T_{\lambda,\pi}$ be an interval exchange transformation on $I=[0,a)$ where $\pi$ is an irreducible labeled permutation on an alphabet $\mathcal{A}$ with $n$ letters and $\lambda \in \RR_+^{\mathcal{A}}$ satisfies $\sum_{\alpha \in \mathcal{A}} \lambda_\alpha = a$. For $i=1,\ldots,n+1$ we set $\displaystyle x_i = \sum_{j = 1}^{i-1} \lambda_{\pi_t^{-1}(j)}$ the discontinuities of $T$ and $\displaystyle y_i = \sum_{j = 1}^{i-1} \lambda_{\pi_b^{-1}(j)}$ the ones  of $T^{-1}$. We have $x_1 = y_1 = 0$ and $x_{n+1} = y_{n+1} = a$. Let $J=\left[0, \max(x_{n},y_{n})\right) \subset I$. The Rauzy induction of $T$, denoted by $R(T)$, is the interval exchange $T' = T_{\lambda',\pi'}$ obtained as the first returned map of $T$ on $J$. The type of $T$ is \emph{top} if $\lambda_{\pi_t^{-1}(n)} > \lambda_{\pi_b^{-1}(n)}$ and \emph{bottom} if $\lambda_{\pi_b^{-1}(n)} > \lambda_{\pi_t^{-1}(n)}$. In the case $\lambda_{\pi_t^{-1}(n)} = \lambda_{\pi_b^{-1}(n)}$ there is no Rauzy induction defined. When $T$ is of type top (resp. bottom) the label $\pi_t^{-1}(n)$ (resp. $\pi_b^{-1}(n)$) is called the \emph{winner} and $\pi_b^{-1}(n)$ (resp. $\pi_t^{-1}(n)$) the \emph{loser}. Let $\alpha \in \mathcal{A}$ (resp. $\beta \in \mathcal{A}$) be the winner (resp. loser) of $T$. The vector $\lambda'$ of interval lengths of $T'$ is given by
\begin{align*}
\lambda'_\alpha &= \lambda_\alpha - \lambda_\beta, \\
\lambda'_\nu &= \lambda_\nu \quad \text{for all $\nu \in \mathcal{A} \backslash \{\alpha\}$}.
\end{align*}
The permutation $\pi' = R_\varepsilon(\pi)$ is defined as follows, where $\varepsilon \in \{t,b\}$ is the type of $T$ ($t$ for top and $b$ for bottom). Let $\alpha_t = \pi_t^{-1}(n)$ (resp. $\alpha_b = \pi_b^{-1}(n)$) the label on the right of the top line (resp. bottom line). As $\pi$ is irreducible, the position $m=\pi_b(\alpha_t)$ of $\alpha_t$ on the bottom line (resp. $m = \pi_t(\alpha_b)$ of $\alpha_b$ in the top line) is different from $n$. We obtain $\pi'$ from $\pi$ by moving $\alpha_b$ (resp. $\alpha_t$) from position $n$ to position $m+1$ in the bottom interval (resp. top interval). The operations $R_t$ and $R_b$ are formally defined by \vspace{0.3cm}
\\
\begin{tabular}{p{7cm}@{\hspace{1cm}}|@{\hspace{1cm}}p{7cm}}
$R_t(\pi_t,\pi_b) = (\pi_t,\pi_b')$ where if $m = \pi_b(\alpha_t)$ we have
&
$R_b(\pi_t,\pi_b) = (\pi_t', \pi_b)$ where if $m = \pi_t(\alpha_b)$ we have
\\
$
{\pi'}_b^{-1}(j) = \left\{
\begin{array}{ll}
\pi_b^{-1}(j) & \text{if $j \leq m$,} \\
\pi_b^{-1}(n) & \text{if $j = m+1$,} \\
\pi_b^{-1}(j-1) & \text{otherwise.}
\end{array} \right.
$
&
$
{\pi'}_t^{-1}(j) = \left\{
\begin{array}{ll}
\pi_t^{-1}(j) & \text{if $j \leq m$,} \\
\pi_t^{-1}(n) & \text{if $j = m+1$,} \\
\pi_t^{-1}(j-1) & \text{otherwise.}
\end{array} \right.
$ 
\end{tabular} \vspace{0.3cm} \\
The map $R_t$ and $R_b$ are called \emph{Rauzy moves}. An example of a Rauzy induction of an interval exchange transformation is shown in Figure~\ref{fig:alternatives_of_Rauzy_induction}. The Rauzy moves on reduced permutations are defined using the section $\pi \mapsto (id,\,\pi)$ and the projection $(\pi_t,\pi_b) \mapsto \pi_b \circ \pi_t^{-1}$ introduced in Section~\ref{subsubsection:labeled_permutations}.

We consider one more operation called \emph{inversion} and denoted by $s$ which reverses the top and the bottom and the left and the right
\[
s \left(
\begin{array}{llll}
a_1 & a_2 & \ldots & a_n \\
b_1 & b_2 & \ldots & b_n
\end{array}
\right)
=
\left(
\begin{array}{llll}
b_n & \ldots & b_2 & b_1 \\
a_n & \ldots & a_2 & a_1
\end{array}
\right).
\]

The following is standard.
\begin{lemma}
The Rauzy moves $R_t$, $R_b$ and the inversion $s$ preserve irreducible permutations. The Rauzy moves and the symmetry restricted to the set of irreducible permutations are bijections. 
\end{lemma}

\begin{definition}
\label{def:rauzy_class_combinatorial}
Let $\pi$ be an irreducible permutation. The orbit of $\pi$ under the action of $R_t$ and $R_b$ (resp. $R_t$, $R_b$ and $s$) is called the \emph{Rauzy class} (resp. \emph{extended Rauzy classes}) of $\pi$ and note it $\mathcal{R}(\pi)$. The \emph{Rauzy diagram} (resp. \emph{extended Rauzy diagram}) of $\pi$ is the labeled oriented graph with vertices $\mathcal{R}(\pi)$ and edges corresponding to the action of $R_t$ and $R_b$ (resp. $R_t$, $R_b$ and $s$).

Let $\pi$ is a reduced (resp. labeled) permutation than the Rauzy class of $\pi$ is called a \emph{reduced Rauzy class} (resp. \emph{labeled Rauzy class}).
\end{definition}

The standard permutations play a central role in Rauzy classes in particular we have.
\begin{proposition}[\cite{Rauzy79}]
Every Rauzy class contains a standard permutation.
\end{proposition}

\begin{proof}
Let $\mathcal{R}$ be a Rauzy class of permutations on $n$ letters and let $\pi \in \mathcal{R}$. Let $\alpha_t = \pi_t^{-1}(n)$ and $\alpha_b = \pi_b^{-1}(n)$ be the labels of the right extremities. Let $n_b = \pi_b(\alpha_t)$ and $n_t = \pi_t(\alpha_b)$.

If $n_t = \min(n_t,n_b) \not= 1$ then by irreducibility, in the set $\pi_t \circ \pi_b^{-1}(\{n_t+1,\ldots,n\})$ the minimal element is less than $n_t$. Let $n'_b$ be this minimum and $\alpha'_b$ be the letter for which the minimum is reached. Applying $R_t$ we can move the letter $\alpha'_b$ at the right extremity of the bottom line. After this first step the quantity $n'_b = \min(n_t,n'_b)$ is lesser than $n_t = \min(n_t,n_b)$. For, the case $n_b = \min(n_t,n_b)$, we use $R_b$ to decrease the quantity $\min(n_t,n_b)$.

Iterating succesively $R_t$ or $R_b$ as in the above step, we obtain a permutation such that either $n_t = 1$ or $n_b = 1$. Applying one more time a Rauzy move, we obtain both equal to $1$.
\end{proof}

\begin{figure}[ht!]
\centering
\subfloat[The Rauzy diagram of $(4321)$]{
\picinput{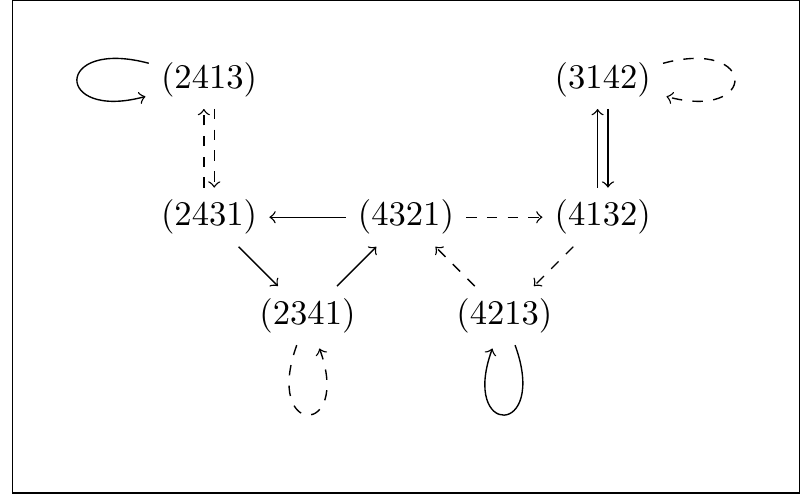}
\label{fig:rauzy_diagram_2}}
\subfloat[The Rauzy diagram of $(4231)$]{
\picinput{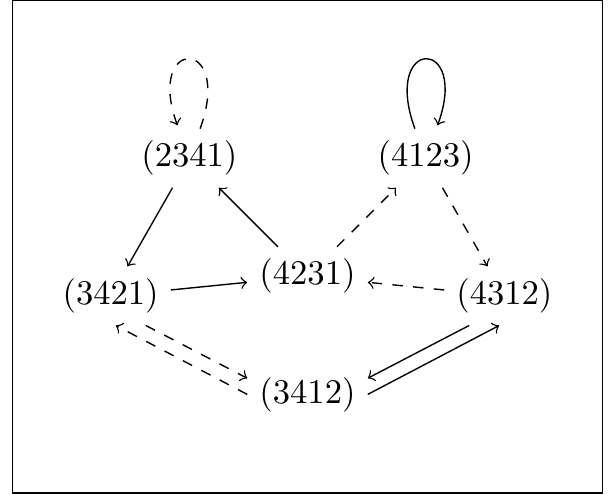}
\label{fig:rauzy_diagram_0_0}}
\caption{The two Rauzy diagrams of $S_4^o$.}
\label{fig:rauzy_diagrams_of_perms_of_length_4}
\end{figure}

There are only two standard permutations of length $4$, $(4\,3\,2\,1)$ and $(4\,2\,3\,1)$, which define two Rauzy classes. Their Rauzy diagrams are presented in Figure~\ref{fig:rauzy_diagrams_of_perms_of_length_4}.

The labeled rauzy diagrams are coverings of reduced rauzy diagrams (the covering map is the projection $(\pi_t,\pi_b) \mapsto \pi_b \circ \pi_t^{-1}$). The degree of the covering which gives the multiplicative coefficient between the cardinality of reduced Rauzy classes and labeled Rauzy classes and its computation involves geometric methods which are developed in \cite{Boissy10}.

\subsubsection{Examples of Rauzy diagrams}
\label{subsubsection:Rauzy_diagrams_of_rotations}
\label{subsubsection:symmetric_permutation}

We denote by $\pi_n^{sym}$ the \emph{symmetric permutation} on $n$ letters defined by $\pi_n^{sym}(k) = n-k+1$ for $k=1,\ldots,n$. In our notation, $\pi_n^{sym}$ writes
\begin{equation}
\pi_n^{sym} =\left( \begin{array}{cccc}
1 & 2 & \ldots & n \\
n & n-1 & \ldots & 1
\end{array} \right).
\label{eq:definition_permutation_symmetric}
\end{equation}
The permutation $\pi_n^{sym}$ has a Rauzy class which is described in \cite{Rauzy79} (see also \cite{Yoccoz05} p. 53).

\begin{proposition} \label{prop:rauzy_class_of_symmetric_permutation}
The Rauzy class $\mathcal{R}^{sym}_n$ of $\pi_n^{sym}$ coincide with its extended Rauzy class. It contains $2^{n-1}-1$ permutations and among them only $\pi_n^{sym}$ is a standard permutation.
\end{proposition}
In that case we remark that the labeled Rauzy class coincide with the reduced one.

We now describe an other class. Let $n$ be a positive integer and let
\begin{equation}
\pi_n^{rot} = \left( \begin{array}{cccccc}
1 & 2 & 3 & \ldots & n-1 & n \\
n & 2 & 3 & \ldots & n-1 & 1
\end{array} \right).
\label{eq:definition_of_permutation_of_rotation}
\end{equation}
The permutation $\pi_n^{rot}$ is called of \emph{rotation class}. Any interval exchange transformation with permutation $\pi_n^{rot}$ is a first return map of a rotation. We denote by $\mathcal{R}_n^{rot}$ the Rauzy diagram of $\pi_n^{rot}$

\begin{figure}[!ht]
\centering
\subfloat[$n=2$]{
\picinput{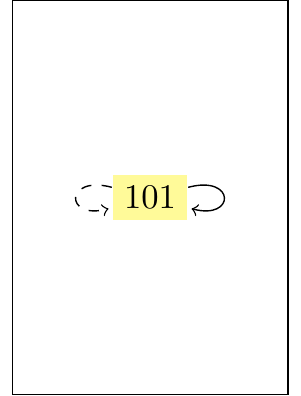}
\label{fig:rauzy_diagram_torus1.tex}}
\subfloat[$n=3$]{
\picinput{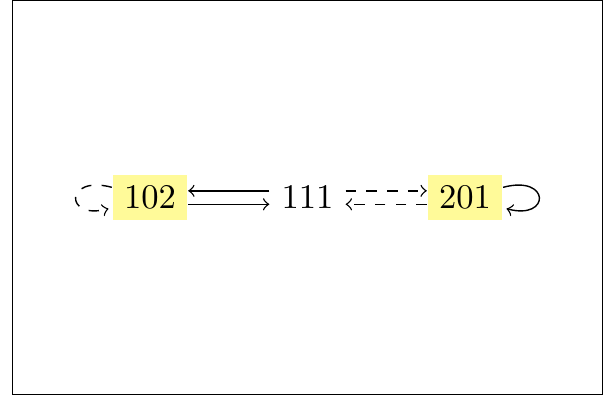}
\label{fig:rauzy_diagram_torus2.tex}}
\subfloat[$n=4$]{
\picinput{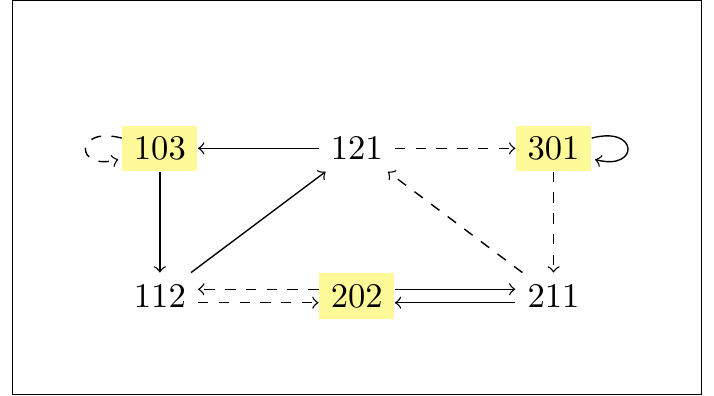}
\label{fig:rauzy_diagram_torus3.tex}}

\subfloat[$n=5$]{
\picinput{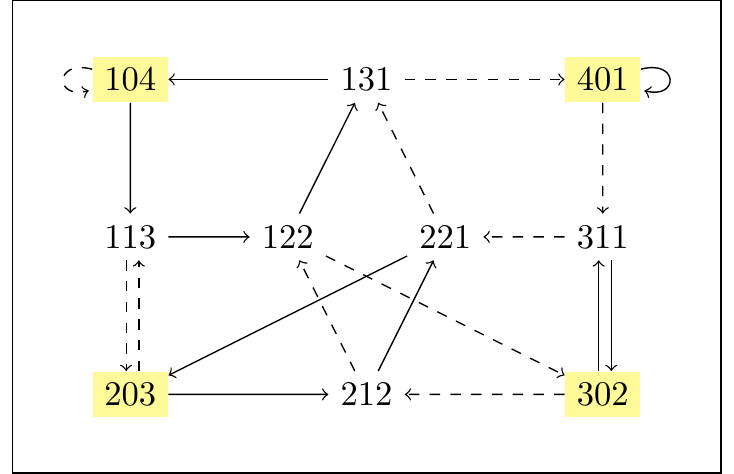}
\label{fig:rauzy_diagram_torus4.tex}}
\hspace{1cm}
\subfloat[$n=6$]{
\picinput{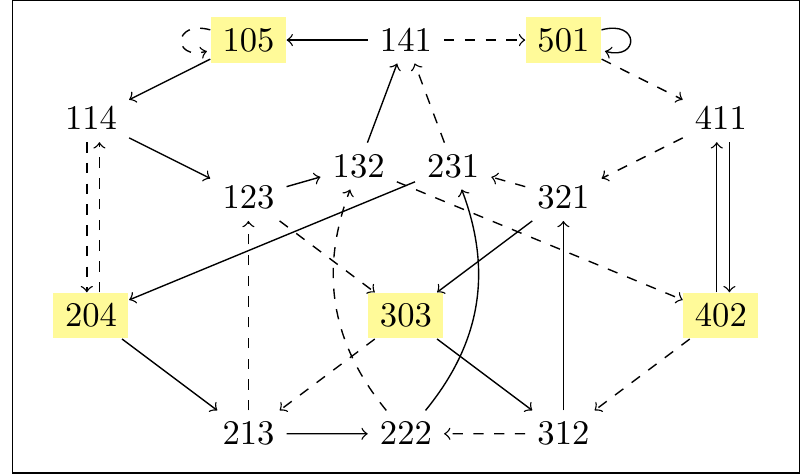}
\label{fig:rauzy_diagram_torus5.tex}}
\caption{The graphs $\mathcal{G}_n$ for $n=2,3,4,5,6$.}
\label{fig:rauzy_diagrams_torus}
\end{figure}

We now build a graph $\mathcal{G}_n$. Let $V_n = \{(a,b,c) \in \NN^3;\ a,c\geq 1,\ b\geq 0\ \text{and}\ a+b+c = n\}$. From a triple $(a,b,c) \in V_n$ we define the permutation
\begin{equation} \label{eq:definition_pi_abc}
\pi(a,b,c) = 
\left( \begin{array}{ccc|ccc|ccc}
1 & \ldots &  a          & a+1 & \ldots & a+b  & a+b+1 & \ldots & a+b+c \\
a+b+1 & \ldots & a+b+c   & a+1 & \ldots & a+b  & 1 & \ldots & a
\end{array} \right).
\end{equation}
Let $\mathcal{G}_n$ be the oriented labeled graph with vertices $V_n$ and edges are of two types
\begin{itemize}
\item the \textit{left edges} are $(a,0,c) \to (1,a-1,c)$ and if $b \not= 0$, $(a,b,c) \to (a+1,b-1,c)$,
\item the \textit{right edges} are $(a,0,c) \to (a,c-1,1)$ and if $b \not= 0$ $(a,b,c) \to (a,b-1,c+1)$.
\end{itemize}
From the rules that define the edges, we see that each vertex has exactly one incoming and one outgoing edge of each type. Moreover, in each cycle made by left edges (resp. right edges) there is exactly one element of the form $(a,0,c)$. The number $a$ (resp. $c$) is the length of the cycle. In Figure~\ref{fig:rauzy_diagrams_torus} we draw examples of such graphs.

\begin{proposition}
\label{thm:rauzy_diagram_torus}
The graph $\mathcal{G}_n$ is isomorphic to the Rauzy diagram $\mathcal{R}_n^{rot}$ under the map $(a,b,c) \mapsto \pi(a,b,c)$. The left edges (resp. right edges) in $\mathcal{G}^{rot}_n$ correspond to top (resp. bottom) Rauzy moves in $\mathcal{R}_n^{rot}$.

\noindent Moreover the extended Rauzy diagram of $\pi_n^{rot}$ has the same set of vertices as $\mathcal{R}_n^{rot}$. The action of $s$ in the extended Rauzy class corresponds to $(a,b,c) \mapsto (c,b,a)$ in $\mathcal{G}_n$.
\end{proposition}
We remark that for $\pi_n^{rot}$ the ratio between the cardinalities of labeled and reduced Rauzy classes is $(n-1)!$. This result is a particular instance of a theorem of \cite{Boissy10}.

\begin{proof}
The permutation $\pi_n^{rot}$ corresponds to the triple $(1,n-2,1) \in V_n$. From the definition of $R_t$ and $R_b$ it can be easily checked that the edges of $\mathcal{G}_n$ corresponds to Rauzy moves on $\pi(a,b,c)$. Hence, the set of permutations associated to $V_n$ is invariant under Rauzy induction. As the graph $\mathcal{G}_n$ is connected, this set is exactly the Rauzy class of $\pi_n^{rot}$.

The inversion $s$ exchanges the three parts of the permutation $\pi(a,b,c)$ delimited by the bars in (\ref{eq:definition_pi_abc}). The structure of the permutation in three blocks is preserved and we get that $s \cdot \pi(a,b,c) = \pi(c,b,a)$.
\end{proof}

\begin{proposition}
The Rauzy class $\mathcal{R}^{rot}_n$ of $\pi_n^{rot}$ coincide with its extended Rauzy class. It contains $\binom{n}{2} = \frac{n(n-1)}{2}$ permutations and among them only $\pi_n^{rot}$ is a standard permutation.
\end{proposition}

\subsection{From permutations to translation surfaces}
\label{subsection:from_permutations_to_surfaces}
\subsubsection{Translation surface}
\label{subsubsection:translation_surfaces}
Let $S$ be a compact oriented connected surface. A \emph{translation structure} on $S$ is a flat metric defined on $S-\Sigma$ where $\Sigma \subset S$ is a finite set of points which has trivial holonomy (the parallel transport along a loop is trivial). The latter condition implies that at any point $P \in \Sigma$ the metric has a conical singularity of angle an integer multiple of $2\pi$ : the length of a circle centered at a conic point of angle $2\pi m$ with small radius $r$ will not measure $2 \pi r$ but $2 \pi m r$. More concretely, a translation surface can be built from gluing polygons. Let $P_1, \ldots, P_f \subset \RR^2$ be a finite collection of polygons and $\tau$ a pairing of their sides such that each pair is made of two sides which are parallel, with the same length and opposite normal vectors. We define the equivalence relation $\sim_\tau$ on the union $P = \cup P_i$: $x_1 \sim_\tau x_2$ if $x_1$ and $x_2$ are, respectively, on two sides $s_1$ and $s_2$ which are paired by $\tau$ and differ by the unique translation that maps $s_1$ onto $s_2$. The quotient $S = S(P_i,\tau) := P / \sim_\tau $ is a translation surface for which the metric and the vertical direction are induced from $\RR^2$. We call the couple $(P_i,\tau)$ a \emph{polygonal representation} of the translation surface $S$. Reciprocally, any translation surface admits a geodesic triangulation which gives a polygonal representation of the surface.
\begin{figure}[h!]
\centering
\subfloat[A regular octogon]{
\picinput{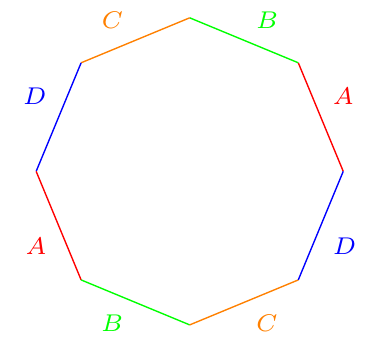}
\label{fig:polygons_gluing_octogon}}
\hspace{1cm}
\subfloat[An L-shaped polygon]{
\picinput{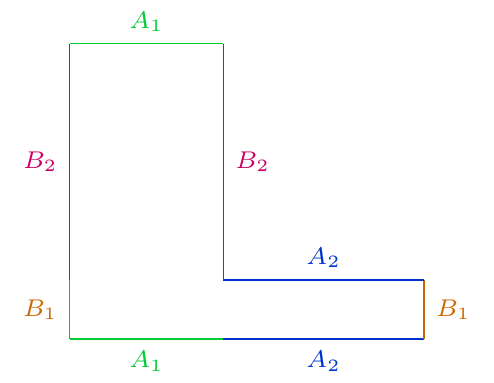}
\label{fig:polygons_gluing_l_shaped}}
\hspace{1cm}
\subfloat[Two pentagons]{
\picinput{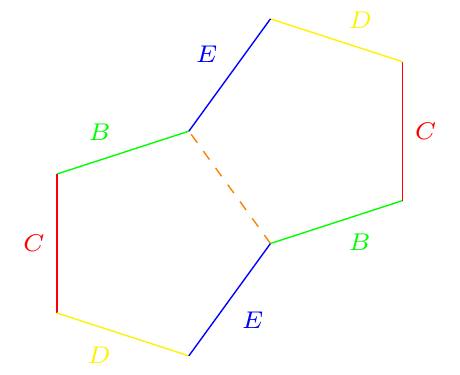}
\label{fig:polygons_gluing_bi_pentagon}}
\caption{Three surfaces built from polygons. The pairings are defined by colors and labels.}
\label{fig:polygons_gluing}
\end{figure}

Let $S$ be a translation surface and $(2\pi n_1, 2\pi n_2, \ldots, 2\pi n_k)$ the list of angles of its conical singularities. The genus $g$ of the surface satisfies
\begin{equation}
  \label{eq:genus_and_conic_points}
  2g - 2 = \sum_{i=1}^k (n_i - 1).
\end{equation}
The integer partition $p = (n_1,n_2,\ldots,n_k)$ is the \emph{profile} of the translation surface $S$ and Equation~(\ref{eq:genus_and_conic_points}) resumes to $s(p)-l(p)=2g-2$ where $s(p) = n_1 + \ldots +n_l$ is the sum of the terms of $p$ and $l(p)=k$ its length. As a consequence the number of even terms in $p$ is even. This is the unique obstruction for a profile of a flat surface: for any integer partition $p$ such that the number of even terms is even there exists a translation surface $S$ with profile $p$. 

The genus, related to the collection of angles in Equation~(\ref{eq:genus_and_conic_points}), can also be deduced from the way the polygons are glued together. Let $f$ (for faces) be the number of polygons. Each pair of sides gives an embedded geodesic segment in the surface, let $e$ (for edges) be the number of those pairs. The vertices of the polygons are identified in a certain number of classes depending on the combinatorics of the pairing $\tau$, let $v$ (for vertices) be the number of classes. Then we have
\begin{equation}
\label{eq:genus_maps}
2 - 2g = v - e + f \qquad \text{where $g$ is the genus of the surface}.
\end{equation}

Consider the example of Figure~\ref{fig:polygons_gluing_octogon}, the surface obtained from the octogon has four edges and one vertex, thus $2-2g = 1 - 4 + 1 = - 2$, therefore its genus is $2$. On other hand, the angle at the unique conic point of the surface is $6\pi$. The two other examples of Figure~\ref{fig:polygons_gluing} have the same profile.


If a translation surface has a conical angle of $2\pi$ then, from the viewpoint of the metric, the singularity is removable: there exists a unique continuous way to extend the metric at this point. To a surface $S$ with profile $(n_1,n_2,\ldots,n_l,1,\ldots,1)$ with $k$ parts equal to $1$ we associate a surface $\overline{S}$ with profile $(n_1,n_2,\ldots,n_l)$. We say that surface $\overline{S}$ is obtained from $S$ by marking $k$ points.

\subsubsection{Moduli space of translation surfaces}
\label{subsubsection:moduli_space}
Two translation surfaces $S_1$ and $S_2$ are \emph{isomorphic} if there exists an orientation preserving isometry between $S_1$ and $S_2$ which maps the vertical direction of $S_1$ on the vertical direction of $S_2$. Let $\Omega \mathcal{M}(n_1-1, n_2-1, \ldots, n_k-1)$ be the collection of isomorphism classes of flat surfaces for whose profile is $(n_1, n_2, \ldots, n_k)$. The notation $\Omega \mathcal{M}(\kappa)$ comes from algebraic geometry where $\Omega \mathcal{M}$ is the tangent bundle to the moduli space of complex curves $\mathcal{M}_g$. In this settings, translation surfaces are considered as Riemann surfaces together with an Abelian differential. A conical singularity of angle $2 \pi m$ for the flat metric corresponds to a zero of degree $m-1$ of the Abelian differential (see \cite{Zorich06} for more details about the relations between flat structure and Abelian differential).

\begin{figure}[!ht]
\centering
\picinput{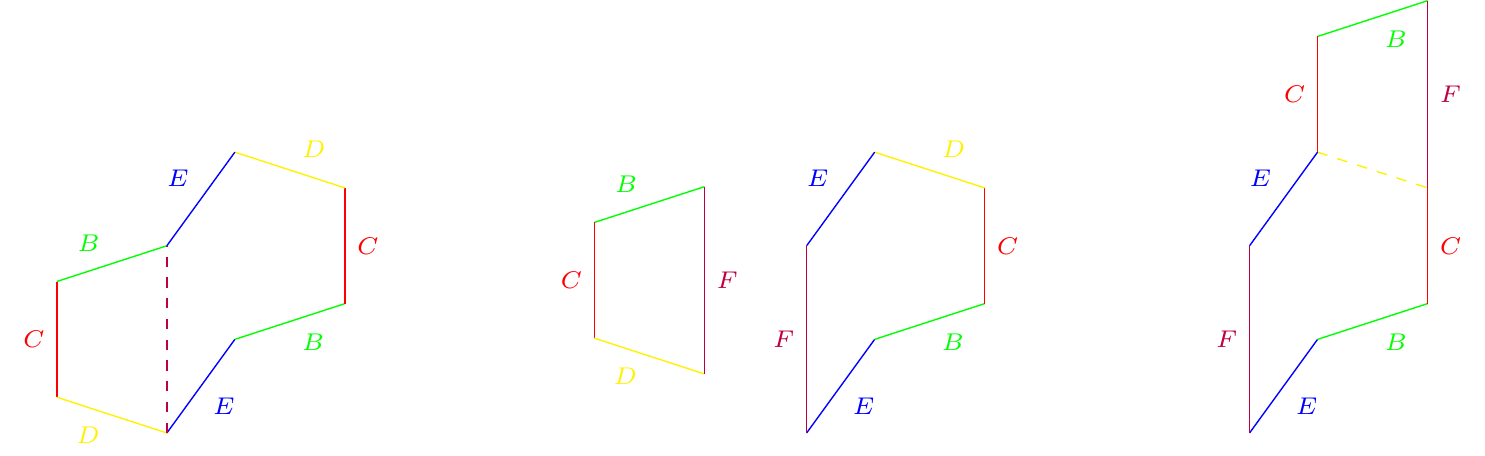}
\caption{An example of one cut followed by one paste.}
\label{fig:cut_and_paste_operations}
\end{figure}

We now define a topology on $\Omega \mathcal{M}$ using the construction with polygons. We first remark that given the combinatorics of polygons (e.g. the cyclic order of the edges in each polygon, and the pairing $\tau$), the set of vectors that are admissible as sides for the polygons forms an open set in $(\RR^2)^{e-f+1} = (\RR^2)^{2g+s-1}$ where as before $v$, $e$ and $f$ denote the number of vertices, edges and faces in the polygon. On other hand, two different polygonal representations may give isomorphic translation surfaces. We consider, on polygonal representations, the following operations (see also Figure~\ref{fig:cut_and_paste_operations})
\begin{itemize}
\item The \emph{cut operation} consists in the creation of a new pair of edges between two vertices (if it is possible). This operation creates an edge and the number of faces increases by $1$. 
\item The \emph{paste operation} consists in pasting two polygons along two edges which are paired. This operation delete an edge and the number of faces decreases by $1$.
\end{itemize}
We have the following.
\begin{proposition}[\cite{Masur82},\cite{Veech93}]
The isomorphism class of a surface $S(P_i,\tau)$ built from polygons is invariant under cut and paste operations of the polygonal representation $(P_i,\tau)$. Moreover, if $P$ and $P'$ are two polygonal representations of the same surface $S$ then there exists a sequence $P_0 = P$, $P_1$, \ldots, $P_n= P'$ of polygonal representations such that $P_{i+1}$ is obtained from $P_i$ either by a cut or a paste operation.
\end{proposition}
The above proposition states that the space $\Omega \mathcal{M}$ can be considered as a quotient of a finite union of open sets of $(\RR^2)^{2g-2+v}$ by the action of cutting and pasting. The topology of $\Omega \mathcal{M}$ is by definition the quotient topology. As the action of cut and paste operations is discrete, the local system of neighborhood in $\Omega \mathcal{M}$ are open sets in vector spaces. Hence, two translation surfaces are near if they admit decompostions in polygons which have the same combinatorics and roughly the same shape.

We drafted a construction of the moduli space of translation surfaces $\Omega \mathcal{M}$ which is a quotient of the tangent bundle of a Teichm\"uller space (which corresponds to polygonal representation) by the mapping class group (which corresponds to cut and paste operations). See \cite{Masur82} and the textbooks \cite{Ahlfors66},  \cite{Nag88}, \cite{ImayoshiTaniguchu92} or \cite{Hubbard06}.

\subsubsection{Suspension of a permutation and Rauzy-Veech induction}
\label{subsubsection:suspensions}
\begin{figure}[h!]
\centering
\picinput{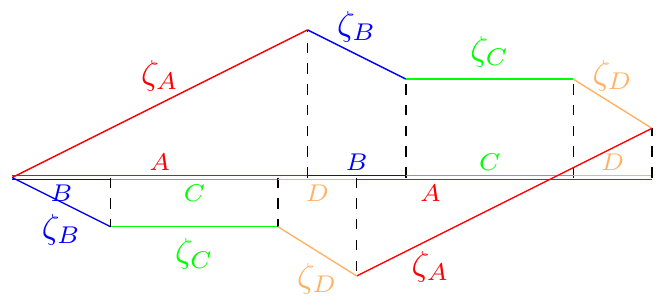}
\caption{A suspension of the permutation $\binom{A B C D}{ B C D A}$ and the first return map of the vertical linear flow on its canonical transverse segment.}
\label{fig:suspension}
\end{figure}

We recall the method in \cite{Veech82} for building a translation surface from a permutation. The version for labeled permutations comes from \cite{MarmiMoussaYoccoz05} and \cite{Bufetov06}.
Let $\pi = (\pi_t,\pi_b)$ be an irreducible labeled permutation, $\mathcal{A}$ its alphabet and $n = |\mathcal{A}|$. A \emph{suspension datum} for $\pi$ is a collection of vectors $\zeta = (\zeta_\alpha)_{\alpha \in \mathcal{A}} = ((\lambda_\alpha, \tau_\alpha))_{\alpha \in \mathcal{A}} \in (\RR_+\times \RR)^\mathcal{A}$ such that 
\[
\forall 1 \leq k \leq n-1,\
\sum_{\pi_t(\alpha) \leq  k} \tau_\alpha > 0 
\quad \text{and} \quad
\sum_{\pi_b(\alpha) \leq k} \tau_\alpha < 0.
\]
To each suspension datum $\zeta$ we associate a translation surface $S = S(\zeta,\pi)$ in the following way. Consider the broken lines $L_t$ (resp. $L_b$) in $\RR^2$ starting at the origin and obtained by the concatenation of the vectors $\zeta_{\pi_t^{-1}(j)}$ (resp. $\zeta_{\pi_b^{-1}(j)}$) $j=1,\ldots,n$ (in this order). If the broken line $L_t$ and $L_b$ have no intersection other than the endpoints, we can construct a translation surface $S$ from the polygon bounded by $L_t$ and $L_b$. The pairing of the sides associate to the side $\zeta_\alpha$ of $L_t$ the side $\zeta_\alpha$ of $L_b$ (see Figure~\ref{fig:suspension}). Note that the lines $L_t$ and $L_b$ might have some other intersection points. But in this case, one can still define a translation surface using the \emph{zippered rectangle construction} due to \cite{Veech82}. In the suspension $S=S(\pi,\zeta)$ there is a canonical embedding of the segment $I = \left[0, |\lambda|\right)$. The first return map on $I$ of the translation flow of $S$ is the interval exchange map $T$ with permutation $\pi$ and vector of lengths $\lambda$ (see Figure~\ref{fig:suspension}).
The Rauzy induction can be extended to suspensions and will be still denoted by $R$. If $\zeta = (\lambda,\tau)$ is a suspension data for $\pi$, then $R(\zeta,\pi)$ is the suspension $(\zeta',\pi') = ((\lambda',\tau'),\pi')$ where
\begin{itemize}
\item $\pi' = R_\varepsilon(\pi)$ where $\varepsilon \in \{t,b\}$ is the type of $T$,
\item $\zeta'_\alpha = \zeta_\alpha - \zeta_\beta$ where $\alpha$ (resp. $\beta$) is the winner (resp. loser) for $T$.
\end{itemize}
This extension is known as the \emph{Rauzy-Veech induction}, and is used as a discretization of the Teichm\"uller flow.

\begin{figure}[!ht]
\begin{center}
\picinput{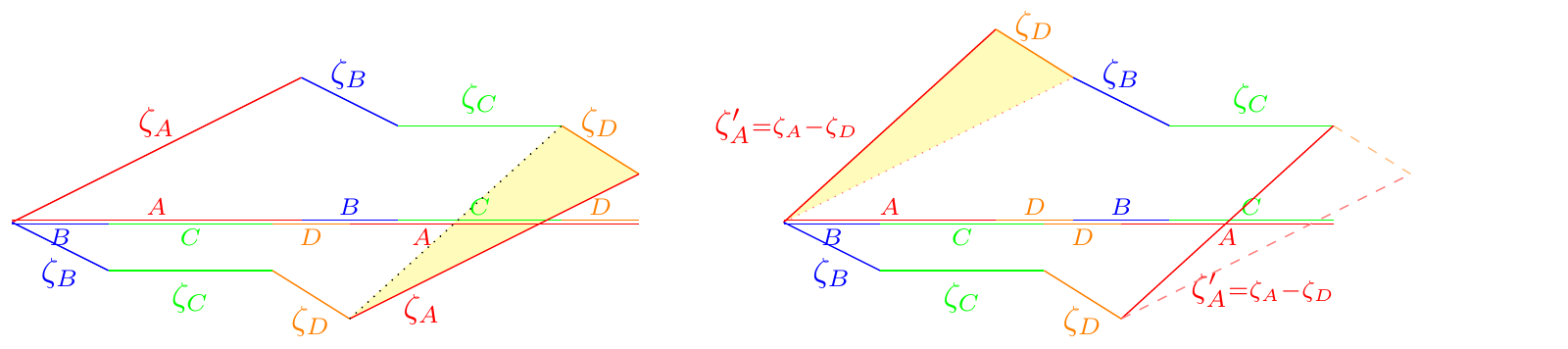}
\end{center}
\caption{The (bottom) Rauzy-Veech induction on a suspension of $\pi=\binom{A B C D}{B C D A}$}
\label{fig:Rauzy_Veech_induction}
\end{figure}

By construction the surfaces $S_{\zeta,\pi}$ and $S_{\zeta',\pi'}$ are isomorphic: the Rauzy-Veech induction corresponds to one cut followed by one paste operations (see Figure~\ref{fig:Rauzy_Veech_induction}). In particular, by the definition of Rauzy class (Definition~\ref{def:rauzy_class_combinatorial}), we have the following proposition which is a key ingredient in the correspondance between Rauzy classes and moduli space of translation surfaces.
\begin{proposition}[\cite{Veech82}]
\label{prop:associated_connected_component_of_a_rauzy_class}
Let $\mathcal{R}$ be a Rauzy class or an extended Rauzy class. Then, the set of suspensions obtained from permutations in $\mathcal{R}$ is open and connected in $\Omega \mathcal{M}$.
\end{proposition}
The case of extended Rauzy class in the above proposition follows from the fact that the involution $s$ on permutations (see Section~\ref{subsubsection:combinatorial_rauzy_classes}) can be seen as a central symmetry of the suspension $S(\pi,\zeta)$.

\subsection{Permutation invariants of Rauzy classes} \label{subsection:invariants_for_Rauzy_classes}
We now define the three invariants of permutations that lead to a classification of Rauzy classes.

\subsubsection{Interval diagram and profile}
\label{subsubsection:profile}
\begin{figure}[!h]
\centering
\picinput{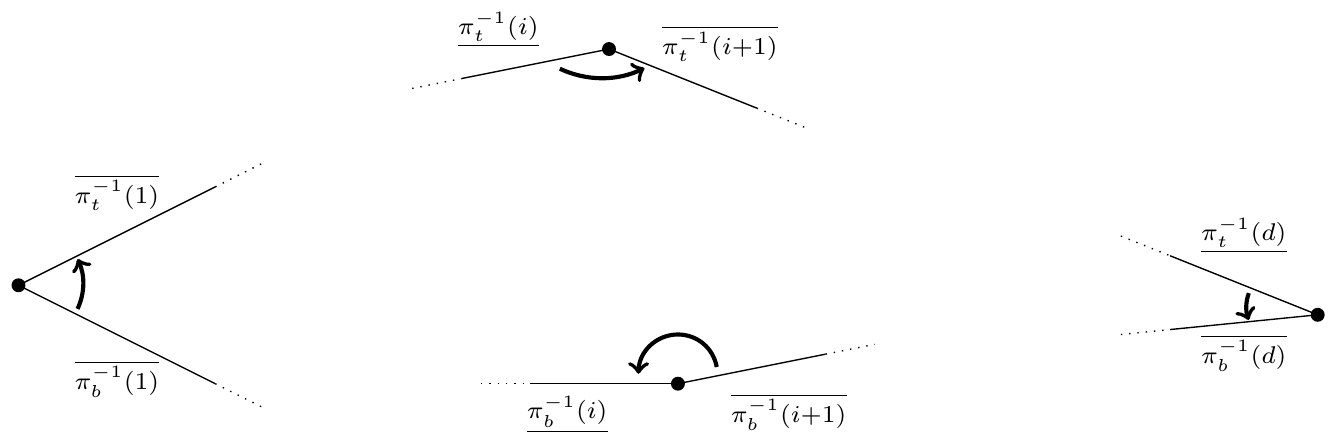}
\caption{The permutation $\tilde{\sigma_\pi}$ and turning around vertices of a suspension of $\pi$.}
\label{fig:interval_diagram}
\end{figure}

Let $\pi$ be a labeled permutation with alphabet $\AAA$. We consider a refinement of the permutation $\sigma$ introduced in \cite{Veech82} which take care of the labels of $\pi$. Let $\tilde{\sigma}$ be the permutation on the set $\outgos{\AAA} \cup \incoms{\AAA} = \{\outgos{a};\ a\in \AAA\} \cup \{\incoms{a};\ a \in \AAA\}$ defined by
\[
\tilde{\sigma}(\outgos{a}) = 
\left\{ \begin{array}{c@{,\quad}c}
\outgos{\pi_t^{-1}(1)} & \text{if $\pi_b(a) = 1$} \\
\incoms{\pi_b^{-1}(\pi_b(a)-1)} & \text{if $\pi_b(a) \not= 1$}
\end{array} \right.
\quad \text{and} \quad
\tilde{\sigma}(\incoms{a}) =  
\left\{ \begin{array}{c@{,\quad}c}
\outgos{\pi_t^{-1}(\pi_t(a)+1)} & \text{if $\pi_t(a) \not= n$} \\
\incoms{\pi_b^{-1}(n)} & \text{if $\pi_t(a) = n$} 
\end{array} \right. .
\]
Assume that the permutation $\pi$ is irreducible and consider a suspension $S$ of $\pi$. We identify $\outgos{a}$ (resp. $\incoms{a}$) to the left-half (resp. right-half) of the edge labeled $a$ in $S$. The permutation $\tilde{\sigma}$ corresponds to the sequence of half-edges that we cross by turning around vertices of $S$ (see Figure~\ref{fig:interval_diagram}).

Let $\pi$ be a labeled permutation on $\AAA$. We define $\outgos{\AAA_\pi}$ (resp. $\incoms{\AAA_\pi}$) to be the quotient of $\outgos{\AAA}$ (resp. $\incoms{\AAA}$) in which $\outgos{\pi_b^{-1}(1)}$ and $\outgos{\pi_t^{-1}(1)}$ (resp. $\incoms{\pi_t^{-1}(n)}$ and $\incoms{\pi_b^{-1}(n)}$) are identified.
\begin{definition}
\label{def:interval_diagram}
The \emph{interval diagram} of $\pi$ is the permutation $\sigma = \sigma_\pi$ on the set $\mathcal{A}_\pi = \outgos{\mathcal{A}_\pi} \cup \incoms{\mathcal{A}_\pi}$ defined by
\[
\sigma_\pi(a) = \left\{
\begin{array}{ll}
\tilde{\sigma}(\pi_t^{-1}(1)) & \text{if $a = \outgos{(\pi_b^{-1}(1),\pi_t^{-1}(1))}$}, \\
\tilde{\sigma}(\pi_b^{-1}(n)) & \text{if $a = \incoms{(\pi_t^{-1}(n),\pi_b^{-1}(n))}$}, \\
\tilde{\sigma}(a) & \text{otherwise}.
\end{array} \right.
\]
\end{definition}
As an example, on the permutation $\left(\begin{smallmatrix}A&B&C&D\\B&C&D&A\end{smallmatrix}\right)$ in Figures~\ref{fig:suspension} and \ref{fig:Rauzy_Veech_induction} the interval diagram is
\[
\sigma_\pi = \left(\outgos{(B,A)},\incoms{(D,B)}\right)
\left(\outgos{C},\incoms{B}\right)
\left(\outgos{D},\incoms{C}\right)
\]
The interval diagram $\sigma_\pi$ exchanges $\outgos{\mathcal{A}_\pi}$ and $\incoms{\mathcal{A}_\pi}$. In particular, the permutation $\sigma_\pi^2$ can be written as a product of two permutations $\outgos{\sigma_\pi}$ and $\incoms{\sigma_\pi}$ on respectively $\outgos{\mathcal{A}_\pi}$ and $\incoms{\mathcal{A}_\pi}$.

We recall that conjugacy class of permutations of a set with $n$ elements are in bijection with integer partition of $n$. To a permutation $\sigma$ we associate the length of the cycles in the disjoint cycle decomposition of $\sigma$.
\begin{lemma}
Let $\pi$ be an irreducible permutation and $S$ a suspension of $\pi$. The profile of $S$ is the integer partition associated to the conjugacy class of the permutation $\outgos{\sigma_\pi}$ (or $\incoms{\sigma_\pi}$).
\end{lemma}

\begin{proof}
Following \cite{Veech82} and \cite{Boissy10}, the permutation $\outgos{\sigma_\pi}$ (resp. $\incoms{\sigma_\pi}$) can be seen as the crossing of the horizontal direction. In particular each cycle corresponds to a conical singularity of the suspension $S$ and its length $k$ equals the angle divided by $2 \pi$.
\end{proof}

\subsubsection{The spin parity}
\label{subsubsection:spin_invariant}
Now we define the spin parity of a permutation $\pi$ whose profile $p$ contains only odd numbers. As the spin parity relies on the classification of quadratic forms over the field with two elements $\FF_2$, we first recall this classification in Theorems~\ref{thm:classification_of_non_degenerate_quadratic_binary_forms} and \ref{thm:classification_of_all_quadratic_binary_forms}. For more details about the spin invariant see \cite{Johnson80} and \cite{KontsevichZorich03}.

Let $n \geq 1$ and $V$ a vector space over $\FF_2$. A \emph{quadratic form} on $V$ is a map $q: V \rightarrow \FF_2$ which is an homogeneous polynomial of degree $2$ in any coordinate system of $V$. If $q$ is a quadratic form, then the application $B_q$ defined on $V \times V$ by $B_q(u,v) = q(u+v) - q(u) - q(v)$ is bilinear. The form $q$ is called \emph{nondegenerate} if $B_q$ is nondegenerate. Because the characteristic is two, the form $B_q$ satisfies 
\begin{equation} \label{eq:sym_char2}
B_q(u,v) = B_q(v,u) \quad \text{and} \quad B_q(u,u) = 0.
\end{equation}.
If there exists a non degenerate bilinear form $B$ on $V$ which satisifies (\ref{eq:sym_char2}) then the dimension of $V$ is even. We consider from now that the dimension $n=2g$ is even and $V = (\FF_2)^{2g}$. On $V$, there is only one linear equivalence class of nondegenerate bilinear form $B$ that satisfies (\ref{eq:sym_char2}). The \emph{standard nondegenerate bilinear} form on $V=\FF_2^{2g}$ is the bilinear form $B_0$ given in coordinates $v=(x_1,y_1,\ldots,x_g,y_g) \in V$, $v'=(x'_1,y'_1,\ldots,x'_g,y'_g)$ by
\[
B_0(v,v') = \sum_{i=1}^g (x_i y'_i + x'_i y_i).
\]
By the above remark, any non degenerate quadratic form is linearly equivalent to one whose associated bilinear form is $B_0$. In order to classify quadratic form up to linear equivalence, we assume that $q$ is such that $B_q = B_0$. In other words the quadratic form $q$ writes in terms of the coordinates of $v$ as
\begin{equation}
\label{eq:quad_form_over_F2}
q(v) = \sum_{i=1}^n \left( a_i x_i^2 + b_i y_i^2 + x_i y_i \right),
\end{equation}
where $t = ((a_i,b_i))_{i=1,\ldots,n} \in \left(\FF_2\right)^{2n}$. We denote by $q_t$ the quadratic form (\ref{eq:quad_form_over_F2}).

\begin{theorem}
\label{thm:classification_of_non_degenerate_quadratic_binary_forms}
Let $V = (\FF_2)^{2n}$ with $n \geq 1$. There are two equivalence classes of non-degenerate quadratic forms over $V$. They are identified by their Arf invariant $\Arf(q) \in \FF_2$ which is defined by
\[
\#\{v \in V; q(v)=0\} - \#\{v \in V; q(v)=1\} = (-1)^{\Arf(q)}\ 2^{n-1}.
\]
The Arf invariant of the form $q_t$ defined in (\ref{eq:quad_form_over_F2}) is the number of indices $i \in \{1,\ldots,n\}$ such that $(a_i,b_i) = (1,1)$ modulo $2$.
\end{theorem}

\begin{proof}
The proof follows from the cases of $n=1$ and $n=2$. For $n=1$, the form $x^2 + xy + y^2$ is invariant under $Sp(B_0) = GL_2(\FF_2)$ whereas the three other forms $xy$, $x^2 + xy$ and $xy + y^2$ are linearly equivalent. We denote $U_0 = \{(0,0),(0,1),(1,0)\}$ and $U_1 = \{(1,1)\}$ and consider the case of $n=2$. The case $n=1$ implies that the forms $q_t$ with $t \in U_0 \times U_0$ are equivalent and using symmetries of coordinates the forms $q_t$ with $t \in U_0 \times U_1 \cup U_1 \times U_0$ are equivalent. There is a linear transformation that maps $q_{(1,1,1,1)}$ to $q_{(0,0,0,0)}$, namely
\[
q_{(1,1,1,1)}(x_1+x_2+y_2,y_1+x_1,x_2+x_1+y_1,y_2+x_2) = q_{(0,0,0,0)}(x_1,y_1,x_2,y_2).
\]
Hence there are at most two equivalence classes. The fact that we have at least two classes follows from the formula relating the Arf invariant to the number of solutions of $q(v)=1$. The general case follows by recurrence.
\end{proof}
The formula in Theorem~\ref{thm:classification_of_non_degenerate_quadratic_binary_forms} states that the Arf invariant of a quadratic form $q$ is the majority value assumed by $q$ on $V$ among $0$ and $1$. We now states a theorem about the classifcation theorem of all quadratic forms.
\begin{theorem}
\label{thm:classification_of_all_quadratic_binary_forms}
Let $V = \left(\FF_2\right)^n$ with $n \geq 2$. There are three linear equivalence classes of quadratic forms on $V$ of rank $2g$ with $0 < g < n$:
\begin{itemize}
\item $\{q;\ q|_{\ker(B_q)} \not= 0\}$,
\item $\{q;\ q|_{\ker(B_q)} = 0$ and $Arf(\overline{q}) = 0\}$ where $\overline{q}$ is $q$ on the quotient $V/\ker(B_q)$,
\item $\{q;\ q|_{\ker(B_q)} = 0$ and $Arf(\overline{q}) = 1\}$.
\end{itemize}
\end{theorem}

Now, we define the spin parity of a permutation. Let $\pi = (\pi_t,\pi_b)$ be a labeled permutation on the alphabet $\mathcal{A}$ with $n$ elements. Let $V := \left(\FF_2\right)^\mathcal{A}$ and $e_\alpha$ be the elementary vector for which the only non zero coordinate is in position $\alpha$. The \emph{intersection form} of $\pi$ is the bilinear form $\Omega = \Omega_\pi$ on $V$ defined by
\[
\Omega_{\alpha,\beta} = \Omega(e_\alpha,e_\beta) = \left\{
\begin{array}{ll}
1 & \text{if $(\pi_t(\alpha) - \pi_t(\beta))(\pi_b(\alpha) - \pi_b(\beta)) < 0$}, \\
0 & \text{else}.
\end{array}
\right.
\]
The matrix $\Omega$ corresponds to crossings: the entry $(\alpha,\beta)$ of the matrix is $1$ if and only if the order of $(\pi_t(\alpha),\pi_t(\beta)$ is the opposite of $(\pi_b(\alpha),\pi_b(\beta))$.

\begin{figure}
\begin{center}
\picinput{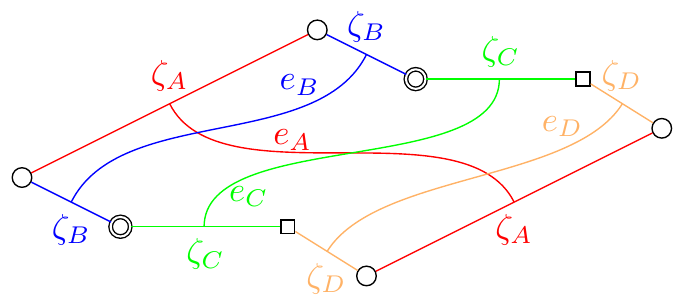}
\end{center}
\caption{Canonic basis of $H_1(S,\Sigma;\ZZ/2)$ and $H_1(S-\Sigma;\ZZ/2)$ of a suspension of $\pi = \binom{A B C D}{B C D A}$.}
\label{fig:homology_of_suspension}
\end{figure}

Let $S=S(\pi,\zeta)$ be a suspension of $\pi$. The sides $(\zeta_\alpha)_{\alpha \in \mathcal{A}}$ of $S$ form a basis of the relative homology $H_1(S,\Sigma; \ZZ/2)$. The elements $(e_\alpha)$ can be considered as its dual basis in $H_1(S-\Sigma;\ZZ/2)$ (see Figure~\ref{fig:homology_of_suspension}).
The intersection form on $H_1(S;\ZZ/2)$ is well defined on $H_1(S-\Sigma;\ZZ/2)$ by composition of the natural morphism $H_1(S-\Sigma;\ZZ/2) \rightarrow H_1(S;\ZZ/2)$ obtained from the inclusion $S - \Sigma \rightarrow S$. The matrix $\Omega$ corresponds to the the intersection matrix of the vectors $(e_\alpha)_{\alpha \in \mathcal{A}}$ viewed as elements of $H_1(S-\Sigma; \ZZ/2)$. In particular the rank of $\Omega_\pi$ is $2g$ where $g$ is the genus of the suspension.

We remark that $\Omega$ only depends on the topological structure of $S(\pi,\zeta)$ and not on the flat metric. Now, we define a quadratic form $q_\pi$. For any closed curve $\gamma: [0,1] \rightarrow S$ there is an associated winding number (relative to the flat metric) which is an integer multiple of $2\pi$. We denote by $w(\gamma)$ this integer modulo $2$ and extends it by linearity to $H_1(S-\Sigma;\ZZ/2)$. We may notice that any linear form on $\FF_2$ can be canonically transformed into a totally degenerate quadratic form without changing its values as $0^2=0$ and $1^2=1$. The quadratic form $q_\pi$ on $H_1(S-\Sigma;\ZZ/2)$ is
\[
q(x) = w(x) + \#(\text{components of x}) + \# (\text{self intersections of x}).
\]

\begin{proposition}
Let $\pi$ be a permutation. The quadratic form $q_\pi$ is such that the restriction to $\ker(B_{q_\pi})$ is null if and only if the profile of $\pi$ as only odd parts.
\end{proposition}

\begin{proof}
Let $q$ be the quadratic form of $\pi$ and $B$ its associated bilinear form. The vector space $\ker(B)$ is generated by small loops around the singularities (each loop around a singularity is non trivial in $H_1(S - \Sigma; \ZZ/2\ZZ)$ and becomes trivial in $H_1(S;\ZZ/2\ZZ)$). Let $\gamma$ be a simple curve around a singularity of angle $k 2\pi$. The winding number of $\gamma$ is $w(\gamma) = k$ and hence $q_\pi(\gamma) = k+1$.
\end{proof}

\begin{definition}
Let $\pi$ be a permutation such that its profile has only odd parts. The \emph{spin parity} of $\pi$ is the Arf invariant of the quadratic form $q_\pi$.
\end{definition}

As an example the permutations 
\[
\pi_0 = \binom{0\ 1\ 2\ 3\ 4\ 5\ 6\ 7}{7\ 4\ 3\ 2\ 1\ 6\ 5\ 0}
\quad \text{and} \quad
\pi_1 = \binom{0\ 1\ 2\ 3\ 4\ 5\ 6\ 7}{7\ 2\ 1\ 4\ 3\ 6\ 5\ 0}
\]
have both profiles $(7)$ but the spin parity are, respectively, $0$ and $1$. The permutations $\pi_0$ and $\pi_1$ hence belong to two different Rauzy classes. This fact can be checked by explicit computation of Rauzy classes but is fastidious as the cardinality of Rauzy classes are respectively $5209$ and $2327$.

\subsubsection{Hyperellipticity}
\label{subsubsection:hyperellipticity}
A translation surface $S$ is \emph{hyperelliptic} if there exists a morphism of degree two from $S$ to the Riemann sphere $\PP^1 \CC$ such that the flat structure of $S$ comes from a quadratic differential on $\PP^1 \CC$.

\begin{proposition}[\cite{KontsevichZorich03}]
\label{prop:hyperelliptic_components}
In the strata $\Omega \mathcal{M}(2g-1)$ (resp. $\Omega \mathcal{M}(g,g)$) there exists a connected component $\Omega \mathcal{M}^{hyp}(2g-1)$ (resp. $\Omega \mathcal{M}^{hyp}(g,g)$) such that each surface in the component is hyperelliptic. These two families are the only connected components of strata without marked point with this property.
\end{proposition}
For strata $\Omega \mathcal{M}(2g-1,1^k)$ and $\Omega \mathcal{M}(g,g,1^k)$ which contain marked points, there is also a connected components which comes from the hyperelliptic ones in $\Omega \mathcal{M}(2g-1)$ and $\Omega \mathcal{M}(g,g)$. We will call them \emph{hyperelliptic} as well.

\begin{proposition}[\cite{KontsevichZorich03}] \label{prop:hyperelliptic_components_and_hyperelliptic_permutations}
Let $\pi$ be an irreducible permutations with profile $(2g-1)$ or $(g,g)$ and $S$ a suspension of $\pi$. Then $S$ is in an hyperelliptic component of $\Omega \mathcal{M}$ defined in Proposition~\ref{prop:hyperelliptic_components} if and only if $\pi$ is in the Rauzy class of a symmetric permutation
\[
\pi_n^{sym} = \binom{1\ 2\ \ldots\ n}{n\ n-1\ \ldots \ 1}.
\]
\end{proposition}

\subsection{Definition of Rauzy classes in terms of invariants}
\label{subsection:Rauzy_classes_in_terms_of_invariants}
As we have seen in Proposition~\ref{prop:associated_connected_component_of_a_rauzy_class}, we can associate to each Rauzy class and each extended Rauzy class a connected component of a stratum $\Omega \mathcal{M}(p_\pi)$. In this section we recall the results of \cite{Veech82} and \cite{Boissy09} which prove how this association can be turned into a one to one correspondance. Next, we explain the classification of connected components of strata of \cite{KontsevichZorich03} and deduce a classification of Rauzy classes.

\subsubsection{Connected components of moduli space and Rauzy classes}
\label{subsubsection:connected_components_and_Rauzy_classes}
In order to get a correspondance between Rauzy classes and connected components of moduli space of translation surfaces, we need to encode a combinatorial data which corresponds to the fact that the Rauzy induction fixes the left endpoint of the interval. Let $\Omega \mathcal{M}(p)$ be a stratum and $m_l \in p$. Let $p' = p \backslash \{m_l\}$. We denote by $\Omega \mathcal{M}(m_l; p')$ the moduli space of translation surfaces $\Omega \mathcal{M}(p)$ with a choosen singularity of degree $m_l$.

If $\pi$ is a permutation, we denote by $m_l(\pi)$ the angle of the singularity on the left of $\pi$. It corresponds to the length of the cycle of the interval diagram which contains the element $\outgos{\pi_b^{-1}(1),\pi_t^{-1}(1)}$ (see Section~\ref{subsubsection:profile}). To an irreducible permutation we associate a connected component with a choosen singularity of degree $m_l$.

\begin{theorem}[\cite{Veech82},\cite{Boissy09}]
\label{thm:extended_rauzy_classes_and_connected_components}
\label{thm:rauzy_classes_and_connected_components}
The association $\pi \mapsto \Omega \mathcal{M}(p_\pi)$ induces a bijection between extended Rauzy classes of irreducible permutations and connected components of strata of moduli spaces $\Omega \mathcal{M}(p)$.

\noindent The association $\pi \mapsto \Omega \mathcal{M}(m_l(\pi); p'_\pi)$ induces a bijection between Rauzy classes and connected components of strata of moduli spaces with a chosen fixed degree.
\end{theorem}

\begin{corollary}[\cite{Boissy09}]
Let $\mathcal{R}$ be an extended Rauzy class associated to a connected component $\mathcal{C}$ of a stratum $\Omega \mathcal{M}(p)$. Then $\mathcal{R}$ is the union of $r$ Rauzy classes where $r$ is the number of distinct elements of $p$.
\end{corollary}
If $\mathcal{R}$ is an extended Rauzy class, we denote by $\mathcal{R}(m_l)$ the Rauzy class which consist of permutations for which $m_l(\pi) = m_l$. Note that $r$ is not the number of singularities, we have $r=1$ for any connected component of $\Omega \mathcal{M}(2,2,2,2)$.

There is a map from a component with chosen fixed degree to the one without: $\Omega \mathcal{M}(m_l;p') \rightarrow \Omega \mathcal{M}(p)$. At the level of Rauzy classes this corresponds to a disjoint union: the extended Rauzy class corresponding to a permutation $\pi$ is the union of the Rauzy classes associated to the possible degrees associated to the left endpoint. As an example there is one extended Rauzy $\mathcal{R}$ class with $2638248$ elements associated to the connected stratum $\Omega \mathcal{M}(4,3,2,1)$ which is the union of four Rauzy classes $\mathcal{R}(4)$, $\mathcal{R}(3)$, $\mathcal{R}(2)$ and $\mathcal{R}(1)$ with respectively $1060774$, $792066$, $538494$ and $246914$ elements.

The labeled Rauzy classes also have a geometric interpretation in terms of moduli space of translation surfaces. If $\pi = (\pi_t,\pi_b)$ is a labeled permutation, then the permutation $\outgos{\sigma_\pi}$ deduced from the Rauzy diagram (see Section~\ref{subsubsection:profile}) is invariant under Rauzy induction which implies a bijection as Theorem~\ref{thm:rauzy_classes_and_connected_components} between labeled Rauzy classes and a moduli space of translation surfaces with combinatorial data. In this case, the combinatorial data consist in a label for each horizontal outgoing separatrices of the surface. The classification of connected component of this moduli space is done in \cite{Boissy10}. In particular, he establishes a formula that relates the cardinality of a labeled Rauzy class of a permutation $(\pi_t,\pi_b)$ and the cardinality of the reduced Rauzy class of the associated reduced permutation $\pi_b \circ \pi_t^{-1}$. But we emphasize that there is no known relation between labeled extended  Rauzy classes and moduli space of translation surfaces.

\subsubsection{Kontsevich-Zorich classification of connected components}
\label{subsubsection:classification_of_connected_components}
The strata of moduli spaces of translation surfaces $\Omega \mathcal{M}(p)$ are not connected in general. The three invariants above (profile, spin, and hyperellipticity) as proved in \cite{KontsevichZorich03} are enough to give a complete classification.
\begin{theorem}[\cite{KontsevichZorich03}] \label{thm:classification_of_connected_components}
The connected components of a stratum with marked points $\Omega \mathcal{M}(n_1,n_2,\ldots,n_k,1^l)$ are in bijection with connected components of the stratum $\Omega \mathcal{M}(n_1,n_2,\ldots,n_k)$.

The classification of connected components of stratum whose profile does not contains any $1$ are given by the classification below. For genus $g \geq 4$ we have
\begin{itemize}
\item The strata $\Omega \mathcal{M}(2g-1)$ and $\Omega \mathcal{M}(g,g)$ with $g$ odd have three components: a hyperelliptic component associated to the symmetric permutations on respectively $2g$ and $2g+1$ letters. A component with odd spin parity and a component with even spin parity.
\item The other strata with only odd parts $\mathcal{H}(2m_1+1, 2m_2+1, \ldots, 2m_n+1)$ have two connected components which are distinguished by their spin parities.
\item $\Omega \mathcal{M}(g,g)$ for $g$ even has two components: one hyperelliptic and an other one (called the non-hyperelliptic component).
\item Any other stratum is connected.
\end{itemize}
For small genera, the preceding classification holds but there are empty components:
\begin{itemize}
\item genus $1$ and $2$: the strata $\Omega \mathcal{M}(1)$, $\Omega \mathcal{M}(3)$ and $\Omega \mathcal{M}(2,2)$ are non empty and connected.
\item genus $3$: $\Omega \mathcal{M}(5)$ and $\Omega \mathcal{M}(3,3)$ have two connected components one hyperelliptic and one odd. The other strata of $\Omega \mathcal{M}_3$ are connected.
\end{itemize}
\end{theorem}

By the above theorem, Theorem~\ref{thm:extended_rauzy_classes_and_connected_components} and Theorem~\ref{thm:rauzy_classes_and_connected_components} we obtain the following classification of Rauzy classes.

\begin{theorem}
\label{thm:classification_of_Rauzy_classes}
Let $p = (n_1, \ldots, n_k)$ be an integer partition such that $s(p) + l(p) \equiv 0 \mod 2$. Then the set of permutations $\pi$ with profile $p$ is the union of $1$, $2$ or $3$ extended Rauzy classes depending on the number of connected components of $\Omega \mathcal{M}(p)$ given by Theorem~\ref{thm:classification_of_connected_components}. Each extended Rauzy class is the union of $r$ Rauzy classes where $r$ is the number of distinct part in $p$.
\end{theorem}

Recall from the introduction that if $p$ be a partition such that $s(p) + l(p) \equiv 0 \mod 2$ we denote by $\gamma^{irr}(p)$ the number of irreducible permutations with profile $p$. Moreover, if $p$ has only odd terms we denote $\delta^{irr}(p) = \gamma^{irr}_1(p) - \gamma^{irr}_0(p)$ where $\gamma^{irr}_s(p)$ is the number of irreducible permutations with profile $p$ ans spin paruty $s$.

The below corollary is a direct consequence of Theorem~\ref{thm:classification_of_Rauzy_classes}.
\begin{corollary} \label{corollary:formula_for_the_cardinality_of_a_Rauzy_class}
Let $p$ be an integer partition such that $s(p)+l(p) \equiv 0 \mod 2$ and $\Omega \mathcal{M}(p)$ the stratum of the moduli space of translation surfaces with profile $p$.

If $\Omega \mathcal{M}(p)$ is connected then the only Rauzy class $\mathcal{R}$ which consists of irreducible permutations with profile $p$ satisfies $|\mathcal{R}| = \gamma^{irr}(p)$.

If $\Omega \mathcal{M}(p)$ is a union of an odd and an even component then there are two Rauzy classes $\mathcal{R}^{odd}$ and $\mathcal{R}^{even}$ with profile $p$ which satisfy $|\mathcal{R}^{odd}| = \frac{\gamma^{irr}(p) + \delta^{irr}(p)}{2}$ and $|\mathcal{R}^{even}| = \frac{\gamma^{irr}(p) - \delta^{irr}(p)}{2}$.

If $p=(g,g,1^k)$ with $g$ even, then there are two Rauzy classes $\mathcal{R}^{hyp}$ and $\mathcal{R}^{nonhyp}$ with profile $p$ which satisfy $|\mathcal{R}^{nonhyp}|  = \gamma^{irr}(p) - |\mathcal{R}_{hyp}|$.

If $p=(2g-1,1^k)$ or $p=(g,g,1^k)$ with $g$ odd, then there are three Rauzy classes $\mathcal{R}^{hyp}$, $\mathcal{R}^{odd}$ and $\mathcal{R}^{even}$ associated respectively to the hyperelliptic, the odd spin and even spin components of $\Omega \mathcal{M}(2g-1,1^k)$ (resp. $\Omega \mathcal{M}(g,g,1^k)$ with $g$ odd). Then, if $g \equiv 1,2 \mod 4$ then
\[
|\mathcal{R}^{odd}| = \frac{\gamma^{irr}(p)  + \delta^{irr}(p)}{2}\ -\ |\mathcal{R}_{hyp}|
\quad \text{and} \quad
|\mathcal{R}^{even}| = \frac{\gamma^{irr}(p) - \delta^{irr}(p)}{2}.
\]
And if $g \equiv 0,3 \mod 4$ then
\[
|\mathcal{R}^{odd}| = \frac{\gamma^{irr}(p)  + \delta^{irr}(p)}{2}
\quad \text{and} \quad
|\mathcal{R}^{even}| = \frac{\gamma^{irr}(p) - \delta^{irr}(p)}{2}\ -\ |\mathcal{R}^{hyp}|.
\]
\end{corollary}

As an example, the $461$ irreducible permutations on six letters is the union of seven Rauzy classes (respectively five extended Rauzy classes) as below:
\begin{itemize}
\item two Rauzy classes (two extended) associated to $\Omega \mathcal{M}_3(5) = \Omega \mathcal{M}_3^{hyp}(5) \cup \Omega \mathcal{M}_3^{odd}(5)$ with respectively $31$ and $134$ permutations,
\item two Rauzy classes (one extended) associated to $\Omega \mathcal{M}_2(3; 1,1)$ and $\Omega \mathcal{M}_2(1; 3,1)$ with respectively $105$ and $66$ permutations,
\item two Rauzy classes (one extended) associated to $\Omega \mathcal{M}_2(2; 2,1)$ and $\Omega \mathcal{M}_2(1; 2,2)$ with respectively $90$ and $20$ permutations,
\item one Rauzy class (one extended) associated to $\Omega \mathcal{M}_1(1,1,1,1,1)$ with $15$ elements.
\end{itemize}
Corollary~\ref{corollary:formula_for_the_cardinality_of_a_Rauzy_class} can be formulated as well for Rauzy classes introducing natural notations $\gamma^{irr}(m,p')$ and $\delta^{irr}(m,p')$.

\section{Enumerating labeled standard permutations}
\label{section:counting_standard_permutations}
In this section we are interested in the number of standard permutations (Definition~\ref{def:standard_permutation}) in any Rauzy class (Definition~\ref{def:rauzy_class_combinatorial}) which is the starting point to enumerate the whole class. Recall that the conjugacy classes of $S_n$ are in bijection with integer partition of $n$. To a permutation $\sigma$ we associate the integer partition $(n_1,\ldots,n_k)$ whose parts are the lengths of the cycles in the disjoint cycle decomposition. As the bijection is canonic we identify conjugacy classes of $S_n$ and integer partition of $n$. 

Let $p$ be an integer partition and $\sigma \in S_n$ a permutation whose conjugacy class is $p$. We establish in Proposition~\ref{prop:bijection_standard_permutations_constellations} a bijection between the solutions $(\tau_t,\tau_b)$ of the equation
\begin{equation}
\label{eq:cyclic_constellation_equation}
\sigma = \tau_t \ \tau_b^{-1} \qquad \text{where $\tau_t$, $\tau_b$ are $n$-cycles of $S_n$},
\end{equation}
and the labeled permutations $(\pi_t,\pi_b)$ with profile $p$ and fixed labels on outgoing separatrices (see Section~\ref{subsubsection:profile}). We denote by $c(p)$ the number of solutions of~(\ref{eq:cyclic_constellation_equation}) as it does not depend on the choice of $\sigma$ with conjugacy class $p$. We remark that when $p$ statisfies $l(p)+s(p) \not\equiv 0 \mod 2$ then there is no labeled permutation with profile $p$ (because $s(p)+l(p) = 2g-2$ where $g$ is the genus of a suspension of $\pi$, see (\ref{eq:genus_and_conic_points}) in Section~\ref{subsection:from_permutations_to_surfaces}). On the other hand, the signature $\tau$ of a permutation with conjugacy class $p$ is $\varepsilon(\tau) = (-1)^{s(p)+l(p)}$. Hence, if there is a solution $(\tau_t,\tau_b) \in S_n \times S_n$ of~(\ref{eq:cyclic_constellation_equation}) the signature of $\sigma$ is necessarily $1$.

If $p$ has only odd parts (in which case the condition $s(p) + l(p) \equiv 0 \mod 2$ is automatic), we denote by $c_1(p)$ (resp. $c_0(p)$) the number of labeled permutations with spin parity $1$ (resp. $0$) and set $d(p) = c_1(p) - c_0(p)$. Using geometrical analysis, we prove recurrence formulas for $c$ and $d$ (Theorems~\ref{thm:std_perm_strata_recurrence} and \ref{thm:std_perm_spin_recurrence}) and then provide explicit formulas for both (Theorems~\ref{thm:Boccara_formula_for_c} and \ref{thm:formula_for_d}).

\subsection{Standard permutations and equations in the symmetric group}
\label{subsection:standard_permutation_to_cyclic_constellations}
The particular form of a standard permutation allows the construction of a surface which is no more built from a polygon but from a cylinder. We explain this construction which can be found in \cite{KontsevichZorich03}, \cite{Zorich08} and \cite{Lanneau08}. Instead of considering a standard permutation $\pi$ as a double ordering $\pi_t,\pi_b$ of the alphabet $\mathcal{A}$, we describe it as a triple of permutations $(\tau_t, \tau_b, \sigma) \in S_{\outgos{\AAA}_\pi} \times S_{\outgos{\AAA}_\pi} \times S_{\outgos{\AAA}_\pi}$ with the following properties
\begin{itemize}
\item $\tau_t$ and $\tau_b$ are $n-1$ cycles,
\item $\sigma = \tau_t\ \tau_b^{-1}$ is the permutation $\outgos{\sigma_\pi}$,
\end{itemize}
where the notation $\outgos{\AAA_\pi}$ and $\outgos{\sigma_\pi}$ were defined in Section~\ref{subsubsection:profile}.

Given $(\tau_t, \tau_b, \sigma) \in S_n \times S_n \times S_n$ with $\tau_t \tau_b^{-1} = \sigma$ we develop the method of \cite{Boccara80} which consists in defining another triple $(\tau_t', \tau_b', \sigma') \in S_{n-1} \times S_{n-1} \times S_{n-1}$ in order to relate the solutions of (\ref{eq:cyclic_constellation_equation}) in $S_n$ to the ones on $S_{n-1}$.

\subsubsection{Cylindric suspension and equation $\sigma = \tau_t\ \tau_b^{-1}$}
\label{subsubsection:standard_permutations_cyclic_constellations}
\begin{figure}[ht!]
\centering
\picinput{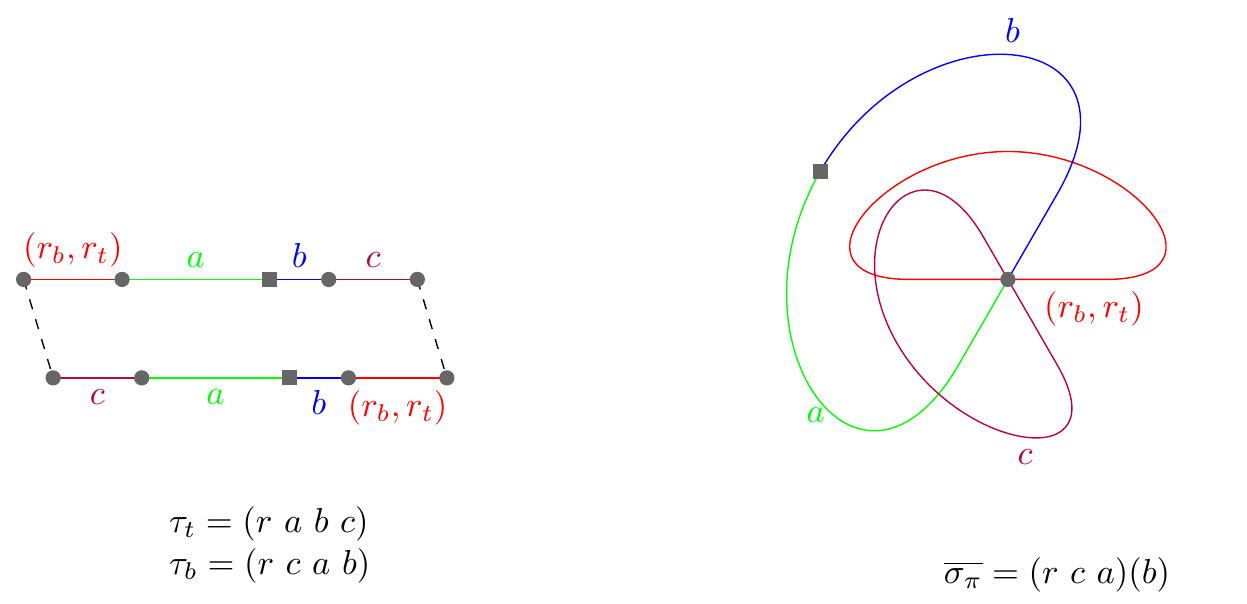}
\caption{Cylindric suspension of $\pi = \binom{r_t\ a\ b\ c\ r_b}{r_b\ c\ a\ b\ r_t}$ and its interval diagram.}
\label{fig:cylindric_suspension}
\end{figure}

Let $\pi = (\pi_t,\pi_b)$ be a labeled standard permutation on the alphabet $\mathcal{A}$ of cardinality $n+1$. Let $r_t = \pi_t^{-1}(1) = \pi_b^{-1}(n)$ and $r_b =\pi_t^{-1}(n) = \pi_b^{-1}(1)$. Let $\zeta \in \CC^{\mathcal{A}}$ be such that
\begin{itemize}
\item $Im(\zeta_{r_b}) < 0$ and $Re(\zeta_{r_b}) > 0$,
\item for all $\alpha \not= r_b$, $Im(\zeta_\alpha) = 0$ and $Re(\zeta_\alpha) > 0$.
\end{itemize}
Therefore, the vector $\zeta$ is not a suspension data as in Section~\ref{subsubsection:suspensions}. However, using the same construction with broken lines $L_t$ and $L_b$, we get a surface which we call a \emph{cylindric suspension} of $\pi$ (see Figure~\ref{fig:cylindric_suspension}). If we glue together the vertical associated to $r_b$ on $L_t$ to the one on $L_b$ we obtain an horizontal cylinder. Its boundary consists of two circles cut in $n$ intervals.

There is an arbitrary choice between $r_t$ and $r_b$ as vertical edge. To take care of this flexibility, we label the top and bottom circles with the alphabet 
\[
\mathcal{A}_\pi = \{(r_b,r_t)\} \cup \{\alpha \in \mathcal{A}: \alpha \not= r_b \ \text{and} \ \alpha \not=r_t\},
\]
instead of $\mathcal{A} \backslash \{r_b\}$ (see Figure~\ref{fig:cylindric_suspension}). Remark that the labelization of the two circles coincide with the interval diagram defined in Section~\ref{subsubsection:profile}. Recall that the interval diagram $\sigma_\pi$ of $\pi$ is a permutation defined on the alphabet $\outgos{\mathcal{A}_\pi} \cup \incoms{\mathcal{A}_\pi}$ which consists in two copies of $\mathcal{A}_\pi$ above. The interval diagram $\sigma_\pi$ exchanges $\outgos{\mathcal{A}_\pi}$ and $\incoms{\mathcal{A}_\pi}$. The square of $\sigma_\pi$ decomposes as a product of two permutations $\outgos{\sigma_\pi}$ and $\incoms{\sigma_\pi}$  on respectively $\outgos{\mathcal{A}_\pi}$ and $\incoms{\mathcal{A}_\pi}$.

\begin{proposition}
\label{prop:bijection_standard_permutations_constellations}
Let $\mathcal{A}$ be a finite alphabet and $r_t,r_b$ two distinct elements of $\mathcal{A}$. Set $\mathcal{A}' = \{(r_b,r_t)\} \cup \mathcal{A} \backslash \{r_b,r_t\}$. Let $\sigma \in S_{\mathcal{A}'}$, then there is a bijection between the set of labeled standard permutations $\pi$ on the alphabet $\mathcal{A}$ such that $\outgos{\sigma_\pi} = \sigma$ and the set of solutions $(\tau_t,\tau_b) \in S_{\mathcal{A}'} \times S_{\mathcal{A}'}$ such that $\tau_t \tau_b^{-1} = \sigma$.
\end{proposition}

\begin{proof}
Let $n$ be the cardinality of $\mathcal{A}$. The proof follows directly from the definition of the interval diagram (Definition~\ref{def:interval_diagram}). Let $\pi$ be a standard permutation on $\mathcal{A}$. We associate to $\pi$ the two $n$-cycles that consists of the top and bottom lines
\[
\tau_t = ((r_b,r_t)\ \pi_t^{-1}(2)\ \pi_t^{-1}(3)\ \ldots\ \pi_t^{-1}(n-1))
\quad \text{and} \quad
\tau_b = ((r_b,r_t)\ \pi_b^{-1}(2)\ \pi_b^{-1}(3)\ \ldots\ \pi_b^{-1}(n-1)).
\]
The fact that $\tau_t$, $\tau_b$ and $\sigma$ satisfies Equation~(\ref{eq:cyclic_constellation_equation}) can be resumed in the following picture (Figure~\ref{fig:cyclic_constellation_relation}) which represents a vertex of a suspension $S(\pi,\zeta)$ of $\pi$ together as the action of $\tau_t$, $\tau_b$ and $\sigma$ as permutation.
\begin{figure}[!ht]
\centering
\picinput{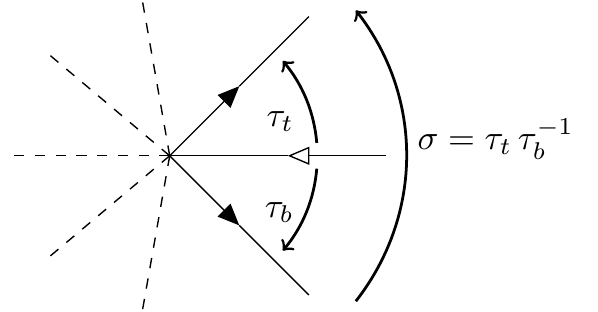}
\caption{The relation $\sigma = \tau_t\, \tau_b^{-1}$ on the level of the interval diagram of $\pi$.}
\label{fig:cyclic_constellation_relation}
\end{figure}
\end{proof}
An example is shown in Figure~\ref{fig:cylindric_suspension}.

Counting labeled stantard permutations is now expressed in a group theoritical way. Let $X$ $Y$, $Z$ be three conjugacy classes of a finite group $G$, we want to count the number of solutions of an equation $x\, y\, z = 1$ where $x \in X$, $y \in Y$ and $z \in Z$. This problem is known to be equivalent to a formula involving characters called the Frobenius formula.
\begin{proposition}[Frobenius formula]
Let $G$, $X$, $Y$ and $Z$ as above. Let $N_{X,Y,Z}$ be the number of triples $(x,y,z) \in X \times Y \times Z$ such that $x\ y\ z = 1$. Then
\[
N_{X,Y,Z} = \frac{|X|\ |Y|\ |Z|}{|G|} \sum_{\chi \in \hat{G}} \frac{\chi(X)\ \chi(Y) \chi(Z)}{\chi(1)},
\]
where $\hat{G}$ denotes the set of irreducible characters of $G$.
\end{proposition}
The proof of Frobenius formula can be found for example in Section 7.2. of \cite{Serre92}. For the numbers $c(p)$ we deduce from Frobenius formula the following expression
\begin{equation}
c(p) = (n-1)! \sum_{\chi \in \hat{S_n}} \frac{\chi{(n)}^2 \chi(p)}{\chi(1)}.
\label{eq:equation_and_characters}
\end{equation}
It is a hard to pass from expression~(\ref{eq:equation_and_characters}) which involves characters to a formula which involves numbers. The recursive construction we adopt does not use Frobenius formula. However there are some works, for example \cite{GoupilSchaeffer98} (see Theorem~\ref{thm:GoupilSchaeffer_formula_for_c}), that from Frobenius formula obtain formulas for the value of $c(p)$. The conjugacy class of $\sigma$ encodes the stratum associated to the suspension of $\pi$. For the numbers $d(p)= c_1(p) - c_0(p)$, there is still an approach using Group Theory. The spin parity can be viewed as a refinement of the signature of a permutation in the Sergeev group \cite{EskinOkounkovPandharipande08}.

\subsubsection{Recursive construction}
\label{subsection:recursive_construction}
In order to obtain formulas for the numbers $c(p)$ and $d(p)$ we follow an approach of \cite{Boccara80}. Let $\AAA$ be an alphabet of size $n$ and $\sigma \in S_\AAA$ a permutation. Let $(\tau_t, \tau_b, \sigma)$ be a solution of Equation~(\ref{eq:cyclic_constellation_equation}) and $x \in \AAA$. Starting from a triple $(\tau_t,\tau_b,\sigma) \in S_\AAA$ of equation~(\ref{eq:cyclic_constellation_equation}), we choose a letter $x \in \AAA$, then we remove $x$ in both cycles $\tau_t$ and $\tau_b$ and get two $(n-1)$-cycles $\tau_t'$ and $\tau_b'$ on $\AAA' = \AAA \backslash \{x\}$. Set $\sigma' = \tau_t'\ {\tau_b'}^{-1}$, we want to know the relation between $\sigma$ and $\sigma'$.

The $(n-1)$-cycles $\tau_t'$ and $\tau_b'$ obtained are formally given by
\[
\tau'(y) = \left\{ \begin{array}{ll}
  \tau(x) & \text{if $y = \tau^{-1}(x)$}, \\
  \tau(y) & \text{otherwise}. \\
\end{array} \right.
\qquad \text{where $\tau$ equals $\tau_t$ or $\tau_b$.}
\]
The operation $\tau \mapsto \tau'$ can be obtained as a multiplication by a transposition, where we consider $\tau'$ as a permutation on $\mathcal{A}$ which fixes $x$. More precisely
\begin{equation}
\tau_t' = (x \ \tau_t(x)) \ \tau_t \quad \text{and} \quad {\tau_b'}^{-1} = \tau_b^{-1} \ (x\ \tau_b(x)).
\label{eq:RemovingAInTau}
\end{equation}
To $\tau_t'$ and $\tau_b'$ which are $(n-1)$-cycles on $\mathcal{A}'$ we associate the permutation $\sigma'$ by the formula $\sigma' = \tau_t' \ \tau_b'^{-1}$. Using formulas~(\ref{eq:RemovingAInTau}) we write $\sigma'$ as a product involving $(\tau_t, \tau_b, \sigma)$ and the letter $x$
\begin{equation}
\sigma' = \tau_t'\,  {\tau_b'}^{-1} = (x\ \tau_t(x))\, \sigma \, (x\ \tau_b(x)).
\label{eq:RemovingAInSigma}
\end{equation}
The conjugacy class of $\sigma'$ depends only of the positions of $x$ and $\tau_t(x)$ in the cycle decomposition of $\sigma$. If $p = (n_1,\ldots,n_k)$ and $p'=(n'_1,\ldots,n'_{k'})$ are integer partitions we denote $p \uplus p' = (n_1,\ldots,n_k,n'_1,\ldots n'_{k'})$ their disjoint union. If $m$ is an integer we write $m \in p$ if $m$ is a part of $p$ and if $q$ is an integer partition we write $q \subset p$ if there exists $p'$ such that $p = q \uplus p'$. In which case $p'$ is denoted $p \backslash q$.
\begin{definition}[\cite{Boccara80}]
\label{def:splitting_and_collapsing_of_partition}
Let $p$ be an integer partition of $n$. Let $m \in p$ and $a \in \{1,2,\ldots,m-2\}$, we call the \emph{splitting} of $m$ in $p$ by $a$ the integer partition
\[ p_{m|a} = (a,\, m - a - 1)\ \uplus\ p \backslash \{m\}.\]
Let $(m_l,m_r) \subset p$ we call the \emph{collapsing} of $m_l$ and $m_r$ in $p$ the integer partition
\[ p_{m_l \odot m_r} = (m_l+m_r-1)\ \uplus\ p \backslash (m_l,m_r). \]
\end{definition}
Remark that if $p$ is a partition of $n$ then both $p_{m|a}$ and $p_{m_l \odot m_r}$ are partitions of $n-1$.
\begin{proposition}[\cite{Boccara80}]
\label{prop:collapse_or_split_combinatoric}
Let $(\tau_t,\,\tau_b,\,\sigma)$, $p$ the conjugacy class of $\sigma$, $x \in \mathcal{A}$ and $(\tau_t',\,\tau_b',\,\sigma')$ be as above. If $x$ and $\tau_t(x)$ are in the same cycle of $\sigma$ with length $m$, then the conjugacy class of $\sigma'$ is $p_{m|a}$ where $a$ is the smallest number such that $\sigma^{a} (\tau_t(x)) = x$. If $x$ and $\tau_t(x)$ are in different cycles of $\sigma$ of length, respectively, $m_l$ and $m_r$, then the profile of $\sigma'$ is $p_{m_l \odot m_r}$. 
\end{proposition}
We remark that $x$, $\tau_t(x)$ and $\tau_b(x)$ belong to the same cycle $c$ of $\sigma$. More precisely, $\tau_b(x)$ and $\tau_t(x)$ are successive letters in $c$, as by definition $\sigma (\tau_b(x)) = \tau_t(x)$.
\begin{proof}[Proof of \ref{prop:collapse_or_split_combinatoric}]
By~(\ref{eq:RemovingAInTau}) and (\ref{eq:RemovingAInSigma}), the differences between $\sigma$ and $\sigma'$ occur for $\sigma^{-1}(x)$ and $\tau_b(x)$ for which we have
\begin{equation}
\sigma'(\tau_b(x)) =  \sigma(x) \qquad \text{and} \qquad \sigma'(\sigma^{-1}(x)) = \tau_t(x).
\label{eq:expression_for_sigma_sigmap}
\end{equation}
We first prove the first part of the proposition. We assume that $x$ and $\tau_t(x)$ belong to different cycles $c_l$ and $c_r$ of $\sigma$ whose lengths are, respectively, $m_l$ and $m_r$. We write $c_l = \left(x\ A\ \sigma^{-1}(x)\right)$ and $c_r = \left(\tau_t(x)\ B\ \tau_b(x)\right)$ where $A$ and $B$ are two blocks of labels which may be empty. The cycles $c_l$ and $c_r$ collapse in $\sigma'$ in a unique cycle c =$ \left(\tau_t(x)\ B\ \tau_b(x)\ A\ \sigma^{-1}(x)\right)$. Because $x$ is removed the length of $c$ is $m_l+m_r-1$.

Now, consider the second part of the proposition. We assume that $x$ and $\tau_t(x)$ are in the same cycle $c$ of $\sigma$ of length $m$. Because $\sigma(\tau_b(x)) = \tau_t(x)$, the cycle of $\sigma$ containing $x$ writes $c = (\tau_t(x)\ A_t\ \sigma^{-1}(x)\ x\ \sigma(x)\ A_b\ \tau_b(x))$, with $\sigma(x) \not= \tau_b(x)$ and $\sigma^{-1}(x) \not= \tau_t(x)$. As before, $A_t$ and $A_b$ are two blocks which may be empty. Now $\sigma'$ has the same cycle decomposition as $\sigma$ but the cycle containing $x$ splits into two cycles $c_t=(\tau_t(x)\ A_t\ \sigma^{-1}(x))$ and $c_b=(\sigma(x)\ A_b\ \tau_b(x))$. The lengths $a_t$ and $a_b$ of the cycles $c_t$ and $c_b$ can be defined symmetrically by $\sigma^{a_t}(\tau_t(x)) = x$ and $\sigma^{-a_b}(\tau_b(x)) = x$. Therefore, as the label $x$ is removed, those lengths satisfy the expression $a_t + a_b = m - 1$.
\end{proof}
Proposition~\ref{prop:collapse_or_split_combinatoric} is the heart of the recurrence formula for the numbers $c(p)$ (Theorem~\ref{thm:std_perm_strata_recurrence}).
\subsection{Spin parity}
\label{subsubsection:geometric_operation}
Let $(\tau_t,\,\tau_b\,\sigma) \in S_\mathcal{A} \times S_\mathcal{A} \times S_\mathcal{A}$ be a solution of~(\ref{eq:cyclic_constellation_equation}) and $x \in \mathcal{A}$. The suppression of the label $x$ in the cycle decomposition of $\tau_t$ and $\tau_b$ studied in the preceding section leads to a solution $(\tau_t',\,\tau_b',\,\sigma')$ on $\mathcal{A} \backslash \{x\}$. Let $S$ be a cylindric suspension of $(\tau_t, \tau_b, \sigma)$. The geometric operation associated to the suppression of $x$ corresponds to remove a cylinder associated to the edge $\zeta_x$ in $S$ (see Figure~\ref{fig:removing_a_cylinder_in_cylindric_suspension}). The operation leads to a cylindric suspension $S'$ of $(\tau_t',\tau_b',\sigma')$. Proposition~\ref{prop:collapse_or_split_combinatoric} can be interpreted as an answer to the stratum behavior of the operation $(S,\zeta_x) \mapsto S'$ (see Proposition~\ref{prop:collapse_or_split_geometric}). In this section, we analyze the geometric operation and get a relation between the spin parities of $S$ and $S'$.

\subsubsection{Removing a cylinder in a translation surface}
\label{subsubsection:removing_a_cylinder}
\begin{figure}[!ht]
\begin{center}
\picinput{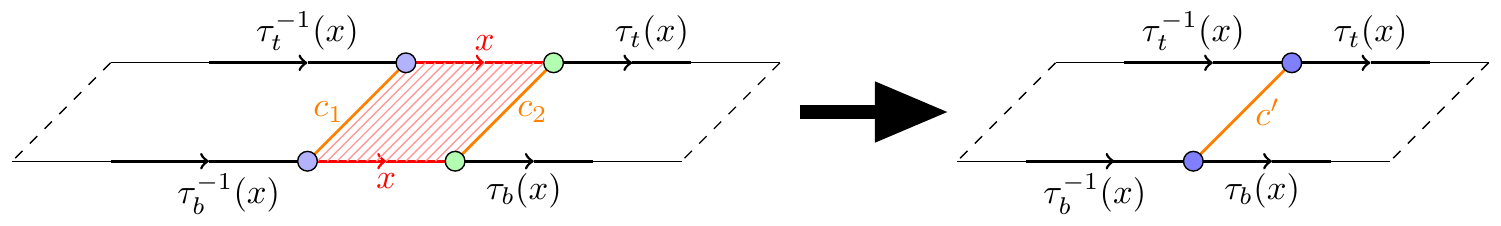}
\end{center}
\caption{Removing the cylinder associated to $\zeta_x$ in a cylindric suspension.}
\label{fig:removing_a_cylinder_in_cylindric_suspension}
\end{figure}

In a cylindric suspension $S$ of a triple $(\tau_t,\tau_b,\sigma) \in S_\AAA \times S_\AAA \times S_\AAA$, a label $x \in \AAA$ corresponds to a horizontal geodesic $\zeta_x$ in $S$ which join two singularities (possibly the same). More generally, let $S$ be a translation surface and $\zeta$ a geodesic segment between two singularities of $S$. We assume that $\zeta$ contains no singularity in its interior. Such a segment is called a \emph{saddle connection}.
\begin{definition}
Let $S$ be a translation surface and $\zeta$ a saddle connection in $S$. A geodesic cylinder which contains $\zeta$ in its interior and each of its boundary circle contains an endpoint of $\zeta$ and no other singularity is called a \emph{cylinder associated to $\zeta$}.
\end{definition}
In the case of cylindric suspension each edge $\zeta_x$ is a saddle connection and there are several cylinders that are associated to $\zeta_x$ but we emphasize that in general given a saddle connection in a translation surface there is no associated cylinder. Let $S$ be a cylindric suspension whose permutations are defined on the alphabet $\mathcal{A}$. The cylinders associated to an edge $\zeta_x$ which are of interest for our purpose are cylinders for which the boundary circles are obtained by a straight line in the polygonal representation joining the endpoints of $\zeta_x$ in the bottom circle to the endpoints of $\zeta_x$ in the bottom circle as in the left part of Figure~\ref{fig:removing_a_cylinder_in_cylindric_suspension}.

Let $S$ be a translation surface, $\zeta$ a saddle connection in $S$, $C$ a cylinder associated to $\zeta$ and $c_1$, $c_2$ its boundary circles. Denote by $S'$ the surface which is obtained from $S$ by removing the interior of $C$ and identifying $c_1$ and $c_2$ under the unique isometry $f: c_1 \rightarrow c_2$ that maps the endpoint of $\zeta$ in $c_1$ to the endpoint of $\zeta$ in $c_2$. 
In the surface $S'$ there is a saddle connection $c'$ which corresponds to the identified boundary circles $c_1$ and $c_2$ in $S$. The operation $(S,C) \rightarrow S'$ is invertible as soon as we know the saddle connection $c'$ in $S'$ and the parameters of the cylinder $C$ which is removed in $S$, namely its height $h \in (0,\infty)$ and a twist parameter $\theta \in S^1$. The converse operation $(S',c',h,\theta) \rightarrow S$ is called \emph{bubbling a handle} in \cite{KontsevichZorich03} and \emph{a figure eight operation} in \cite{EskinMasurZorich03}.

Consider a triple $(\tau_t,\,\tau_b,\,\sigma) \in S_\AAA \times S_\AAA \times S_\AAA$ satisfying~(\ref{eq:cyclic_constellation_equation}) and an associated cylindric suspension $S$. Let $\zeta_x$ be the edge in $S$ associated to a label $x$ in $\mathcal{A}$ and $C$ an associated cylinder whose boundary circles are straight line in the polygonal representation as in Figure~\ref{fig:removing_a_cylinder_in_cylindric_suspension}. The surface $S'$ obtained by removing the cylinder $C$ is still a cylindric suspension but of the triple $(\tau_t',\tau_b',\sigma')$ obtained by removing $x$ in the cycle decomposition of $\tau_t'$ and $\tau_b'$ as defined Section~\ref{subsection:recursive_construction}. While the choice of a cylinder associated to $x$ is not unique, the surface $S'$ is.  With our convention, the set of outgoing edges of each singularity $P$ of $S$ is invariant under the permutation $\sigma$. The cycle $c$ of $\sigma$ containg $x$ corresponds to the startpoint of $\zeta_x$ while the endpoint of $\zeta_x$ corresponds to the cycle of $\sigma$ containing $\tau_t(x)$. Proposition~\ref{prop:collapse_or_split_combinatoric} can then be rephrased in terms of translation surfaces, cylinders and strata.
\begin{proposition}
\label{prop:collapse_or_split_geometric}
Let $S \in \Omega \mathcal{M}(n_1,\ldots,n_k)$ be a translation surface, $\zeta$ a saddle connection in $S$ and $C$ a cylinder associated to $\zeta$. Let $S'$ be the surface obtained by removing the cylinder $C$ in $S$. If the endpoints of $\zeta$ corresponds to the same singularity of degree $\kappa_1$ in $S$ and the start and end of $\zeta$ are separated by an angle $(2a+1)\pi$ then the stratum of $S'$ is $\Omega \mathcal{M}(a,n_1-a-1,n_2,\ldots,n_k)$. If the endpoints of $\zeta$ corresponds to two different singularities of $S$ of degrees respectively $\kappa_1$ and $\kappa_2$ the the stratum of $S'$ is $\Omega \mathcal{M}(n_1+n_2-1,n_3, \ldots, n_k)$.
\end{proposition}

Let $S$ be a translation surface, $\Sigma \subset S$ its singularities, $\zeta$ a saddle connection in $S$ and $C$ a cylinder of $S$ associated to $\zeta$. Let $S'$ be the translation surface obtained by removing $C$ from $S$, $\Sigma' \subset S'$ its singularities, and $c'$ the saddle connection in $S'$ which corresponds to the identified boundary circles $c_1$ and $c_2$ of $C$. We define a map $\Psi: H_1(S' \backslash \Sigma';\ZZ/2) \rightarrow H_1(S \backslash \Sigma;\ZZ/2)$ which will be used to compare the spin parities of $S$ and $S'$.

Recall that the surgery operation $S \mapsto S'$, does not affect $S \backslash C$. Hence, if $\xi \subset S$ is a closed curve disjoint from the cylinder $C$, it defines a curve $\xi' \subset S'$. Let $\xi' \subset S'$ is a closed curve which intersects $c'$. We assume that the intersection is transverse. Let $\xi \subset S$ be the curve which co\"incides with $\xi'$ outside of $C$ and, for each intersection $P'$ of $\xi'$ and $c'$, we replace $P'$ by the unique geodesic segment in $C$ which joins the preimages $P_1 \in c_1$ and $P_2 \in c_2$ of $P'$ and do not intersect $\zeta$.

\begin{lemma}
\label{lemma:shrinking_map_at_homology_level}
Let $S$ , $\Sigma$, $S'$ and $\Sigma'$ as above. Then the map $\xi' \mapsto \xi$ defines a map $\Psi: H_1(S' \backslash \Sigma'; \ZZ/2) \rightarrow H_1(S \backslash \Sigma; \ZZ/2)$. Moreover $\Psi$ is injective, preserves the intersection forms and the winding numbers.
\end{lemma}

\begin{proof}
The map $\xi' \mapsto \xi$ is well defined on homology because it preserves boundaries. Let $\xi' \subset S' \backslash \Sigma'$ be a simple closed curve such that $[\xi] = 0 \in H_1(S \backslash \Sigma; \ZZ/2)$. Then there is a disc $D \subset S \backslash \Sigma$ such that $\xi = \partial D$. The disc $D$ goes down to a disc in $S'$ and shows that $[\xi'] = 0$.

If $\xi'$ is disjoint from $c'$, it is clear that the intersection with $\xi'$ is preserved and $w'(\xi') = w(\xi)$. Now if $\xi'$ is transverse to $c'$ then the pieces added to build $\xi$ are all parallel and in particular do not intersect and has no winding. As the preceding case, the intersection with $\xi'$ is preserved and $w'(\xi') = w(\xi)$.
\end{proof}

\subsubsection{Spin parity in the collapsing case}
\label{subsubsection:spin_behavior_collapsing}
Let $S$ be a translation surface with spin parity. Depending on the alternative of Proposition~\ref{prop:collapse_or_split_geometric}, the behavior of the spin structure is different. Let $\zeta$ be a saddle connection in $S$ whose endpoints are two different singularities of $S$, $C$ a cylinder associated to $\zeta$ and $S'$ the surface obtained by removing the cylinder $C$ in $S$. The genus of $S$ is the same as the one of $S'$ and we have the following result.
\begin{lemma}
\label{lemma:collapsing}
Let $S$, $\zeta$, $C$ and $S'$ as above. Then $S'$ has a spin parity and is the same as the one of $S$.
\end{lemma}

\begin{proof}
From Proposition~\ref{prop:collapse_or_split_geometric}, we know that if $S$ has a spin structure (meaning that all its singularities have degrees even multiples of $2\pi$) then $S'$ has also one. Recall that the spin structure of, respectively, $S$ and $S'$ are given by Arf invariants of quadratic forms $q_S$ and $q_{S'}$ on $H_1(S;\ZZ/2)$ and $H_1(S';\ZZ/2)$ (see Section~\ref{subsubsection:spin_invariant}).

Let $\Psi: H_1(S' \backslash \Sigma'; \ZZ/2) \rightarrow H_1(S \backslash \Sigma; \ZZ/2)$ be the map of Lemma~\ref{lemma:shrinking_map_at_homology_level}. As all singularities of $S$ and $S'$ are conical angles of odd multiple of $2\pi$ the winding numbers $w: H_1(S \backslash \Sigma; \ZZ/2) \rightarrow \ZZ/2$ and $w': H_1(S' \backslash \Sigma'; \ZZ/2) \rightarrow \ZZ/2$ are well defined on $H_1(S;\ZZ/2)$ and $H_1(S';\ZZ/2)$. In the collapsing case, the genus of $S$ equals the genus of $S'$ and hence the vector spaces $H_1(S;\ZZ/2)$ and $H_1(S';\ZZ/2)$ have the same dimension. 

As $\Psi$ is injective, it is an isomorphism. $\overline{\Psi}$ preserves the intersection form and the winding number, thus $q_{S'} = q_S \circ \overline{\Phi}$ and the Arf invariant of $q_{S'}$ and $q_S$ are equal. This proves that $S$ and $S'$ have the same spin parity.
\end{proof}

\subsubsection{Spin parity in the splitting case}
\label{subsubsection:spin_behavior_splitting}
We now consider the case of a translation surface $S$ with a saddle connection $\zeta$ which has the same singularity $P \in S$ as endpoints. Let $C$ be a cylinder associated to $\zeta$. By Proposition~\ref{prop:collapse_or_split_geometric}, removing $C$ in $S$ gives a surface $S'$ whose genus is the one of $S'$ minus $1$. The start and the end of the geodesic $\zeta$ form an angle at the point $P$ which is an odd multiple of $\pi$ that we denote $(2a+1) \pi$ (see Proposition~\ref{prop:collapse_or_split_combinatoric} and Proposition~\ref{prop:collapse_or_split_geometric}). In order to get the recurrence for the numbers $d(p)$, we have two cases to treat:
\begin{itemize}
\item $S$ and $S'$ have a spin parity, which corresponds to $a$ odd (Lemma~\ref{lemma:splitting_odd}),
\item $S$ has a spin parity but $S'$ has not, which corresponds to $a$ even (Lemma~\ref{lemma:splitting_even}).
\end{itemize}
Similarly to Lemma~\ref{lemma:collapsing}, we have.
\begin{lemma}
\label{lemma:splitting_odd}
Let $S$ and $C$ as above. We assume that $S$ has a spin parity and that $a$ is odd. Then $S'$ obtained by removing $C$ in $S$ has a spin parity and is the same as the one of $S$.
\end{lemma}

\begin{figure}[!ht]
\centering
\picinput{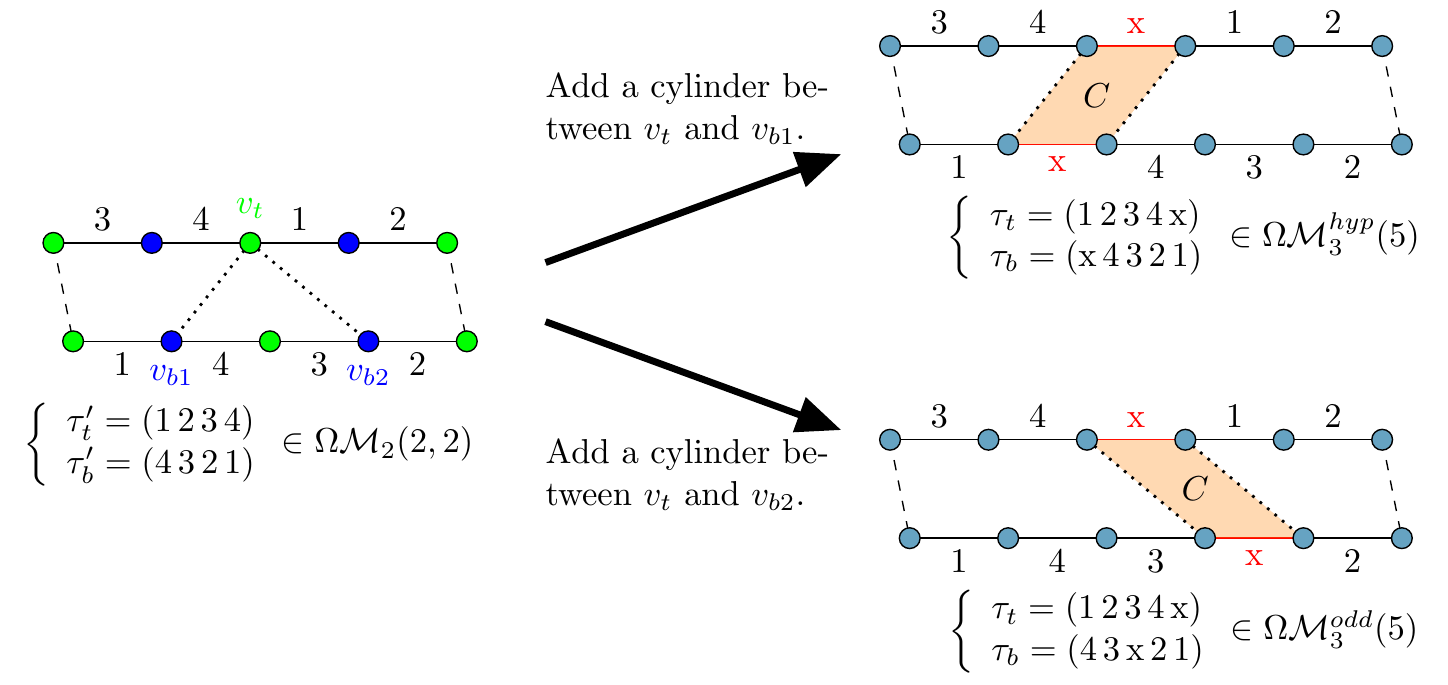}
\caption{The two ways of adding a cylinder to a cylindric suspension in $\Omega \mathcal{M}(2,2)$.}
\label{fig:collapsing_and_spin}
\end{figure}

\begin{proof}
We consider the maps $\Psi$ and $\overline{\Psi}$ of Lemma~\ref{lemma:shrinking_map_at_homology_level}. The map $\overline{\Psi}$ identifies a subspace of codimension $2$ of $H_1(S;\ZZ/2)$ with $H_1(S';\ZZ/2)$. Let $c$ be a circumference of the cylinder $C$. Then, the symplectic complement of $H_1(S';\ZZ/2)$ in $H_1(S;\ZZ/2)$ is the subspace $M = \ZZ/2\ [c] \oplus \ZZ/2\ [\zeta]$. Hence $q_S \simeq q_S' \oplus q_S|_M$ and, as the Arf invariant is additive, to compare the Arf invariant of $q_S$ and $q_{S'}$ we compute the Arf invariant of $q_S|_M$.

As $\zeta$ is geodesic and its start and end are separated by an angle $(2a+1) \pi$ we have $w(\zeta) = a \mod 2$ and hence $q_{S}([\tilde{\zeta}]) = a + 1 = 0 \mod 2$. On other hand $q_{S}([c]) = 1$, and from Theorem~\ref{thm:classification_of_non_degenerate_quadratic_binary_forms} we get that $Arf(q_S|_M) = 0$. Thus $q_S$ and $q_{S'}$ have the same Arf invariant which proves that $S$ and $S'$ have the same spin parity.
\end{proof}
Now, we treat the case of $a$ even. The surface $S'$ obtained after removing the cylinder has no spin but the surface $S$ can have one. In the following lemma the surface $S'$ is fixed and we count how many surfaces of each spin parity we get by the procedure of adding a cylinder. Let $(\tau_t',\tau_b',\sigma')$ the combinatorial datum associated to a cylindric suspension $S'$. We assume that the profile $p$ of $S'$ contains only odd numbers excepted two, $m_t$ and $m_b$ and we write $p = (m_t,m_b) \uplus q$. Let $P_t$ and $P_b$ be the two singularities of $S'$ of conical angles respectively $m_t$ and $m_b$. We fix a vertex $v_t$ corresponding to $P_t$ in the top circle of $S'$. Consider all saddle connections that joins $v_t$ to a vertex associated to $P_b$ in the bottom line of the circle of $S'$ (see Figure~\ref{fig:collapsing_and_spin}). The following is similar to Lemma~14.4 of \cite{EskinMasurZorich03}.
\begin{lemma}
\label{lemma:splitting_even}
Let $S'$, $P_t,P_b \in S'$, $m_t$, $m_b$ and $v_t$ as above. Then there are $m_b$ vertices in the bottom circle of $S'$ associated to $P_b$. Amongst the $m_b$ cylindric suspension obtained by adding a cylinder to $(S',[vt,v_b])$ where $v_b$ is a vertex associated to $P_b$ in the bottom line, half of them have an odd spin parity and half of them have an even spin parity.
\end{lemma}

\begin{proof}
There are exactly $m_b$ vertices associated to $P_b$ in the bottom circle as the conical angle at $P_b$ is $2 \pi m_b$. We fix $v_b$ associated to $P_b$ in the bottom cylinder. We use the same strategy as in Lemma~\ref{lemma:splitting_odd}, we use a map $H_1(S';\ZZ/2) \rightarrow H_1(S;\ZZ/2)$ and then look at the symplectic complement of its range.

Consider a small neighborhoods $V_b$ of $P_b$ in $S'$ and $c$ the saddle connection that joins $v_t$ to $v_b$. Any other saddle connection between $v_t$ and a representative of $P_b$ in the bottom circle can be obtained by adding to $c$ an arc of circle contained in $V_b$.
Hence each curves that joins $v_t$ to a representative of $P_b$ in the bottom line can be numeroted with respect to the angle from $c$. We denote them by $c_0=c$, $c_1$, \ldots, $c_{m_b-1}$. Let $S_j$, $j=0,\ldots,m_b-1$, be the surface obtained by adding a cylinder corresponding to $c_j$ and $q_j$ its associated quadratic form. The contribution of the module $M_j = \ZZ/2 [c_j] \oplus \ZZ/2 [x_j] \subset H_1(S;\ZZ/2)$ to the spin structure is $q_j(c_j) = q(c) + j \mod 2$ and $q_j(x_j) = 1$. In particular $Arf(q_j) = Arf(q) + j \mod 2$ which proves the lemma.
\end{proof}

\subsection{Formulas for $c(p)$ and $d(p)$}
\label{subsection:explicit_formulas}
In this section we prove formulas for the numbers $c(p)$ and $d(p)$. We will use two notations for partitions of an integer $n$. Either $p = (n_1,n_2,\ldots,n_k)$ where $n_1$, \ldots, $n_k$ are positive integers whose sum are $n$. Or $p = (1^{e_1}, 2^{e_2}, \ldots, n^{e_n})$ where $e_i$ denotes the number of times $i$ occurs in $p$. The numbers $e_i$ satisifies $\sum e_i\ i = n$.

\subsubsection{Marked points}
\label{subsubsection:labeled_standard_permutations_marked_points}
We first consider the presence of $1$ in the integer partition $p$. They correspond to marked point in the associated cylindric suspension $S$. See for example Figure~\ref{fig:cylindric_suspension} where the vertex represented by a square with outgoing edge $b$ is a marked point.

\begin{proposition}
\label{prop:marked_points_labeled_standard_permutations}
We have $c((1^n)) = d((1^n)) = (n-1)!$ and, more generally, if $p$ is a partition of the integer $n$ and $k$ is a non negative integer then
\[
c(p \uplus (1^k)) = \frac{(n+k-1)!}{(n-1)!}\ c(p).
\]
If moreover $p$ has only odd parts, then 
\[
d(p \uplus (1^k)) = \frac{(n+k-1)!}{(n-1)!}\ d(p).
\]
\end{proposition}

\begin{proof}
The identity permutation $1 \in S_n$ is the only element with profile $(1^n)$. On the other hand, the solutions of the form $(\tau_t,\tau_b,1)$ of Equation~(\ref{eq:cyclic_constellation_equation}) are given by $(c,c,1)$ where $c$ is any $n$-cycles. Thus $c((1^n)) = (n-1)!$ is the number of $n$-cycles in $S_n$. As the partition $(1^n)$ corresponds to a torus (a surfaces with genus~$1$), it is well known that the spin is odd. Hence $d((1^n)) = c((1^n))$. More generally, adding marked points in a surface do not modify the spin parity.

We denote by $C_n$ the set of $n$-cycles in $S_n$. Let $p$ be a partition of $n$ and $\sigma' \in S_n$ whose conjugacy class is $p$. Let
\[
E' = \{(\tau_t',\tau_b') \in C_n \times C_n|\ \tau'_t \ {\tau'_b}^{-1} = \sigma'\}.
\]
Let $\sigma \in S_{n+1}$ be the permutation which equals $\sigma'$ on $\{1,\ldots,n\}$ and such that $\sigma(n+1)=(n+1)$ and
\[
E = \{(\tau_t,\tau_b) \in C_{n+1} \times C_{n+1}|\ \tau_t \ \tau_b^{-1} = \sigma\}.
\]
We claim that there is a canonic bijection $E \rightarrow E' \times \{1,\ldots,n\}$. The conclusion of the lemma follows from the claim which we prove now.

The map $E \rightarrow E'$ on the first factor correspond to remove $(n+1)$ in the cycles $\tau_t$ and $\tau_b$ as in Section~\ref{subsection:recursive_construction}. The map $E \rightarrow \{1,\ldots,n\}$ on the second factor is $(\tau_t,\tau_b) \mapsto {\tau_t}^{-1}(n+1)$. As $\sigma(n+1) = n+1$, we have ${\tau_t}^{-1}(n+1) = {\tau_b}^{-1}(n+1)$. The preimage $(\tau_t,\tau_b)$ of the element $(\tau'_t,\tau'_b, x) \in E' \times \{1,\ldots,n\}$ is given by
\[
\tau(i) = \left\{
\begin{array}{ll}
n+1 & \text{if $i=x$}, \\
\tau(x) & \text{if $i=n+1$}, \\
\tau(i) & \text{otherwise},
\end{array}
\right.
\quad \text{for $\tau = \tau_t$ or $\tau = \tau_b$}.
\]
\end{proof}

\subsubsection{Two formulas for $c(p)$}
We first give a recurrence formula for the number $c(p)$ of labeled standard permutations in the stratum associated to $p$. The initialization $c((1)) = 1$ of the recurrence can be considered as a particular case of Proposition~\ref{prop:marked_points_labeled_standard_permutations}.
\begin{theorem}[\cite{Boccara80} prop. 4.2.]
\label{thm:std_perm_strata_recurrence}
Let $p = (n_1,\,\dots,\,n_k)$ be a partition of an integer $n \geq 2$, then
\[
c(p) = \sum_{i = 2}^k n_i\ c(p_{n_1 \odot n_i}) + \sum_{a=2}^{n_1-1} c(p_{n_1 | a}).
\]
\end{theorem}

\begin{proof}
Let $\sigma \in S_n$ whose conjugacy class is $p$ such that the length of the cycle containing $n$ is $n_1$. As in the proof of Proposition~\ref{prop:marked_points_labeled_standard_permutations} we set
\[
E(\sigma) = \{(\tau_t,\tau_b) \in C_n \times C_n; \tau_t \ \tau_b^{-1} = \sigma \}.
\]
To an element $(\tau_t,\tau_b) \in E$ we associate $(\tau_t',\tau_b',\tau_t(n)) \in C_{n-1} \times C_{n-1} \times \{1,\ldots,n-1\}$ where $(\tau_t',\tau_b')$ is obtained from $(\tau_t,\tau_b)$ by removing $n$ in their cycle decomposition (see Section~\ref{subsection:recursive_construction}). The map $E \rightarrow C_{n-1} \times C_{n-1} \times \{1,\ldots,n-1\}$ is injective. As proved in Proposition~\ref{prop:collapse_or_split_combinatoric}, the conjugacy class of $\sigma' = \tau_t' \ {\tau'_b}^{-1}$ depends on the nature of the cycle of $\sigma$ that contains $\tau_t(n)$. The formula of the theorem follows by summing over all possibilities for $\tau_t(n)$. The first sum corresponds to the cases where $\tau_t(n)$ is in a different cycle from the one of $n$. The second sum corresponds to the cases where $n$ and $\tau_t(n)$ are in the same cycle.
\end{proof}
Boccara in \cite{Boccara80} find an explicit formula from the recurrence of Theorem~\ref{thm:std_perm_strata_recurrence} using an identity involving a polynom and integration.
\begin{theorem}[\cite{Boccara80}]
\label{thm:Boccara_formula_for_c}
Let $p=(n_1,n_2,\ldots,n_k)$ be a partition of the integer $n$. Then, we have
\[
c(p) = \frac{2 (n-1)!}{n+1} \left( \sum_{q \subset (n_2,n_3,\ldots,n_k)} (-1)^{s(q)-l(q)} \binom{n}{s(q)}^{-1} \right).
\]
\end{theorem}
\noindent From the theorem, we deduce several explicit values
\begin{corollary}
\label{corollary:c_explicit_formulas}
Let $n=2k+1$ then
\[
c((n)) = \frac{2\ (n-1)!}{n+1}.
\]
Let $n = n_1 + n_2 \equiv 0 \mod 2$, then
\[
c((n_1,n_2)) = \frac{2\ (n-1)!}{n+1} \ \frac{\binom{n}{n_1,n_2} + (-1)^{n_1+1}}{\binom{n}{n_1,n_2}}.
\]
\end{corollary}
\noindent We also have
\begin{proposition}
\label{prop:c_explicit_value_2k}
Let $k$ be a positive integer and $n=2k$ then
\[
c((2^k)) = \frac{2 \ (n-1)!}{n+2}.
\]
\end{proposition}
Using the representation theory of the symmetric group A. Goupil and G. Schaeffer \cite{GoupilSchaeffer98} gave an explicit formula for more general numbers than $c(p)$. Their formula has the advantage of containing only positive numbers. In our particular case we get
\begin{theorem}[\cite{GoupilSchaeffer98}]
\label{thm:GoupilSchaeffer_formula_for_c}
Let $p$ be a partition of the integer $n$ with length $k$. We set $g = (n-k)/2$. Then we have
\[
c(p) = \frac{z_p}{2^{2g}} \sum_{g_1+g_2=g} \frac{(2g_1)!}{2g_1+1} \binom{n-1}{2g_1} S_{k,g_2}(p),
\]
where $S_{k,g} \in \QQ[x_1,x_2,\ldots,x_k]$ is the symmetric polynomial
\[
S_{k,g}(x_1,x_2,\ldots,x_k) = (k+2g-1)! \sum_{(p_1,p_2,\ldots,p_k) \models g} \ \prod_{i=1}^k \frac{1}{2 p_i + 1} \binom{x_i-1}{2p_i},
\]
where $(p_1,\ldots,p_k) \models g$ design the set of $k$-tuples $(p_1,\ldots,p_k)$ of non-negative integers whose sum is $g$. And $z_p$ is the cardinality of the centralizer of any permutation in the conjugacy class associated to $p$. Writing $p = (1^{e_1},2^{e_2},\ldots,n^{e_n})$ in exponential notation we have
\[
z_p = \prod_{i=1}^n e_i!\ i^{e_i}.
\]
\end{theorem}
In \cite{Walkup79}, D. Walkup made a conjecture about the asymptotic behavior of the numbers $c$ which was proved few years later by R. Stanley in \cite{Stanley81}.
\begin{theorem}[\cite{Walkup79},\cite{Stanley81}]
\label{thm:gamma_asymptotic}
Let $(p_i)_{i \geq 0}$ be a sequence of partition of integers $(n_i)_{i \geq 0}$ such that $n_i$ tends to infinity and the number of $1$ in $p_i$ is $O(\log(n_i))$ then 
\[
c(p_i) \sim 2 (n_i-2)! (1 + o(1)).
\]
\end{theorem}
The asymptotic behavior of the above theorem proves that in Boccara's formula (Theorem~\ref{thm:Boccara_formula_for_c}) the only contribution comes from the factor $\frac{2(n-1)!}{n+1}$ and the sum in parentheses is asymptotically $(1+o(1))$. For the particular cases in Corollary~\ref{corollary:c_explicit_formulas} this fact is clear.

\subsubsection{A formula for $d(p)$}
For an integer partition $p$ whose parts are odd numbers, recall that $c_1(p)$ and $c_0(p)$ denote the number of standard permutations with fixed labels and respectively odd and even spin parity. We have $c(p) = c_1(p) + c_0(p)$ and $d(p) = c_1(p) - c_0(p)$. As for $c$, we first prove a recurrence formula and then solve it explicitely.

The recurrence formula is similar to Theorem~\ref{thm:std_perm_strata_recurrence}.
\begin{theorem}
\label{thm:std_perm_spin_recurrence}
Let $p=(n_1,\ldots,n_k)$ be an integer partitions with odd parts then
\[
d(p) = \sum_{i = 2}^k n_i\ d(p_{n_1 \odot n_i}) + \sum_{a = 1 \atop a \equiv 1 \mod 2}^{n_1-2} d(p_{n_1|a}). 
\]
\end{theorem}

\begin{proof}
The proof is identic to the one of Theorem~\ref{thm:std_perm_strata_recurrence}. We fix a permutation $p$ and an element $\sigma \in S_n$ such that the conjugacy class of $\sigma$ is $p$. We assume that the cycle containing $n$ has length $n_1$.

Let $E_s$ be the set of standard permutations $(\tau_t,\tau_b)$ with labels $\sigma$ and spin parity $s \in \{0,1\}$. According to the position of $\tau_t^{-1}(n)$ we separate $E_s$ in different subsets. 

If $\tau_t^{-1}(n)$ and $n$ are in different cycles, then we apply the Lemma~\ref{lemma:collapsing} and we get that their number is
\[
\sum_{j \not= i}^k n_j c_s(p_{n_i \odot n_j}).
\]
If $n$ and $\tau_t^{-1}(n)$ are in the same cycle, then we differentiate the case $a$ odd and $a$ even (see Section~\ref{subsubsection:removing_a_cylinder}). For $a$ even, Lemma~\ref{lemma:splitting_odd} gives that the total number of such standard permutation is
\[
\sum_{a \equiv 1 \, [2]} c_s(p_{n_1|a}).
\]
For $a$ odd, Lemma~\ref{lemma:splitting_even} implies that their number is
\[
\frac{1}{2}\ \sum_{a \equiv  0 \, [2]}  c(p_{n_1|a}).
\]
As this last term does not depend on the spin parity $s$, it cancels in the difference $c_1(p) - c_0(p)$.
\end{proof}
\noindent The formula for the numbers $d(p)$ is given by the following.
\begin{theorem}
\label{thm:formula_for_d}
Let $p$ be an integer partition with only odd parts, then the number $d(p)$ depends only on the sum $n$ and the length $k$ of $p$. Set $g=(n-k)/2$ then
\[
d(p) = \frac{(n-1)!}{2^g}.
\]
\end{theorem}

\begin{proof}
Set $\tilde{d}(n,k) := (n-1)!/2^{(n-k)/2}$. Those numbers satisfy the recurrences
\[
\tilde{d}(n+1,k+1) = n \ \tilde{d}(n,k) \quad \text{and} \quad \tilde{d}(n+1,k-1) = \frac{n}{2} \ \tilde{d}(n,k).
\]
On the other hand if $i,j \in \{1,\ldots,k\}$ and $a \in \{1,\ldots,n_i-2\}$, we have for the sum $s(p_{m_i \odot m_j}) = s(p_{m_i|a})=s(p)-1$ and for the length $l(p_{n_i \odot n_j}) = l(p)-1$ and $l(p_{n_i|a}) = l(p)+1$. It is then straightforward to check that $\tilde{d}$ satisfies the same recurrence as the formula given in Theorem~\ref{thm:std_perm_spin_recurrence}. The initial value needed to start the recurrence is the one for the only partition of $1$ which is $p=(1)$. But $\tilde{d}(1,1) = 1 = d((1))$. Hence $d(p) = \tilde{d}(s(p),l(p))$ for all partitions with odd parts.
\end{proof}

\section{From standard permutations to cardinality of Rauzy classes}
\label{section:from_standards_to_all}
We now prove a recurrence formula for the numbers $\gamma^{irr}(p)$ (resp. $\delta^{irr}(p) = \gamma_1^{irr}(p) - \gamma_0^{irr}(p)$) in terms of the number of standard permutations $\gamma^{std}(p)$ (resp. $\delta^{std}(p) = \gamma_1^{std}(p) - \gamma_0^{std}(p)$). We relate the latter ones to the numbers $c(p)$ and $d(p)$ computed in the preceding section. The recurrence formula is based on the construction of suspensions for any permutation (non necessarily irreducible) and a geometrical analysis of the concatenation of permutations.

\subsection{Irreducibility, concatenation and non connected surfaces}
\label{subsection:concatenation_irreducibility}
\subsubsection{Concatenation and irreducible permutations}
Let $\pi_1$ (resp. $\pi_2$) be a labeled permutation on the alphabet $\mathcal{A}_1$ (resp. $\mathcal{A}_2$). The \emph{concatenation} $\pi_1 \cdot \pi_2$ is the labeled permutation on the disjoint union $\mathcal{A}_1 \sqcup \mathcal{A}_2$ defined by
\[
\left(\begin{array}{lll}
  a_1 & \ldots & a_{n_1} \\
  b_1 & \ldots & b_{n_1}
\end{array}\right)
\cdot
\left(\begin{array}{lll}
  a'_1 & \ldots & a'_{n_2} \\
  b'_1 & \ldots & b'_{n_2}
\end{array}\right)
=
\left(\begin{array}{llllll}
  a_1 & \ldots & a_{n_1} & a'_1 & \ldots & a'_{n_2} \\
  b_1 & \ldots & b_{n_1} & b'_1 & \ldots & b'_{n_2} 
\end{array}\right).
\]
The concatenation of two reduced permutations can be defined from the section $\pi \mapsto (id,\pi)$ and projection $(\pi_t, \pi_b) \mapsto \pi_b \circ \pi_t^{-1}$ (see Section~\ref{subsubsection:labeled_permutations}). More precisely, let $\pi_1$ and $\pi_2$ be two reduced permutations of lengths $n_1$ and $n_2$. The concatenation $\pi = \pi_1 \cdot \pi_2$ is the permutation of length $n_1+n_2$ defined by
\[
\pi(i) = \left\{
\begin{array}{ll}
  \pi_1(i) & \text{if $1 \leq i \leq n_1$}, \\
  \pi_2(i-n_1) + n_1 & \text{if $n_1 + 1 \leq i \leq n_1+n_2$}.
\end{array}
\right.
\]

One has the following elementary.
\begin{proposition} \label{prop:UniqueFactorizationOfPermutation}
A permutation $\pi \in S_n$ is irreducible if and only if it can not be written as a non trivial concatenation.

Each (reduced or labeled) permutation has a unique decomposition in irreducible permutations.
\end{proposition}
As an example, we write in the table below the decomposition of the reducible permutations of length $4$. We call \emph{class} of a permutation $\pi$ the ordered list of the lengths of the irreducible components of $\pi$ (which is a \emph{composition} of $4$, e.g. an ordered list of positive integers whose sum is sum $4$).
\[
\begin{array}{|c|c|c|}
\hline
\text{permutation} & \text{decomposition} & \text{class} \\
\hline
(1234) &  (1) \cdot (1) \cdot (1) \cdot (1) & [1,1,1,1] \\
\hline
(1243) & (1) \cdot (1) \cdot (21) & [1,1,2] \\
(1324) & (1) \cdot (21) \cdot (1) & [1,2,1] \\
(2134) & (21) \cdot (1) \cdot (1) & [2,1,1] \\
\hline
(2143) & (21) \cdot (21) & [2,2] \\
(1342) & (1) \cdot (231) & [1,3] \\
(1423) & (1) \cdot (312) & [1,3] \\
(1432) & (1) \cdot (321) & [1,3] \\
(2314) & (231) \cdot (1) & [3,1] \\
(3124) & (312) \cdot (1) & [3,1] \\
(3214) & (321) \cdot (1) & [3,1] \\
\hline
\end{array}
\]

As a corollary, we get a formula relating factorial numbers $n! = |S_n|$ to $p(n) = |S_n^o|$.
\begin{corollary}
Let $f(n)$ be the number of irreducible permutations in $S_n$. Then
\begin{equation}
n! = \sum_{k=1}^n \ \sum_{c_{1}+\ldots+c_k=n} \ f(c_1)\, f(c_2)\, \ldots\, f(c_k),
\label{eq:irreducible}
\end{equation}
\end{corollary}

\subsubsection{Suspensions of reducible permutations}
\label{subsubsection:suspensions_of_reducible_permutations}
Let $\pi_1$ and $\pi_2$ be two labeled permutations on the alphabet $\mathcal{A}$ of lengths respectively $n_1$ and $n_2$ and $\pi = \pi_1 \cdot \pi_2 = (\pi_t,\pi_b)$ their concatenation of  length $n = n_1 + n_2$. If $\zeta \in \CC^{\mathcal{A}}$ then
\begin{equation}
\sum_{1 \leq j \leq n_1} \zeta_{\pi_t^{-1}{j}} = \sum_{1 \leq j \leq n_1} \zeta_{\pi_b^{-1}{j}}.
\label{eq:relation_for_zeta}
\end{equation}
Thus there is no suspension data for $\pi$ (see Section~\ref{subsubsection:suspensions}). But if $\pi_1$ and $\pi_2$ are irreducible, we can assume that $n_1$ is the only index such that (\ref{eq:relation_for_zeta}) holds.

\begin{definition}
Let $\pi$  be a labeled permutation on the alphabet $\mathcal{A}$ and $\pi_1 \cdot \pi_2 \cdot \ldots \cdot \pi_k$ its decomposition in irreducible permutations. Let $\mathcal{A}_j$ be the alphabet of $\pi_j$. A \textit{suspension data} for $\pi$ is a vector in $\zeta \in \CC^{\mathcal{A}}$ such that each $\left(\zeta_\alpha\right)_{\alpha \in \mathcal{A}_j}$ is a suspension data for the irreducible permutation $\pi_j$.
\end{definition}
In the case of the irreducible permutation on $1$ letter $\pi = \binom{A}{A}$, the suspension datum is an element $\zeta_A \in \RR_+ \times i\RR \subset \CC$.

\begin{figure}[!ht]
\centering
\picinput{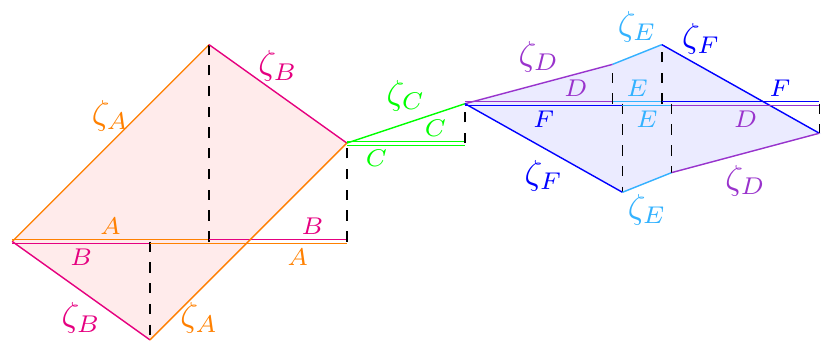}
\caption{A suspension of the reducible permutation $\binom{A B C D E F}{B A C F E D} = \binom{A B}{B A} \cdot \binom{C}{C} \cdot \binom{D E F}{F E D}$.}
\label{figure:reducible_suspension_iet}
\end{figure}

Let $\pi$ and $\zeta$ as in the above definition. Then, as for suspension of irreducible permutations in Section~\ref{subsubsection:suspensions}, we build two broken lines $L_t$ and $L_b$ made, respectively, of the concatenation of the vectors $\zeta_{\pi_t^{-1}(j)}$ and $\zeta_{\pi_b^{-1}(j)}$. The surface obtained by identifying the side $\zeta_\alpha$ on $L_t$ with the side $\zeta_\alpha$ on $L_b$ is a sequence $S_1$, $S_2$, \ldots, $S_k$ of translation surfaces such that $S_j$ and $S_{j+1}$ are connected at a singularity. In the case of the degenerate permutation $\binom{A}{A}$ the surface associated to $\zeta_A \in \RR_+ \times \RR$ corresponds to a (degenerate) sphere with two conical singularities of angle $0$. We take as convention that the stratum of $\binom{A}{A}$ is $\Omega \mathcal{M}(0,0)$ (see Figure~\ref{figure:reducible_suspension_iet}).


\subsubsection{Marking of a permutation}
\label{subsubsection:marking}
Let $\pi_1$ and $\pi_2$ be two permutations. We want to deduce the profile of the permutation $\pi = \pi_1 \cdot \pi_2$ as defined in Section~\ref{subsubsection:profile} from the profiles of $\pi_1$ and $\pi_2$. We first look at an example with the following permutations
\begin{equation}
\label{eq:permutations_product_profiles}
\pi_1 = \binom{1\,2\,3\,4\,5}{3\,5\,4\,2\,1}
\quad \text{and} \quad 
\pi_2 = \binom{1\,2\,3\,4\,5}{2\,5\,4\,1\,3}.
\end{equation}
Both permutations have have profile $(3,1)$ but the products $\pi_1 \cdot \pi_1$, $\pi_1 \cdot \pi_2$, $\pi_2 \cdot \pi_1$ and $\pi_2 \cdot \pi_2$ have respectively profiles $(5,3,1)$, $(7,1,1)$, $(3,3,3)$ and $(5,3,1)$. In a product $\pi_1 \cdot \pi_2$ the permutations are glued at the right of $\pi_1$ and the left of $\pi_2$. To keep track of left and right, we consider profile of permutation with an additional data which encodes the configuration of the two singularities at both extremities of the permutation. In the introduction, we defined markings in term of suspension. We give here a more combinatorial version based on the interval diagram of a permutation defined in Section~\ref{subsubsection:profile}.
\begin{definition}
\label{def:markings}
Let $\pi$ be a permutation, $\Gamma$ its interval diagram and $c_l$ (resp. $c_r$) be the cycle in $\Gamma$ that corresponds to the left (resp. right) endpoint of $\pi$.

If $c_l = c_r$, let $m$ be the length of $c_l$ and $2a$ be the number of edges in $\Gamma$ between the outgoing edge on the left of $\pi$ and the incoming edge on the right of $\pi$. The \emph{marking} of $\pi$ is the couple $(m,a)$ which we call a marking of the \emph{first type} and denote by $m|a$.

If $c_l \not= c_r$, let $m_l$ and $m_r$ be respectively the lengths of $c_l$ and $c_r$. The \textit{marking} of $\pi$ is the couple $(m_l,m_r)$ which we call a marking of the \emph{second type} and denote by $m_l \odot m_r$. 
\end{definition}
The notation similar to the one in Definition~\ref{def:splitting_and_collapsing_of_partition} is explained by the Corollaries~\ref{corollary:c_and_gamma} and \ref{corollary:d_and_delta} below.

For the permutations $\pi_1$ and $\pi_2$ defined in~(\ref{eq:permutations_product_profiles}) the interval diagrams are respectively
\[
\Gamma(\pi_1) = (\outgos{(3,1)}\ \incoms{2})\ (\incoms{3}\ \outgos{4}\ \incoms{(5,1)}\ \outgos{2}\ \incoms{4}\ \outgos{5}))
\quad \text{and} \quad
\Gamma(\pi_2) =  (\outgos{3}\ \incoms{1}\ \outgos{(2, 1)}\ \incoms{4}\ \outgos{5}\ \incoms{2})\ (\outgos{4}\ \incoms{(5, 3)}).
\]
Hence the markings are respectively $1 \odot 3$ and $3 \odot 1$. Examples of a marking of the first type with profile $(3,1)$ are given by the permutations
\[
\pi_3 = \binom{1\,2\,3\,4\,5}{2\,4\,5\,1\,3}
\quad
\pi_4 = \binom{1\,2\,3\,4\,5}{4\,5\,3\,2\,1}
\quad
\pi_5 = \binom{1\,2\,3\,4\,5}{2\,5\,3\,4\,1}
\]
The interval diagrams and markings of the above permutations are respectively
\begin{align*}
\Gamma(\pi_3) & = (\outgos{(2, 1)}\ \incoms{(5, 3)}\ \outgos{4}\ \incoms{2})\ \outgos{3}\ \incoms{1}) (\outgos{5}\ \incoms{4}) & \text{with marking $3|0$}, \\
\Gamma(\pi_4) & = (\outgos{(4, 1)}\ \incoms{2}\ \outgos{3}\ \incoms{(5, 1)}\ \outgos{2}\ \incoms{3}\ )\ (\outgos{5}\ \incoms{4}) & \text{with marking $3|1$}, \\
\Gamma(\pi_5) & = (\outgos{(2, 1)}\ \incoms{4}\ \outgos{5}\ \incoms{2}\ \outgos{3}\ \incoms{(5, 1)})\ (\incoms{3}\ \outgos{4}) & \text{with marking $3|2$}.
\end{align*}
Let $p$ be an integer partition. The markings that occur in a permutation $\pi$ with profile $p$ are
\begin{itemize}
\item the markings $m|a$ where $m \in p$ and $a \in \{0,\ldots,m-1\}$,
\item the markings $m_1 \odot m_2$ where $m_1,m_2 \in p$ and $m_1 \not= m_2$,
\item the markings $m \odot m$ for $m$ which appears at least twice in $p$.
\end{itemize}
We remark that for a permutation $\pi$ with marking of the first type $m|a$ the number $a$ belongs to $\{0,\ldots,m-1\}$ whereas for a standard permutation $a$ belongs to $\{1,\ldots,m-2\}$.

\begin{definition}
Let $\pi$ be a permutation with profile $p$ and marking $m|a$ (resp. $m_l \odot m_r$). The \textit{marked profile} of $p$ is the couple $(m|a, p')$ (resp. $(m_l \odot m_r, p')$) where $p'$ is the integer partition $p \backslash (m)$ (resp. $p \backslash (m_l,m_r)$).
\end{definition}
We naturally extend the definition of $\gamma$, $\gamma^{irr}$, $\gamma^{std}$, $\delta$, $\delta^{irr}$ and $\delta^{std}$ to marked profiles.

\subsubsection{Profile and spin parity of a concatenation $\pi_1 \cdot \pi_2$}
We now answer to the question asked previously about the profile of a concatenation. The lemma below expresses the marked profile of a concatenation in terms of the marked profiles of its irreducible components.
\begin{lemma}
\label{lemma:profile_concatenation}
Let $\pi_1$ and $\pi_2$ be two permutations and let $\pi = \pi_1 \cdot \pi_2$ be their concatenation. The following array shows how deduce the marked profile of $\pi$ from the marked profiles of $\pi_1$ and $\pi_2$.
\[
\begin{array}{|c|c||c|}
\hline
\text{marked profile for $\pi_1$} &
\text{marked profile for $\pi_2$} & 
\text{marked profile for $\pi$} \\
\hline
(m|a,\,p') & (n|b,\,q')
 & (m+n+1|a+b,\ p' \uplus q') \\
\hline
(m|a,\,p') & (n_l \odot n_r,\,q')
 & (m + n_l + 1 \odot n_r,\,p' \uplus q') \\
\hline
(m_l \odot m_r,\,p') & (n|b,\,q')
 & (m_l \odot m_r+n+1,\,p' \uplus q') \\
\hline
(m_l \odot m_r,\, p') & (n_l \odot n_r,\, q')
 & (m_l \odot n_r,\,p' \uplus q' \uplus (m_r+n_l+1)\\
\hline
\end{array}.
\]
\end{lemma}
In particular, a concatenation $\pi_1 \cdot \pi_2$ has a marking of the first type if and only if both of $\pi_1$ and $\pi_2$ have a marking of the first type.

\begin{proof}
Let $\Gamma$ (resp. $\Gamma_1$ and $\Gamma_2$) be the interval diagram of $\pi$ (resp. $\pi_1$ and $\pi_2$). Let $c_1$ (resp. $c_2$) be the cycle associated to the right of $\pi_1$ (resp. the left of $\pi_2$). The diagram $\Gamma$ is built from the disjoint union of $\Gamma_1$ and $\Gamma_2$ by gluing the cycles $c_1$ and $c_2$. More precisely, let
\[
\pi_1=
\left( \begin{array}{c}
x_1\ \ldots\ x_n \\
y_1\ \ldots\ y_n
\end{array} \right)
\quad \text{and} \quad
\pi_2=
\left( \begin{array}{c}
 x'_1\ \ldots\ x'_{n'} \\
 y'_1\ \ldots\ y'_{n'}
\end{array} \right).
\]
Then $c_1=\left(\incoms{(x_n,y_n)}\ A\right)$ and $c_2=\left(\outgos{(y'_1,x'_1)}\ A'\right)$ where $A$ and $A'$ design blocks of letters. In the concatenation $\pi = \pi_1 \cdot \pi_2$, the cycles $c_1$ and $c_2$ are glued into $c=\left(\outgos{x'_1}\ A'\ \outgos{y'_1}\ \incoms{y_n}\ A\ \incoms{x_n}\right)$. Hence, the length of $|c|$ is $|c_1| + |c_2| + 1$. In particular the profile of $p$ of $\pi$ can be computed from the profiles $p_1$ and $p_2$ of respectively $\pi_1$ and $\pi_2$ as $p = (p_1 \backslash (|c_1|)) \uplus (p_2 \backslash (|c_2|)) \uplus (|c_1|+|c_2|+1)$. We have proved how the profile of a concatenation $\pi = \pi_1 \cdot \pi_2$ can be deduced from the profiles and markings of its components $\pi_1$ and $\pi_2$. We now consider the marking of the permutation $\pi$.

We treat only the case of two markings of type one, the other being similar. We keep the notation $\Gamma_1$, $\Gamma_2$, $c_1$ and $c_2$ as above. The cycle $c_1$ (resp. $c_2$) of $\Gamma_1$ (resp. $\Gamma_2$) which corresponds to the right of $\pi_1$ (resp. the left of $\pi_2$) can be written as $c_1 = \left(\outgos{(y_1,x_1)}\ A\ \incoms{(x_n,y_n)}\ B\right)$ (resp. $c_2 = \left(\outgos{(y'_1,x'_1)}\ A'\ \incoms{(x'_{n'},y'_{n'})}\ B'\right)$) where $A$, $B$, $A'$ and $B'$ are blocks of letters. The angles in the marking are $a = |A|$ and $A'$. Those two cycles become one in $\pi$ which is 
\[
\left(\outgos{(y_1,x_1)}\ A\ \incoms{x_n}\ \outgos{x'_1}\ A'\
\incoms{(x'_{n'},y'_{n'})}\ B'\ \outgos{y'_1}\ \incoms{y_n}\ B\right).
\]
The angle $a$ (resp. $b$) in the marking of $\pi_1$ (resp. $\pi_2$) is the length of $A$ (resp. $A'$) divided by $2$. The structure of the cycle $c$ shows that the angle in the marking of $\pi$ is the length of $A\ x_n\ x'_1\ A'$ divided by $2$ which equals $a+b+1$.
\end{proof}
Now, we consider the spin parity of a permutation whose profile contains only odd parts. We would like to have a lemma similar to Lemma~\ref{lemma:profile_concatenation} which relates profile to the profiles of the irreducible components. But recall that the spin parity (see Section~\ref{subsubsection:spin_invariant}) is only defined when the profile contains only odd numbers. Hopefully Lemma~\ref{lemma:profile_concatenation} implies
\begin{corollary}
Let $p$ be an integer partition which contains only odd terms and $\pi$ a permutation with profile $p$. Then the profile of each irreducible component of $\pi$ contains only odd terms.
\end{corollary}
Hence, if $\pi$ is a permutation with profile $p$ containing only odd numbers, we can discuss about the spin parity of its components. The situation is simpler than the one in Lemma~\ref{lemma:profile_concatenation} as the spin parity does not depend on the structure of the endpoints of each component.
\begin{lemma}
\label{lemma:spin_parity_concatenation}
Let $p$ be a partition with odd parts and $\pi$ a permutation with profile $p$. Then the spin parity of $\pi$ is the sum mod 2 of the spins of the irreducible components of $\pi$.
\end{lemma}

\begin{proof}
Recall that the spin invariant of an irreducible permutation is the Arf invariant of a quadratic form $q_\pi$ on $\FF_2^{\mathcal{A}}$. It is geometrically defined on $H_1(S;\ZZ/2)$ where $S=S(\pi,\zeta)$ is any suspension of $\pi$ by
\[
q_\pi(x) = \left( w(x) + \#(\text{components of $x$}) + \# (\text{self intersection of $x$})\right) \mod 2.
\]
In the above formula, $w(\gamma)$ is the winding number of $\gamma$ which depends on the flat metric of the suspension while the other two are topological. Let $\pi$ be a permutation and $\pi_1 \cdot \pi_2 \cdot \ldots \cdot \pi_k$ its decomposition in irreducible components. Let $S=S(\pi,\zeta)$ be a suspension of $\pi$ and $S_j=S(\pi_j,\zeta_j)$ the associated suspension of each irreducible components. Then
\[
H_1(S;\ZZ/2) \simeq \bigoplus_{j=1}^k H_1(S_j;\ZZ/2)
\quad \text{and} \quad
q_{\pi} = \sum_{j=1}^k q_{\pi_j}.
\]
To complete the proof, we remark that the Arf invariant is additive (which follows from Theorem~\ref{thm:classification_of_non_degenerate_quadratic_binary_forms}).
\end{proof}

\subsection{Removing the ends of a standard permutation}
\label{subsection:degeneration}
Let $\pi$ be a standard permutation on the $n+2$ ordered symbols $\{0,1,\ldots,n+1\}$ (i.e. $\pi(0) = n+1$ and $\pi(n+1) = 0$). Consider the permutation $\tilde{\pi}$ on the $n$ letters $\{1,\ldots,n\}$ obtained by removing $0$ and $n+1$ in $\pi$. We call $\tilde{\pi}$ the \emph{degeneration} of $\pi$. As a permutation, $\tilde{\pi}$ corresponds to the restriction of the domain of $\pi$ from $\{0,1,\ldots,n+1\}$ to $\{1,\ldots,n\}$. The term degeneration comes from geometric consideration. Let $(\zeta^{(t)})_{t > 0}$ be a continuous sequence of suspensions of $\pi$ which converges to a vector $\bar{\zeta} \in \RR^2$ for which $\bar{\zeta}_0 = 0 = \bar{\zeta}_{n+1} = 0$ and $Re(\bar{\zeta}_k) > 0$ for all $0 < k < n+1$ and the imaginary part of $\bar{\zeta}$ satisfies the condition of suspension for $\tilde{\pi}$. Then the limit $\bar{\zeta}$ is a suspension of $\tilde{\pi}$ which do not live in the same stratum $\Omega \mathcal{M}(p)$ as $S$ but is obtained as a limit of a continuous family $(S_t) \subset \Omega \mathcal{M}(p)$ which degenerates for $t \to \infty$.

The degeneration operation is invertible and gives a bijection between the set of permutations on $n$ letters and standard permutation on $n+2$ letters. We emphasise that the irreducibility property is not preserved. For counting permutations in Rauzy classes, as we did in Section~\ref{section:counting_standard_permutations} we analyze the geometric surgery associated to this combinatorial operation.

\subsubsection{Marked profile, relation between $\gamma^{std}$ and $c$}
As in Lemma~\ref{lemma:profile_concatenation}, the profile of the degeneration depends only on the profile of the initial permutation and its marking. The proposition below expresses the profile of the degeneration from the profile of a standard permutation.
\begin{proposition} \label{prop:profile_of_degeneration}
Let $\pi$ be a standard permutation. If $\pi$ has a marked profile of the first type $(m|a,\,p')$, then its degeneration $\tilde{\pi}$ has marked profile $(m-2|m-a-2,\,p')$. If $\pi$ has a marked profile of the second type $(m_l \odot m_r,\,p')$, then its degeneration $\tilde{\pi}$ has a marked profile $(m_l-1) \odot (m_r-1),\,p')$.
\end{proposition}

\begin{proof}
We write the standard permutation $\pi$, in the following form
\[
\pi = \left( \begin{array}{c}
1\  y_1\, \ldots\, x_0\ 0 \\
0\ y_0\, \ldots\, x_1\ 1 \\
\end{array} \right).
\]
Let $\Gamma$ be the interval diagram of $\pi$. If the marking of $\pi$ is of the first type, let say $m|a$, then the corresponding singluarity in its interval diagram writes $c = (\incoms{x_0}, \outgos{(0,1)}, \incoms{x_1}, A, \outgos{y_0}, \incoms{(0,1)}, \outgos{y_1}, B)$ where $A$ and $B$ are some blocks of lengths respectively $2(m-a-2)$ and $2(a-1)$. Let $\tilde{\pi} = \binom{y_1 \ldots x_0}{y_0 \ldots x_1}$ be the degeneration of $\pi$ and $\tilde{\Gamma}$ its interval diagram. The interval diagram $\tilde{\Gamma}$ is obtained from the one of $\Gamma$ by modifying $c$ as $\tilde{c} = (\incoms{(x_0,x_1)}, A, \outgos{(y_0,y_1)}, B)$ where the blocks $A$ and $B$ have not changed. The angle between the left end point and the right endpoint is $|B|$. Hence, the permutation $\tilde{\pi}$ has a marking of the first type $m-2|m-a-2$.

Now, consider the case of a marking of the second type. By symmetry, it is enough to consider one endpoint of the interval. Let $c_l$ be the cycle of the interval diagram that contains the left end point. It writes $c_l = (\incoms{x_0}, \outgos{(0,1)}, \incoms{x_1}, A)$ and becomes $(\incoms{(x_0,x_1)},A)$ in the degeneration $\tilde{\pi}$ and proves the the proposition.
\end{proof}
From Proposition~\ref{prop:profile_of_degeneration}, we deduce a corollary about the relations between the numbers $c(p)$ of Section~\ref{section:counting_standard_permutations} and the numbers $\gamma^{irr}(p)$ and $\gamma(p)$. For an integer partition $p'$, we denote by $z_{p'}$ the cardinality of the centralizer of any permutation in the conjugacy class associated to $p'$. Let $e_i$ be the number of parts equal to $i$ in $p'$ then
\[
z_{p'} = \prod_{i=1}^n i^{e_i} e_i!.
\]
\begin{corollary}
\label{corollary:c_and_gamma}
Let $p=(m|a,\,p')$ be a marking of the first type then
\[
\gamma(m|a,\,p') = \gamma^{std}(m+2|m+2-a',\,p') 
\quad \text{and} \quad
\gamma^{std} \left(m|a,p'\right) = \frac{c \left(p_{m|a}\right)}{z_{p'}},
\]
Let $p=(m_l \odot m_r,\,p')$ be a marking of the second type then
\[
\gamma(m_l \odot m_r,\,p') = \gamma^{std}((m_l+1) \odot (m_r+1),\, p')
\quad \text{and} \quad
\gamma^{std} \left(m_l \odot m_r,p'\right) =  \frac{c \left(p_{m_l \odot m_r}\right)}{z_{p'}}.
\]
\end{corollary}

\begin{proof}
The two equalities on the left follows from Proposition~\ref{prop:profile_of_degeneration} as the degeneration is a bijection.

Recall that $c(p)$ counts the number of labeled standard permutations while $\gamma^{std}(p)$ counts unlabeled ones. Given a standard permutation $\pi$ the different ways we have to label it with a fixed labelization $\outgos{\sigma_\pi}$ is exactly $z_{p'}$.
\end{proof}

\subsubsection{Spin parity, relation between $\delta^{std}$ and $d$}
\label{subsubsection:spin_behavior}
In order to get a counting formula relative to the spin invariant, we now analyze the relation between the spin parity of a standard permutation $\pi$ and the one of its degeneration.
\begin{proposition}
\label{prop:spin_of_degeneration}
Let $\pi$ be a standard permutation of length $n+2$ and note $\alpha_1 = \pi_t^{-1}(1) = \pi_b^{-1}(n+2)$ and $\alpha_2 = \pi_b^{-1}(1) = \pi_t^{-1}(n+2)$. If $\pi$ has a marking of the first type $m|a$, then $\tilde{\pi}$ has a spin parity which is $\Arf(q_{\tilde{\pi}}) = \Arf(q_\pi)+a+1$ modulo $2$. If $\pi$ has marking of the second type, then the spin parity of $\tilde{\pi}$ is the same as the one of the permutation obtained by removing the letter $\alpha_1$ or $\alpha_2$ in $\pi$.
\end{proposition}

\begin{proof}
Let $\pi$ having a marking of the first type and $\tilde{\pi} = \pi_1 \cdot \pi_2 \cdot \ldots \cdot \pi_k$ the decomposition of $\tilde{\pi}$ into irreducible components. We denote by $S_\pi$ a suspension of $\pi$, $S_{\tilde{\pi}}$ a suspension of $\tilde{\pi}$ and $S_{\pi_j}$ the one induced on each irreducible components. Let $\zeta_j$ for $j=0,\ldots,n+1$ be the sides of the suspension $S_{\tilde{\pi}}$ (see Section~\ref{subsubsection:spin_invariant} and in particular Figure~\ref{fig:homology_of_suspension}). As the marking of $\pi$ is of type one, both intervals labeled $0$ and $n+1$ have the same singularities at both ends. Hence $\zeta_0$ and $\zeta_{n+1}$ are elements of $H_1(S;\ZZ/2)$ and there is a symplectic sum
\[
H_1(S; \ZZ/2) = (\ZZ/2\ [\zeta_0] \oplus \ZZ/2\ [\zeta_{n+1}]) \oplus \bigoplus_{j=1}^k H_1(S_j; \ZZ/2).
\]
The form $q_\pi$ diagonalizes with respect to this decomposition as its bilinear form is $\Omega_\pi$ which is the intersection form in $H_1(S;\ZZ/2)$. We hence only need to compute the restriction of $q_\pi$ to the symplectic module of rank two $M = (\ZZ/2\ [\zeta_0] \oplus \ZZ/2\ [\zeta_{n+1}])$. A direct computation shows that
\[
q_\pi(\zeta_0) = w(\zeta_0) + 1 + 0 = a +1 = q_\pi(\zeta_n) \quad \text{and} \quad q(\zeta_0 + \zeta_1) = 1
\]
Hence $Arf(q_\pi |_M) = a+1$. By additivity of the Arf invariant we get $\sum Arf(q_{\pi_j}) + a + 1 = Arf(q_\pi)$.

Now, we consider a permutation $\pi$ with marked profile $(m_l \odot m_r, p')$. If $\tilde{\pi}$ has a spin parity then both $m_l$ and $m_r$ are even. If we remove the interval labeled $0$ (or $n+1$) the permutation has profile $(m_l+m_r-1|a, p')$. The conservation of the spin statement is a direct consequence of Lemma~\ref{lemma:collapsing} of the preceding section.
\end{proof}
Let $\delta^{std}(p)$ be the difference between the number of odd spin permutations and even spin permutations among standard permutations with profile $p$.
\begin{corollary}
\label{corollary:d_and_delta}
Let $(m|a,p')$ be a marked profil of type one then
\[
\delta(m|a,\,p') = (-1)^{(a+1)} \delta^{std}(m+2|m+2-a',\,p')
\quad \text{and} \quad
\delta^{std} \left(m|a,p'\right) =
\left\{ \begin{array}{ll}
0 & \text{if $a \equiv 0 \mod 2$}, \\
\frac{d \left(p_{m|a}\right)}{z_{p'}} & \text{otherwise}.
\end{array} \right.
\]
Let $(m_l \odot m_r)$ be a marked profile of type two then
\[
\delta(m_l \odot m_r, p') = \frac{d \left( p_{(m_l+1) \odot (m_r+1)} \right)}{z_{p'}}
\quad \text{and} \quad
\delta^{std} \left(m_l \odot m_r,\,p'\right) =  \frac{d \left(p_{m_l \odot m_r}\right)}{z_{p'}}
\]
\end{corollary}

\begin{proof}
The proof is similar to the one of Corollary~\ref{corollary:c_and_gamma}. The left equalities follows from Proposition~\ref{prop:spin_of_degeneration} and the right ones from the definition of $d$.
\end{proof}

\subsection{Counting permutations in Rauzy classes}
\subsubsection{Marked points and hyperelliptic strata}
\label{subsubsection:marked_points_and_hyperelliptic_components}
As we did in Section~\ref{subsubsection:labeled_standard_permutations_marked_points} with labeled standard permutations, we give a relation between cardinalities of a Rauzy diagram and the ones obtained by adding marked points. As a corollary, we get the cardinality of any hyperelliptic Rauzy class.

Let $\mathbf{p}$ be a marked profile which corresponds to an hyperelliptic strata $\Omega \mathcal{M}(2g-1,1^k)$ or $\Omega \mathcal{M}(g,g,1^k)$. We denote by $hyp(\mathbf{p})$ the number of irreducible permutations with marked profile $p$. From the explicit description of the Rauzy class associated to rotation class permutation and hyperelliptic class (Section~\ref{subsubsection:symmetric_permutation}) we get the two following proposition.
\begin{proposition}
We have
\[
\gamma^{irr}(1|0,\,(1^{n-2})) = n-1 \quad \text{and} \quad \gamma^{irr}(1 \odot 1,\, (1^{n-3})) = \frac{(n-1)(n-2)}{2}.
\]

If $n$ is even the profile of $\pi$ is $p = (n-1)$ and the genus of a suspension of $\pi$ is $g=n/2$. In this case for $a \leq g-1$ we have
\[
hyp(2g-1|a) = hyp(2g-1|2g-1-a) = \binom{2g-1}{2a+1}
\]

If $n$ is odd, the profile of $\pi$ is $((n-1)/2,(n-1)/2)$ and the genus of a suspension is $g = (n-1)/2$. In this case for $a \leq g$ we have
\[
hyp(2g-1|a) = hyp(2g-1|2g-a) = \binom{2g-1}{2a+1}
\quad \text{and}
\quad hyp(g \odot g) = \sum_{k=0}^{g-1} \binom{2g}{2k} = 2^{2g-1} - 1
\]
\end{proposition}

Let $\mathcal{C} \subset \Omega \mathcal{M}(\kappa)$ be a connected component of a stratum and $\mathcal{R}$ its associated Rauzy diagram. We assume that the partition $\kappa$ does not contain $1$. Consider $\mathcal{C}_0 \subset \Omega \mathcal{M}(\kappa \uplus 0^k)$ the connected component obtained by marking $k$ points in the surfaces of $\mathcal{C}$. Let $\mathcal{R}_0$ be the extended Rauzy diagram associated to $\mathcal{R}_0$. The following theorem shows that the cardinality of $\mathcal{R}_0$ is a linear combination of the cardinality of $\mathcal{R}$ and the number of standard permutations in $\mathcal{R}$. Recall that $\mathcal{R}(m)$ denotes for $m-1$ an element of $\kappa$ the Rauzy class which correspond to the elements $\pi \in \mathcal{R}$ such that the left end point has an angle $2m\pi$ (see Section~\ref{subsubsection:connected_components_and_Rauzy_classes}).
\begin{theorem} \label{thm:cardinality_of_extended_rauzy_classes_with_marked_points}
Let $\mathcal{R}$, $\mathcal{R}_0$ and $k$ be as above. Le let $d$ be the number of letters in the permutations of $\mathcal{R}$, $r$ the number of standard permutations in $\mathcal{R}$ and $m$ an element of the profile of $\mathcal{R}$. Then
\[
|\mathcal{R}_0(m)| = \binom{d+k}{k} |\mathcal{R}(m)|
\quad \text{and} \quad 
|\mathcal{R}_0(1)| = \binom{d+k}{k-1} |\mathcal{R}| + \binom{d+k}{k-1}\ d\,r
\]
In particular, for the cardinalities of extended Rauzy classes, we get the following relations
\[
|\mathcal{R}_0| = \binom{d+k+1}{k}\ |\mathcal{R}| + \binom{d+k}{k-1}\ d\,r.
\]
\end{theorem}
The proof of the theorem follows from Proposition~\ref{prop:adding_marked_points} below. As a corollary of the theorem, we get an explicit formula for the cardinality of Rauzy diagrams associated to any hyperelliptic component of stratum.
\begin{corollary}
\label{cor:cardinality_of_hyperelliptic_rauzy_diagrams}
Let $\mathcal{R}$ be the extended Rauzy diagram of the hyperelliptic component $\Omega \mathcal{M}_{hyp}(2g-1,1^k)$ (reps. $\Omega \mathcal{M}(g,g,1^k)$) for which $d=2g$ (resp. $d=2g+1$) is the number of intervals in $\Omega \mathcal{M}(2g-1)$ (resp. $\Omega \mathcal{M}(g,g)$). Then the cardinality of the Rauzy diagrams are
\[
|\mathcal{R}(1)| = \binom{d+k}{k-1} (2^{d-1}+d-1)
\quad \text{and} \quad
|\mathcal{R}(d-1)| = \binom{d+k}{k} (2^{d-1}-1).
\]
The cardinality of the extended Rauzy diagram is
\[
|\mathcal{R}| = \binom{d+k+1}{k}(2^{d-1}-1) + \binom{d+k}{k-1} d.
\]
\end{corollary}

We now prove Theorem~\ref{thm:cardinality_of_extended_rauzy_classes_with_marked_points}. As above, let $\mathcal{R}$ be an extended Rauzy class and $\mathcal{R}_0$ the one obtained by adding $k$ marked points. We denote by $p$ the profile of the permutations in $\mathcal{R}$ and we assume that $1 \not\in p$. If $m|a$ (resp. $m_l \odot m_r$) is a marking of the first type (resp. second type) then we denote by $\mathcal{R}(m|a)$ (resp. $\mathcal{R}(m_l \odot m_r)$) the elements of the extended Rauzy class $\mathcal{R}$ which has marking $m|a$ (resp. $m_l \odot m_r$).
\begin{proposition}
\label{prop:adding_marked_points}
Let $\mathcal{R}$ and $\mathcal{R}_0$ be extended Rauzy classes as above, then
\begin{enumerate}
\item $|\mathcal{R}_0(m|a)| = \binom{d+k-1}{k} |\mathcal{R}(m|a)|$,
\item $|\mathcal{R}_0(m_l \odot m_r)| = \binom{d+k-1}{k} |\mathcal{R}(m_l \odot m_r)|$,
\item $|\mathcal{R}_0(m \odot 1)| = |\mathcal{R}_0(1 \odot m)| = \binom{d+k-1}{k-1} |\mathcal{R}(m)|$,
\item $|\mathcal{R}_0(1 \odot 1)| = \binom{d+k-1}{k-2} (|\mathcal{R}| + d\,r)$,
\item $|\mathcal{R}_0(1|0)| = \binom{d+k-1}{k-1}\ d\,r$.
\end{enumerate}
\end{proposition}
\begin{proof}
We first prove equalities \textit{1} and \textit{2}. Let $\pi \in \mathcal{R}$ with marking $m|a$ or $m_l \odot m_r$ and $P_0 \subset \mathcal{R}_0$ the set of permutations $\pi_0$ with the same marking as $\pi$ and such that they are obtained from $\pi$ by adding $k$ zeroes. The marked points of any $\pi_0 \in P_0$ belong inside the $d$ intervals. Hence $|P_0| = \binom{d+k-1}{k}$ is the number of choices of placing $k$ undifferentiated points in $d$ intervals.

Now, we proove equality \textit{3}. Let $\pi \in \mathcal{R}(m)$ and $P_0 \subset \mathcal{R}_0(m \odot 1)$ the set of permutations obtained from $\pi$ by adding $k$ marked points. For any $\pi_0 \in P_0$, because of the marking $m \odot 1$ and $1 \not\in p$ one of the marked point has to go to the right endpoint of the permutation. There is only one way to do this by the following operation
\[
\pi = \binom{\ldots\, y\, \ldots\, x}{\ldots\, x\, \ldots\, y} \mapsto \pi_0 = \binom{\ldots\, y\ c\, \ldots\, x}{\ldots\, x\ c\, \ldots\, y}.
\]
Then, the $k-1$ other marked points belong in the $d+1$ intervals and the number of choices for such operation is $\binom{(d+1)+(k-1)-1}{k-1} = \binom{d+k-1}{k-1}$. Hence $|P_0| = \binom{d+k-1}{k-1}$.

Equality \textit{4} is similar to equality \textit{3} but two of the marked points have to be placed at the extremities. We get the coefficient $\binom{(d+2)+(k-2)-1}{k-2} = \binom{d+k-1}{k-2}$.

We now proove equality \textit{5}. Let $\pi \in \mathcal{R}_0$ be a permutation with marking $1|0$. Then we can write a general form for $\pi_0$ and we see below that removing the marked point of $\pi$ gives a standard permutation.
\begin{equation}
\pi_0 = \left( \begin{array}{cccccc}
a_0 & A & b_1 & a_1 & B & b_0 \\
a_1 & C & b_0 & a_0 & D & b_1
\end{array} \right)
\quad \mapsto \quad
\left( \begin{array}{ccccc}
a & A & c & B & b \\
b & C & c & D & a
\end{array} \right).
\label{eq:remove_mark_point_std}
\end{equation}
Hence, the only way to mark $1$ point on a permutation in $\mathcal{R}$ in order to obtain a marking $1|0$ is that $\pi$ is standard. Starting from a standard permutation $\pi \in \mathcal{R}$ the construction of a permutation $\pi_0$ with marking $1|0$ is as follows. Choose the letter $c$ which will play the role of an intermediate and place it as in (\ref{eq:remove_mark_point_std}). There are $d$ choices for the letter $c$. Then, the other $k-1$ points can be placed inside $d+1$ intervals. We get exactly $\binom{d+k-1}{k-1}\ d$ permutations in $\mathcal{R}_0$ built from $\pi$.
\end{proof}

\subsubsection{The number of irreducible permutations}
\label{subsubsection:the_number_of_irreducible_permutations}
Before counting permutations in Rauzy diagrams, we recall the elementary method to count irreducible permutations. Most of the idea developed here are similar to the one we will use in the next section. As in (\ref{eq:irreducible}), let $f(n) = |S_n^o|$ be the number of irreducible permutations of length $n$. We recall the elementary method for different expression of (\ref{eq:irreducible}) and get an asymptotic development. See the original article \cite{Comtet72} for further details on asymptotics and \cite{FlajoletSedgewick09} for general considerations about the relations between generating series and operations on combinatorial classes.

Let $E(t) = \sum n!\ t^n$ and $F(t) = \sum f(n)\ t^n$ considered as formal serie. Given a permutation, its factorization in irreducible elements is unique. In terms of the generating functions $E$ and $F$, the equation (\ref{eq:irreducible}) can be seen as
\begin{equation*}
E = \frac{1}{1-F} = \sum_{k = 0}^\infty F^k.
\end{equation*}
Using an inclusion/exclusion argument, we get a dual formulation of the equation (\ref{eq:irreducible})
\begin{equation}
f(n) = \sum_{k=1}^n \sum_{c_1 + \ldots + c_k = n} (-1)^{k+1}\, c_1!\, \ldots\, c_k! 
\qquad \text{or} \qquad
F = 1 - \frac{1}{E} = 1 - \sum_{k = 0}^\infty (1-E)^k.
\label{eq:irreducible_dual}
\end{equation}
We can write a simpler relation between factorial numbers $n!$ and the numbers $f(n)$. Any permutation can be decomposed uniquely as the product of an irreducible permutation and a permutation. Hence
\begin{equation}
\sum_{i=1}^{n}\, f(i)\, (n-i)! = n!
\qquad \text{or} \qquad
E F = E - 1.
\label{eq:irreducible_mixed}
\end{equation}
From the equations on generating functions, we see that the formulas (\ref{eq:irreducible}), (\ref{eq:irreducible_dual}), (\ref{eq:irreducible_mixed}) are equivalent. However each one  has its own advantage: equation (\ref{eq:irreducible}) is the most natural, equation (\ref{eq:irreducible_dual}) gives a closed formula and (\ref{eq:irreducible_mixed}) is adapted for explicit computations.

The equation (\ref{eq:irreducible_mixed}) suffices to obtain an equivalent of the number of irreducible permutations. For an asymptotic serie, see \cite{Comtet72}.
\begin{proposition}[\cite{Comtet72}]
\label{prop:number_of_irreducible_permutations}
$f(n)$ is equivalent to $n!$ (e.g. $f(n) = n! (1 + o(1))$).
\end{proposition}

\begin{proof}
Let $g(n) := f(n)/n!$. Those numbers satisfy the inequality $g(n) \leq 1$ and from Equation~(\ref{eq:irreducible_mixed}) we get
\begin{align*}
g(n) &= 1 - \sum_{k=1}^{n-1} g(k) \binom{n}{k}^{-1} \\
   &= 1 - \frac{2}{n} - \sum_{k=2}^{n-2} \binom{n}{k}^{-1} \\
   &\geq 1 - \frac{2}{n} - (n-3)\ \frac{2}{n(n-1)}.
\end{align*}
As the right member of this equation tends to $1$ we get that $g(n)$ tends to $1$ as $n$ tends to $\infty$.
\end{proof}

\subsubsection{Formula for $\gamma^{irr}$ and $\delta^{irr}$, proof of Theorem~\ref{thm0:cardinality}}
\label{subsubsection:formula_for_gamma_irr_and_delta_irr}
Recall that $\gamma(m|a,\,p')$ and $\gamma(m_l \odot m_r,\,p')$ (resp. $\gamma^{irr}(m|a,\,p')$ and $\gamma^{irr}(m_1 \odot m_2,\,p')$) denote the number of permutations (resp. irreducible permutations) with marked profile $(m|a,\,p')$ and $(m_1 \odot m_2,\,p')$. The numbers $\gamma(m|a,\,p')$ and $\gamma(m_l \odot m_r,\,p')$ are related to the number $c(p)$ of Section~\ref{section:counting_standard_permutations} by Corollary~\ref{corollary:c_and_gamma}.

The two formulas in Theorem~\ref{thm0:cardinality} are obtained by an exclusion procedure and are very similar to (\ref{eq:irreducible_mixed}) which gives an explicit formula for the number $f(n)$ of irreducible permutations in $S_n$ in the following form
\[
f(n) = n! - \sum_{k=1}^{n-1} f(k) (n-k)!.
\]
In the above formula, $n!$ corresponds to the cardinality of permutations and the summation corresponds to all reducible ones. Each reducible permutation has to be thought as the concatenation of an irreducible permutation of length $k$ with any permutation of length $n-k$.

We explain the formula for $\gamma^{irr} (m_l \odot m_r, p')$, the other being similar. The set of all permutations (non necessarily reducible) with marked profile $m_l \odot m_r, p')$ is exactly $\gamma^{std}((m_l+1) \odot (m_r+1),\, p')$ (see Section~\ref{subsection:degeneration} and in particular Proposition~\ref{prop:profile_of_degeneration}). Then we have to subtract all irreducible. Recall from Lemma~\ref{lemma:profile_concatenation} that the profile of a reducible permutation can be expressed in terms of its irreducible components. We consider the possible form of a reducible permutation $\pi_1 \pi_2$ with marked profile $(m_l \odot m_r, p')$ where $\pi_1$ is irreducible.
\begin{enumerate}
\item either $\pi_1$ has a marking of the first type and $\pi_2$ a marking of the second type,
\item or $\pi_1$ has marking of the second type and $\pi_2$ of the first type,
\item or $\pi_1$ and $\pi_2$ both have marking of the second types.
\end{enumerate}
The three cases above, correspond to the three summations in the formula $\gamma^{irr} (m_l \odot m_r, p')$ in Theorem~\ref{thm0:cardinality}.

\subsubsection{Explicit formula for profile $(2g-1)$}
We gave in Section~\ref{subsubsection:marked_points_and_hyperelliptic_components} examples of family of Rauzy classes obtained by adding marked points. Theorem~\ref{thm:cardinality_of_extended_rauzy_classes_with_marked_points} gives an explicit formula for the behavior of the cardinality. In those example, the genus was fixed. In this section we consider the family of Rauzy diagrams which are the Rauzy classes associated to the odd and even components of $\Omega \mathcal{M}(2g-1)$. This family of strata are the so called \emph{minimal strata}. Recall that for $g=2$, $\Omega \mathcal{M}(2g-1)$ has only one connected components, for $g=3$ there are $2$ and for $g \geq 4$ there are three. The cardinality of the hyperelliptic component is given in Proposition~\ref{prop:rauzy_class_of_symmetric_permutation}. To get the cardinality of all Rauzy classes, we consider explicit formulas for the numbers $\gamma^{irr}(2g-1)$ and $\delta^{irr}(2g-1)$ in the following proposition.
\begin{proposition}
Let $n=2g-1$ then we have the following formulas for $\gamma^{irr}((n))$ and $\delta^{irr}((n))$
\begin{equation}
\gamma^{std}((n)) = \frac{(n-1)!}{n+1}
\quad \text{and} \quad 
\delta^{std}((n)) = \frac{(n-1)!}{2^{n-1}}\ ,
\end{equation}

\begin{equation}
\gamma^{irr}((n)) = 
\sum_{k=1}^{2n+1} (-1)^{k+1}
\sum_{(c_1,\ldots,c_k) \atop \sum c_j=n+1} \quad \prod_{i=1}^k \frac{(2\ c_i)!}{c_i+1}\ ,
\end{equation}

\begin{equation}
\delta^{irr}((n)) = 
\frac{1}{2^{n+1}} 
\sum_{k=1}^{2n+1} (-1)^{k+1} 
\sum_{(c_1,\ldots,c_k) \atop \sum c_j=n+1} \quad \prod_{i=1}^k (2 c_i)! \ .
\end{equation}
\end{proposition}

\begin{proof}
This is a direct consequence of Corollary~\ref{corollary:c_explicit_formulas} and the explicit formula for $d$ in \ref{thm:formula_for_d}.
\end{proof}

\nocite{*}
\bibliographystyle{alpha}
\bibliography{biblio}

\end{document}